\documentclass[11pt,reqno]{amsart}
\usepackage{amsmath,amssymb,mathrsfs,amscd}
\usepackage{diagbox}
\usepackage{graphicx}
\usepackage{subfigure}
\usepackage{cite}
\usepackage{epstopdf}
\usepackage{color}
\usepackage[colorlinks,linkcolor=red,anchorcolor=blue,citecolor=blue]{hyperref}
\usepackage[mathscr]{euscript}
\usepackage [section]{placeins}
\usepackage{booktabs}
\usepackage{multirow}
\usepackage[ruled,linesnumbered]{algorithm2e}

\newtheorem{thm}{Theorem}[section]

\newtheorem{lem}[thm]{Lemma}
\newtheorem{prop}[thm]{Proposition}
\theoremstyle{definition}
\newtheorem{defn}[thm]{Definition}
\theoremstyle{remark}
\newtheorem{rem}[thm]{Remark}
\numberwithin{equation}{section}

\allowdisplaybreaks[4]

\begin{document}
\title[]{An alternating approach for reconstructing the initial value and source term in a time-fractional diffusion-wave equation}%
\author{Yun Zhang$^{1,*}$, Xiaoli Feng$^{1}$, Xiongbin Yan$^{2}$}%
\email{zhangyun@xidian.edu.cn (Y. Zhang)}%
\email{xiaolifeng@xidian.edu.cn (X. Feng)}%
\email{yanxb2015@163.com (X. Yan)}%
\address{$^*$: Corresponding author}%
\address{1: School of Mathematics and Statistics,~Xidian University, Xi'an~710126, ~PR China}
\address{2: School of Mathematical Sciences,~Shanghai Jiao Tong University, Shanghai~200240, ~PR China}%
%\date{}%
%\dedicatory{}%
%\commby{}%
% ----------------------------------------------------------------
\begin{abstract}
This paper is dedicated to addressing the simultaneous inversion problem involving the initial value and space-dependent source term in a time-fractional diffusion-wave equation. Firstly, we establish the uniqueness of the inverse problem by leveraging the asymptotic expansion of Mittag-Leffler functions. Subsequently, we decompose the inverse problem into two subproblems and introduce an alternating iteration reconstruction method, complemented by a regularization strategy. Additionally, a comprehensive convergence analysis for this method is provided. To solve the inverse problem numerically, we introduce two semidiscrete schemes based on standard Galerkin method and lumped mass method, respectively. Furthermore, we establish error estimates that are associated with the noise level, iteration step, regularization parameter, and spatial discretization parameter. Finally, we present several numerical experiments in both one-dimensional and two-dimensional cases to validate the theoretical results and demonstrate the effectiveness of our proposed method.
\end{abstract} \maketitle


\section{Introduction}\label{sec1}

Diffusion and wave propagation are two prominent physical phenomena in the natural world. The governing equation can be characterized as follows:
\begin{equation}\label{diffwave}
\partial_t^{\alpha}u(x,t)-\Delta u(x,t)=F(x,t),
\end{equation}
where $\alpha=1,2$. When $\alpha=1$, (\ref{diffwave}) represents the standard diffusion equation, which is crucial in describing phenomena such as heat conduction, mass transfer, and particle diffusion. Meanwhile, when $\alpha=2$, (\ref{diffwave}) corresponds to the wave equation, used to depict the propagation and characteristics of waves like sound waves, light waves, and water waves. By restricting $\alpha\in(1,2)$, (\ref{diffwave}) can be interpreted as the time-fractional diffusion-wave equation (TFDWE). TFDWEs capture the physical phenomena involving intermediate diffusion processes and wave propagation. Due to the nonlocality and memory property of fractional derivatives, TFDWEs excel in decribing intricate physical mechanisms such as diffusion processes in viscoelastic media. Furthermore, they are highly significant in various applications spanning physics, chemistry, biology, engineering, and financial markets.

Based on above motivations, we will study the initial-boundary value problem of the following TFDWE in this paper:
\begin{align}\label{P}
\left\{
\begin{aligned}
&\partial_{t}^{\alpha}u(t)-\Delta u(t)=F(t)& &in\quad\Omega\times(0,T),\\
&u(t)=0& &on\quad\partial\Omega\times[0,T],\\
&u(0)=a,~\partial_tu(0)=b& &in\quad\overline{\Omega},
\end{aligned}
\right.
\end{align}
where $\Omega\subset\mathbb{R}^{d}(d=1,2,3)$ is a bounded domain with sufficiently smooth boundary. For fixed $\alpha\in(1,2)$, $\partial_{t}^{\alpha}u(t)$ denotes the left-sided Caputo fractional derivative of order $\alpha$ defined by
\begin{equation*}
\partial_{t}^{\alpha}u(t):=\frac{1}{\Gamma(2-\alpha)}\int_{0}^{t}(t-\tau)^{1-\alpha}\partial_{\tau}^2 u(\tau)d\tau,
\end{equation*}
where $\Gamma(\cdot)$ is the Gamma function.

The direct problem involves finding the solution $u$ in problem (\ref{P}) given the initial values $a$, $b$ and the source term $F(t)$. However, in many real applications, it is difficult or impossible to obtain these parameters directly.  In such instances, we must deduce or infer them from additional measurements, giving rise to what is commonly referred to as the inverse problem.

Over the last two decades, there has been significant progress in the research of inverse problems for time-fractional diffusion equation (TFDE). These problems can be classified based on the parameters involved, including inverse initial value problem (IIVP) \cite{Liu+Yamamoto-2010,Wang+Wei+Zhou-2013,Floridia+Yamamoto-2020}, inverse source problem (ISP) \cite{Zhang+Xu-2011,Wei+Li+Li-2016,Jiang+Li+Liu+Yamamoto-2017,Janno+Kinash-2018}, inverse coefficient problems \cite{Cheng+Nakagawa+Yamamoto+Yamazaki-2009,Li+Zhang+Jia+Yamamoto-2013,Fujishiro+Kian-2016}, and inverse fractional order problems \cite{Li+Yamamoto-2015,Li+Huang+Yamamoto-2020,Yamamoto-2023}. Additionally, there are other types of inverse problems such as inverse random source for TFDE \cite{Feng+Li+Wang-2020,Niu+Helin+Zhang-2020,Lassas+Li+Zhang-2023}, inverse problems for time-space fractional diffusion equations \cite{Jia+Peng+Yang-2017,Helin+Lassas+Ylinen+Zhang-2020} and nonlinear TFDE \cite{Janno+Kasemets-2017,Kaltenbacher+Rundell-2019-2,Ma+Sun-2023}.

Inverse problems are frequently considered ill-posed, which means that the existence, uniqueness, and stability of the solutions may not hold simultaneously. Consequently, various regularization techniques must be employed to tackle the ill-posed nature of these problems. Numerous efficient regularization methods have been developed over the past years for linear and nonlinear inverse problems. For linear inverse problems, Liu and Yamamoto \cite{Liu+Yamamoto-2010} proposed a quasi-reversibility method for the backward problem of TFDE in one dimensional case. Wang et al. \cite{Wang+Wei+Zhou-2013} studied the backward problem for TFDE in general domain by Tikhonov approach. Wei et al. \cite{Wei+Li+Li-2016}  investigated the inverse time-dependent source problem. They provided the uniqueness as well as the stability estimate and reconstructed the source term by conjugate gradient method. Hào et al. \cite{Hao+Liu+Duc+Thang-2019} studied the backward problem by a non-local boundary value method. For nonlinear inverse problems, Zhang et al. \cite{Zhang+Jia+Yan-2018} studied the problem of recovering the fractional order and source term by Bayesian approach.  Kaltenbacher and Rundell \cite{Kaltenbacher+Rundell-2019} reconstructed the potential term by iterative Newton-type methods. Jin and Zhou \cite{Jin+Zhou-2021} investigated the problem of recovering the diffusion coefficient by a regularized output least-squares formulation and provided a fully discrete scheme for numerically solving the inverse problem. Recently, people became concerning on the numerical analysis for the inverse problems. Jiang et al. \cite{Jiang+Liu+Wang-2020} recovered the source term by Tikhonov regularization method and provided a fully discrete scheme for computing the numerical approximated solution and established the convergence analysis. Fan and Xu \cite{Fan+Xu-2021} transformed the ISP in subdiffusion equation into a regularized problem with $L^2$-$TV$ penalty terms. And a discrete scheme based on finite element method is proposed for obtaining the approximate solution. Zhang and Zhou \cite{Zhang+Zhou-2020} studied the backward problem for subdiffusion equation by quasi-boundary value method. They also provided semidiscrete scheme and fully discrete scheme for solving the regularized problem. 

In recent years, there has been growing interest in the simultaneous inversion of multiple parameters in a physical system. This paper focuses on the problem of recovering of initial value and source term simultaneously in problem (\ref{P}). The study of both IIVP and ISP has extensive applications in various fields such as seismology, sound wave propagation, electromagnetic fields, materials science, biomedical engineering and financial engineering. The IIVP for TFDE aims to detect the initial state of an anomalous diffusion process from additional boundary conditions or final time data. On the other hand, the ISP for TFDE refers to inferring the source term in the anomalous diffusion system from observation data. These studies contribute to our understanding of scenarios where the initial state of the system is not directly observable, and assist in identifying the location and strength of the source term through observations of the system’s boundary or interior responses.

In summary, the research on inverse problems for fractional differential equations has seen significant developments in various areas, and the application of regularization methods has played a crucial role in addressing the ill-posed nature of these problems. The simultaneous inversion of multiple parameters in physical systems has also become an area of active research, with potential applications across diverse fields.

When $\alpha\in(0,1)$, problem (\ref{P}) is referred to the TFDE as well as subdiffusion equation. In \cite{Ruan+Yang+Lu-2015}, the authors investigated the uniqueness of the problem of simultaneous inversion of initial value and space-dependent source term, and solved the inverse problem by Tikhonov approach. Wen et al. studied the problem of simultaneous recovering initial value and source term by non-stationary iterative Tikhonov regularization method \cite{Wen+Liu+Wang-2023} and Landweber iteration method \cite{Wen+Liu+Yue+Wang-2022}, respectively. Yu et al. \cite{Yu+Wang+Yang-2023} proposed an exponential Tikhonov regularization method for the simultaneous inversion problem. To the best of the authors' knowledge, there are rare studies on the problem of simultaneously inversion of the initial value and source term when $\alpha\in(1,2)$. 

In this paper, we always assume that $b\equiv0$ and  the source term is time-independent, i.e., $F(t)\equiv f$. Then problem (\ref{P}) becomes
\begin{align}\label{IP}
\left\{
\begin{aligned}
&\partial_{t}^{\alpha}u(t)-\Delta u(t)=f& &in\quad\Omega\times(0,T),\\
&u(t)=0& &on\quad\partial\Omega\times[0,T],\\
&u(0)=a,~\partial_tu(0)=0& &in\quad\overline{\Omega},\\
&u(T_1)=g_1,~u(T_2)=g_2& &in\quad\overline{\Omega}.
\end{aligned}
\right.
\end{align}
The inverse problem considered in this paper is to reconstruct the initial value and source term $(a,f)$ from two noisy terminal measurements $g_1^{\delta}$ and $g_2^{\delta}$ satisfying
\begin{equation}\label{noise}
\left\|g_1^{\delta}-g_1\right\|\leq\delta,\quad \left\|g_2^{\delta}-g_2\right\|\leq\delta,
\end{equation}
where $\|\cdot\|$ denotes the $L^2$-norm, $\delta$ is the noise level and $g_1=u(T_1)$, $g_2=u(T_2)$. Without loss of generality, we always assume that $T_2>T_1>0$. The main contributions of this paper are in the following three aspects:\\
(i) \textbf{The uniqueness of the inverse problem is given.} For $\alpha\in(1,2)$, the inverse problem is fairly complicated compare to the problem with $\alpha\in(0,1)$. In \cite{Wei+Zhang-2018}, Wei and Zhang studied the problem of recovering the initial value from final data. They demonstrated that the uniqueness of the problem may not hold if the terminal time is arbitrarily selected, which is quite different from TFDE and classical diffusion equation. In this study, we prove that the uniqueness holds for the problem of recovering the initial value and source term in a TFDWE, provided that the final time is chosen properly.\\
(ii) \textbf{An efficient alternating regularization method is proposed.} The inverse problem considered in this paper is both ill-posed and coupled. To efficiently solve the problem, we propose an alternating regularization method. First, we decouple the inverse problem into two subproblems of recovering the initial value and source term, respectively. Second,  we apply the quasi-boundary value regularization method for solving each inverse problem. By formulating this process as an iterative scheme, we can efficiently reconstruct the initial value and source term. Lastly, the convergence analysis is given in detail. The proposed method offers two primary advantages. On one hand, the initial value and source term can be recovered through parallel processing, making it suitable for distributed computing environments. On the other hand, the iterative approach typically demands less storage space. Furthermore, this method can also be used for solving a variety of inverse problems for partial differential equations of recovering multiple unknown parameters.\\
(iii) \textbf{Semidiscrete schemes based on standard Galerkin finite element method and lumped mass method are provided for numerically implementing the alternating regularization method.} When deducing the convergence result of the proposed method, the solutions are given in the form of infinite series by means of eigensystem of Laplacian. Thus we cannot numerically solve the inverse problems directly. By applying standard Galerkin method and lumped mass method for discretization the problem in space, we can obtain approximate solutions consisting of finite basis functions and Mittag-Leffler functions. The error estimate between the approximate solution and the exact solution related to the discretization parameter is also derived. Moreover, Garrappa \cite{Garrappa-2015} provided an efficient method for the numerical evaluation of Mittag-Leffler functions. Hence we can directly utilize the code of the corresponding method provided in \cite{Garrappa-2024} for computing the numerical solution. According to above procedure, we can ultimately obtain the numerical solutions of the unknown parameters.

The rest of this paper can be outlined organized as follows. In Section \ref{sec2}, we present some preliminary knowledge and regularity results for the solution of problem (\ref{P}) that are necessary for subsequent analysis. Following this, we introduce the alternating regularization approach and present comprehensive convergence results in Section \ref{sec3}. In Section \ref{sec4}, we introduce the semidiscrete schemes and provide error estimates between the exact solutions and the numerically solutions. Section \ref{sec5} showcases numerical examples in both one-dimensional and two-dimensional cases. Finally, we conclude this paper in Section \ref{sec6}.

\section{Preliminaries and Auxiliary results}\label{sec2}

\begin{defn}\cite{Kilbas-2006,I.Podlubny}
The generalized Mittag-Leffler function is defined by
\begin{equation*}
E_{\alpha,\beta}(z)=\sum_{k=0}^{\infty} \frac{z^k}{\Gamma(\alpha k+\beta)},\quad z\in \mathbb{C},
\end{equation*}
where $\alpha>0$, $\beta\in\mathbb{R}$.
\end{defn}

\begin{lem}\cite{I.Podlubny}\label{lem2.1}
Let $0<\alpha<2$, and $\beta\in\mathbb{R}$ be arbitrary. We suppose that $\mu$ satisfies $\frac{\pi\alpha}{2}<\mu<\min\{\pi,\pi\alpha\}$. Then there exists a constant $C=C(\alpha,\beta,\mu)>0$ such that
\begin{equation*}
|E_{\alpha,\beta}(z)|\leq\frac{C}{1+|z|},\quad\mu\leq|\arg(z)|\leq\pi.
\end{equation*}
\end{lem}

\begin{lem}\cite{I.Podlubny}\label{lem2.3}
For $\lambda>0$ and $1<\alpha<2$, we have 
\begin{equation}\label{lem2.3_eq1}
\partial_t^{\alpha}(E_{\alpha,1}(-\lambda t^{\alpha}))=-\lambda E_{\alpha,1}(-\lambda t^{\alpha}),\quad \partial_t^{\alpha}(tE_{\alpha,2}(-\lambda t^{\alpha}))=-\lambda tE_{\alpha,2}(-\lambda t^{\alpha}),\quad t>0.
\end{equation}
\end{lem}

\begin{lem}\cite{I.Podlubny}\label{lem2.2}
If $1<\alpha<2$, $\beta$ is an arbitrary complex number and $\mu$ is an arbitrary real number such that $\mu\in(\frac{\pi\alpha}{2},\pi)$, then for an arbitrary integer $p\geq1$, the following expansion holds:
\begin{equation}\label{lem22}
E_{\alpha,\beta}(z)=-\sum_{k=1}^{p}\frac{z^{-k}}{\Gamma(\beta-\alpha k)}+O(|z|^{-1-p}),\quad |z|\rightarrow\infty,\quad\mu\leq|\arg(z)|\leq\pi.
\end{equation}
\end{lem}

Let $\{\lambda_n,\varphi_n\}_{n=1}^{\infty}$ be the eigensystem of negative Laplacian defined on the region $\Omega$ with homogeneous Dirichlet boundary condition. Then  $\{\varphi_n\}_{n=1}^{\infty}$ forms an orthonormal basis in $L^2(\Omega)$. Next we define the Hilbert space $D((-\Delta)^p)$ for $p>0$ by
\begin{equation}\label{dlp}
D((-\Delta)^p):=\left\{f\in L^2(\Omega)\mid\sum_{n=1}^{\infty}\lambda_n^{2p}(f,\varphi_n)^2<\infty\right\}.
\end{equation} 
The corresponding norm on $D((-\Delta)^p)$ is defined by $\|f\|_{D((-\Delta)^p)}=\left(\sum_{n=1}^{\infty}\lambda_n^{2p}(f,\varphi_n)^2\right)^{\frac12}$. For general source term, the solution to (\ref{P}) can be expressed by 
\begin{equation}\label{solution1}
u(t)=\mathbb{S}(t)a+\mathbb{T}(t)b+\int_0^t\mathbb{G}(t-\tau)F(\tau)d\tau,
\end{equation}
where the solution operators $\mathbb{S}(t),\mathbb{T}(t)$ and $\mathbb{G}(t)$ are defined by
\begin{eqnarray*}
&&\mathbb{S}(t)v=\sum_{n=1}^{\infty}S_n(t)(v,\varphi_n)\varphi_n:=\sum_{n=1}^{\infty}E_{\alpha,1}(-\lambda_n t^{\alpha})(v,\varphi_n)\varphi_n,\\
&&\mathbb{T}(t)v=\sum_{n=1}^{\infty}T_n(t)(v,\varphi_n)\varphi_n:=\sum_{n=1}^{\infty}tE_{\alpha,2}(-\lambda_n t^{\alpha})(v,\varphi_n)\varphi_n,\\
&&\mathbb{G}(t)v=\sum_{n=1}^{\infty}G_n(t)(v,\varphi_n)\varphi_n:=\sum_{n=1}^{\infty}t^{\alpha-1}E_{\alpha,\alpha}(-\lambda_n t^{\alpha})(v,\varphi_n)\varphi_n.
\end{eqnarray*}
Thus we have the following results for the solution $u(t)$ defined in $(\ref{solution1})$.
\begin{prop}\cite{Zhang+Zhou-2022}\label{direct1}
Let $u(t)$ be defined in (\ref{solution1}). Then the following statements hold. \\
(i) If $a,~b\in D((-\Delta)^p)$ for $p\in[0,1]$ and $F\equiv0$, then $u(t)$ is the solution to problem (\ref{P}), and $u(t)$ satisfies the following estimate
\begin{equation}\label{solu_regab}
\|\partial_t^{(m)}u(t)\|_{D((-\Delta)^q)}\leq Ct^{-m-\alpha(q-p)}\left(\|a\|_{D((-\Delta)^p)}+t\|b\|_{D((-\Delta)^p)}\right),
\end{equation}
for any integer $m\geq0$ and $q\in[p,p+1]$.\\
(ii) If $a=b=0$ and $F\in L^p(0,T;L^2(\Omega))$ with $\frac{1}{\alpha}<p<\infty$, then $u(t)$ is the solution to problem (\ref{P}) such that $u\in C([0,T];L^2(\Omega))$ and
\begin{equation}\label{solu_regF}
\|u\|_{L^p(0,T;D(-\Delta))}+\|\partial_t^{\alpha}u\|_{L^p(0,T;L^2(\Omega))}\leq C\|F\|_{L^p(0,T;L^2(\Omega))}.
\end{equation}
\end{prop}

In the case that the source term is time-independent, the solution $u(t)$ to (\ref{IP}) defined in (\ref{solution1}) can be reduced to the following form: 
\begin{equation}\label{solution2}
u(t)=\mathbb{S}(t)a+\mathbb{Q}(t)f,
\end{equation}
where the solution operator $\mathbb{Q}(t)$ is defined by 
\begin{equation*}
\mathbb{Q}(t)v=\sum_{n=1}^{\infty}Q_n(v,\varphi_n)\varphi_n:=\sum_{n=1}^{\infty}t^{\alpha}E_{\alpha,\alpha+1}(-\lambda_nt^{\alpha})(v,\varphi_n)\varphi_n.
\end{equation*}
Accordingly, we have the following result to the solution $u(t)$ defined in $(\ref{solution2})$. The proof is similar to that of  Proposition \ref{direct1}, so we omit it here.

\begin{prop}\label{direct2}
Let $u(t)$ be defined in (\ref{solution2}). If $a,~f\in D((-\Delta)^p)$ for $p\in[0,1]$, then $u(t)$ is the solution to problem (\ref{IP}), and satisfies the following estimate
\begin{equation}\label{solu_regaf}
\|\partial_t^{(m)}u(t)\|_{D((-\Delta)^q)}\leq Ct^{-m-\alpha(q-p)}\left(\|a\|_{D((-\Delta)^p)}+t^{\alpha}\|f\|_{D((-\Delta)^p)}\right)
\end{equation}
for any integer $m\geq0$ and $q\in[p,p+1]$. 
\end{prop}
For subsequent analysis, we also need the following results.
\begin{lem}\cite{Wei+Zhang-2018,Wei+Luo-2022}\label{SQ_bound}
Let $1<\alpha<2$, then for sufficiently large $T_1>0$, there exist positive constants $\overline{C_1},\underline{C_1},\overline{C_2},\underline{C_2}$ such that 
\begin{equation}\label{SQ_S}
-\frac{\overline{C_i}}{\lambda_n}\leq E_{\alpha,1}(-\lambda_nT_i^{\alpha})\leq-\frac{\underline{C_i}}{\lambda_n},\quad\forall n=1,2,\cdots,
\end{equation}
\begin{equation}\label{SQ_Q}
\frac{\underline{C_i}}{\lambda_n}\leq E_{\alpha,\alpha+1}(-\lambda_nT_i^{\alpha})\leq\frac{\overline{C_i}}{\lambda_n},\quad\forall n=1,2,\cdots,
\end{equation}
where the constants $\overline{C_i}$, $\underline{C_i}$ depend on $\alpha$ and $T_i$, $i=1,2$.
\end{lem}

Let $t=T_1,T_2$ in (\ref{solution2}), respectively. We have 
\begin{equation}\label{eq1}
u(T_1)=\mathbb{S}(T_1)a+\mathbb{Q}(T_1)f,\quad u(T_2)=\mathbb{S}(T_2)a+\mathbb{Q}(T_2)f.
\end{equation}
Then the above equations can be written into the following matrix type:
\begin{equation}\label{solu_map}
\begin{bmatrix}
\mathbb{S}(T_1)&\mathbb{Q}(T_1)\\
\mathbb{S}(T_2)&\mathbb{Q}(T_2)
\end{bmatrix}
\begin{bmatrix}
a\\
f
\end{bmatrix}=
\begin{bmatrix}
g_1\\
g_2
\end{bmatrix}.
\end{equation}

We know that there is a unique pair of solution $(a,f)$ to (\ref{solu_map}) provided the matrix $\begin{bmatrix}
\mathbb{S}(T_1)&\mathbb{Q}(T_1)\\
\mathbb{S}(T_2)&\mathbb{Q}(T_2)
\end{bmatrix}$ is nonsingular. Thus the solution can be represented as 
\begin{equation}\label{exactaf}
\begin{bmatrix}
a\\
f
\end{bmatrix}=\begin{bmatrix}
\mathbb{S}(T_1)&\mathbb{Q}(T_1)\\
\mathbb{S}(T_2)&\mathbb{Q}(T_2)
\end{bmatrix}^{-1}\begin{bmatrix}
g_1\\
g_2
\end{bmatrix}.
\end{equation}

Since $\{\varphi_n\}_{n=1}^{\infty}$ is an orthonormal basis in  $L^{2}(\Omega)$, $g_1,g_2\in L^{2}(\Omega)$ can be expressed by 
\begin{equation*}
g_1=\sum_{n=1}^{\infty}(g_1,\varphi_n)\varphi_n,\quad g_2=\sum_{n=1}^{\infty}(g_2,\varphi_n)\varphi_n.
\end{equation*}
Making use of the Cramer's Rule and spectral decomposition, we have
\begin{equation}\label{afexact}
\begin{bmatrix}
a\\
f
\end{bmatrix}=\sum_{n=1}^{\infty}\frac{1}{R_n(T_1,T_2)}\begin{bmatrix}
Q_n(T_2)&-Q_n(T_1)\\
S_n(T_2)&-S_n(T_1)
\end{bmatrix}\begin{bmatrix}
(g_1,\varphi_n)\varphi_n\\
(g_2,\varphi_n)\varphi_n
\end{bmatrix},
\end{equation}
where $R_n(T_1,T_2)$ are defined as
$$R_n(T_1,T_2):=S_n(T_1)Q_n(T_2)-S_n(T_2)Q_n(T_1).$$

In order to ensure the existence and uniqueness of $(a,f)$, we provide the following lemma, which guarantees the nonsingularity of $\begin{bmatrix}
\mathbb{S}(T_1)&\mathbb{Q}(T_1)\\
\mathbb{S}(T_2)&\mathbb{Q}(T_2)
\end{bmatrix}$.

\begin{lem}\label{SQQS_bound}
Let $1<\alpha<2$, then for sufficiently large $T_1>0$, there exists a constant $d<1$ depending on $\alpha$, $T_1$ and $T_2$ such that
\begin{equation}\label{SQQS1}
0<\frac{Q_n(T_1)S_n(T_2)}{S_n(T_1)Q_n(T_2)}\leq d,\quad n=1,2,\cdots.
\end{equation}
\end{lem}
\begin{proof}
By the asymptotic expansion of Mittag-Leffler function in Lemma \ref{lem2.2}, we have
\begin{eqnarray*}
&&S_n(T_i)=E_{\alpha,1}(-\lambda_nT_i^{\alpha})=\frac{1}{\Gamma(1-\alpha)}\frac{1}{\lambda_nT_i^{\alpha}}+O\left(\frac{1}{\lambda_n^2T_i^{2\alpha}}\right),\\
&&Q_n(T_i)=T_i^{\alpha}E_{\alpha,\alpha+1}(-\lambda_nT_i^{\alpha})=\frac{1}{\lambda_n}+O\left(\frac{1}{\lambda_n^2T_i^{\alpha}}\right).
\end{eqnarray*}
Since $T_2>T_1$, it follows that
\begin{eqnarray}\label{SQQS2}
R_n(T_1,T_2)&=&S_n(T_1)Q_n(T_2)-S_n(T_2)Q_n(T_1)\\\nonumber
&=&\frac{T_2^{\alpha}-T_1^{\alpha}}{\Gamma(1-\alpha)\lambda_n^2T_1^{\alpha}T_2^{\alpha}}+O\left(\frac{1}{\lambda_n^3T_1^{2\alpha}}\right).
\end{eqnarray}
It can be learned from Lemma \ref{lem2.1} that
\begin{equation*}
|S_n(T_i)|\leq\frac{C_1(\alpha)}{1+\lambda_nT_i^{\alpha}},\quad |Q_n(T_i)|\leq\frac{C_2(\alpha)T_i^{\alpha}}{1+\lambda_nT_i^{\alpha}}.
\end{equation*}
By the fact that $\Gamma(1-\alpha)<0$, we know that there exists a constant $M_0>0$ such that
\begin{equation}\label{SQQS3}
S_n(T_1)Q_n(T_2)-S_n(T_2)Q_n(T_1)\leq\frac{T_2^{\alpha}-T_1^{\alpha}}{\max(2,\frac{1}{C_1(\alpha)C_2(\alpha)|\Gamma(1-\alpha)|})\Gamma(1-\alpha)\lambda_n^{2}T_1^{\alpha}T_2^{\alpha}},
\end{equation}
provided $\lambda_n^3 T_1^{2\alpha}>M_0$. Divided by $S_n(T_1)Q_n(T_2)$ on both sides of (\ref{SQQS3}), it can be obtained that
\begin{eqnarray*}
1-\frac{S_n(T_2)Q_n(T_1)}{S_n(T_1)Q_n(T_2)}&\geq&\frac{T_2^{\alpha}-T_1^{\alpha}}{\max(2,\frac{1}{C_1(\alpha)C_2(\alpha)|\Gamma(1-\alpha)|})|\Gamma(1-\alpha)|\lambda_n^{2}T_1^{\alpha}T_2^{\alpha}|S_n(T_1)|Q_n(T_2)}\\
&\geq&\frac{1}{\max(2,\frac{1}{C_1(\alpha)C_2(\alpha)|\Gamma(1-\alpha)|})|\Gamma(1-\alpha)|\cdot\frac{C_1(\alpha)\lambda_nT_1^{\alpha}}{1+\lambda_nT_1^{\alpha}}\cdot\frac{C_2(\alpha)\lambda_nT_2^{\alpha}}{1+\lambda_nT_2^{\alpha}}}\left(1-\frac{T_1^{\alpha}}{T_2^{\alpha}}\right)\\
&\geq&\frac{1}{\max(2,\frac{1}{C_1(\alpha)C_2(\alpha)|\Gamma(1-\alpha)|})\cdot C_1(\alpha)C_2(\alpha)|\Gamma(1-\alpha)|}\left(1-\frac{T_1^{\alpha}}{T_2^{\alpha}}\right).
\end{eqnarray*}
If $\frac{1}{C_1(\alpha)C_2(\alpha)|\Gamma(1-\alpha)|}\geq2$, we have
\begin{equation*}
1-\frac{S_n(T_2)Q_n(T_1)}{S_n(T_1)Q_n(T_2)}\geq 1-\frac{T_1^{\alpha}}{T_2^{\alpha}},
\end{equation*}
which implies that 
\begin{equation*}
\frac{S_n(T_2)Q_n(T_1)}{S_n(T_1)Q_n(T_2)}\leq\frac{T_1^{\alpha}}{T_2^{\alpha}}<1.
\end{equation*}
If $\frac{1}{C_1(\alpha)C_2(\alpha)|\Gamma(1-\alpha)|}<2$, then $2C_1(\alpha)C_2(\alpha)|\Gamma(1-\alpha)|>1$ we have 
\begin{equation*}
1-\frac{S_n(T_2)Q_n(T_1)}{S_n(T_1)Q_n(T_2)}\geq\frac{1}{2 C_1(\alpha)C_2(\alpha)|\Gamma(1-\alpha)|}\left(1-\frac{T_1^{\alpha}}{T_2^{\alpha}}\right),
\end{equation*}
which implies that 
\begin{equation*}
\frac{S_n(T_2)Q_n(T_1)}{S_n(T_1)Q_n(T_2)}\leq 1-\frac{1}{2 C_1(\alpha)C_2(\alpha)|\Gamma(1-\alpha)|}\left(1-\frac{T_1^{\alpha}}{T_2^{\alpha}}\right)<1.
\end{equation*}
Therefore, we set 
\begin{equation}\label{d}
d:=\max\left(\frac{T_1^{\alpha}}{T_2^{\alpha}},1-\frac{1}{2 C_1(\alpha)C_2(\alpha)|\Gamma(1-\alpha)|}\left(1-\frac{T_1^{\alpha}}{T_2^{\alpha}}\right)\right),
\end{equation}
and  the proof is accomplished.
\end{proof}

To end this section, we provide the following conditional stability estimate in order to complete the convergence analysis.
\begin{lem}\label{constab}
For $a,f\in D((-\Delta)^p)$, we have the following two-sided estimate:
\begin{equation}\label{stability}
C_{1}\left(\|g_1\|_{D((-\Delta)^{p+1})}+\|g_2\|_{D((-\Delta)^{p+1})}\right)\leq\|a\|_{D((-\Delta)^p)}+\|f\|_{D((-\Delta)^p)}\leq C_{2}\left(\|g_1\|_{D((-\Delta)^{p+1})}+\|g_2\|_{D((-\Delta)^{p+1})}\right).
\end{equation}
\end{lem}
\noindent The proof is similar to that of Theorem $3.2$ in \cite{Zhang+Zhou-2022}, thus we omit it here.

\section{Alternating approach for reconstructing the initial value and source term}\label{sec3}
In this section, we shall provide the alternating method for reconstructing the initial value and source term in $(\ref{IP})$ and then analyze the convergence results. Firstly, we decouple the inverse problem into two subproblems. That is, by giving $f^{(0)}\in D((-\Delta)^p)$, we can obtain an approximate solution of $a$ by solving the following backward problem:
\begin{align}\label{SP1}
\left\{
\begin{aligned}
&\partial_{t}^{\alpha}v^{(k)}(t)-\Delta v^{(k)}(t)=f^{(k)} & &in\quad\Omega\times(0,T),\\
&v^{(k)}(t)=0& &on\quad\partial\Omega\times[0,T],\\
&v^{(k)}(0)=a^{(k)},~\partial_tv^{(k)}(0)=0& &in\quad\overline{\Omega},\\
&v^{(k)}(T_1)=g_1& &in\quad\overline{\Omega}.
\end{aligned}
\right.
\end{align}
The solution to the backward problem can be obtained as 
\begin{eqnarray}\label{ak}
a^{(k)}&=&\mathbb{S}(T_1)^{-1}\left(g_1-\mathbb{Q}(T_1)f^{(k)}\right)\\\nonumber
&=&\sum_{n=1}^{\infty}\frac{1}{S_n(T_1)}[(g_1,\varphi_n)-Q_n(T_1)(f^{(k)},\varphi_n)]\varphi_n,\quad k=0,1,2,\cdots.
\end{eqnarray}
With obtained $a^{(k)}$ in (\ref{ak}), we can obtain the approximate solution of $f$ by solving the following inverse source problem:
\begin{align}\label{SP2}
\left\{
\begin{aligned}
&\partial_{t}^{\alpha}w^{(k)}(t)-\Delta w^{(k)}(t)=f^{(k)} & &in\quad\Omega\times(0,T),\\
&w^{(k)}(t)=0& &on\quad\partial\Omega\times[0,T],\\
&w^{(k)}(0)=a^{(k-1)},~\partial_tw^{(k)}(0)=0& &in\quad\overline{\Omega},\\
&w^{(k)}(T_2)=g_2& &in\quad\overline{\Omega}.
\end{aligned}
\right.
\end{align}
The solution to the inverse source problem can be obtained as
\begin{eqnarray}\label{fk}
f^{(k)}&=&\mathbb{Q}(T_2)^{-1}(g_2-\mathbb{S}(T_2)a^{(k-1)})\\\nonumber
&=&\sum_{n=1}^{\infty}\frac{1}{Q_n(T_2)}[(g_2,\varphi_n)-S_n(T_2)(a^{(k-1)},\varphi_n)]\varphi_n,\quad k=1,2,\cdots.
\end{eqnarray}
Especially, when $k=0$, we note that 
\begin{equation}\label{a0}
a^{(0)}=\sum_{n=1}^{\infty}\frac{1}{S_n(T_1)}[(g_1,\varphi_n)-Q_n(T_1)(f^{(0)},\varphi_n)]\varphi_n.
\end{equation}
For $k>0$, by substituting (\ref{fk}) into (\ref{ak}), we obtain 
\begin{equation}\label{akk}
a^{(k)}=\sum_{n=1}^{\infty}\left[\frac{Q_n(T_1)S_n(T_2)}{S_n(T_1)Q_n(T_2)}(a^{(k-1)},\varphi_n)+\frac{1}{S_n(T_1)}(g_1,\varphi_n)-\frac{Q_n(T_1)}{S_n(T_1)Q_n(T_2)}(g_2,\varphi_n)\right]\varphi_n.
\end{equation}
Conversely, by substituting (\ref{ak}) into (\ref{fk}), we obtain 
\begin{equation}\label{fkk}
f^{(k)}=\sum_{n=1}^{\infty}\left[\frac{Q_n(T_1)S_n(T_2)}{S_n(T_1)Q_n(T_2)}(f^{(k-1)},\varphi_n)+\frac{1}{Q_n(T_2)}(g_2,\varphi_n)-\frac{S_n(T_2)}{Q_n(T_2)S_n(T_1)}(g_1,\varphi_n)\right]\varphi_n.
\end{equation}

The following Lemma shows that the series $\left(a^{(k)},f^{(k)}\right)$ converge to the exact solutions $(a,f)$.
\begin{lem}\label{SPlem2}
Let $a,f,f^{(0)}\in D((-\Delta)^p)$ for $p\in(0,1]$. $(a,f)$ are exact solutions of the inverse problem defined by $(\ref{afexact})$. Then iterated solutions $(a^{(k)},f^{(k)})$ defined by $(\ref{ak})$ and $(\ref{fk})$ satisfy the following estimate:
\begin{equation}
\left\|a^{(k)}-a\right\|\leq d^k\left\|a^{(0)}-a\right\|,\quad \left\|f^{(k)}-f\right\|\leq d^k\left\|f^{(0)}-f\right\|,
\end{equation}
where $a^{(0)}$ is defined by $(\ref{a0})$, $d<1$ is defined by (\ref{d}).
\end{lem}
\begin{proof}
We can learn from (\ref{afexact}) and (\ref{ak}) that 
\begin{eqnarray*}
\left(a^{(k)}-a,\varphi_n\right)&=&\frac{Q_n(T_1)S_n(T_2)}{S_n(T_1)Q_n(T_2)}\left(a^{(k-1)},\varphi_n\right)+\left[\frac{1}{S_n(T_1)}-\frac{Q_n(T_2)}{R_n(T_1,T_2)}\right]\left(g_1,\varphi_n\right)\\
&&-\left[\frac{Q_n(T_1)}{S_n(T_1)Q_n(T_2)}-\frac{Q_n(T_1)}{R_n(T_1,T_2)}\right]\left(g_2,\varphi_n\right).
\end{eqnarray*}
Thus 
\begin{eqnarray*}
a^{(k)}-a&=&\sum_{n=1}^{\infty}\left\{\frac{Q_n(T_1)S_n(T_2)}{S_n(T_1)Q_n(T_2)}\left(a^{(k-1)},\varphi_n\right)-\frac{Q_n(T_1)S_n(T_2)}{S_n(T_1)R_n(T_1,T_2)}\left(g_1,\varphi_n\right)\right.\\
&&\left.+\frac{Q_n^2(T_1)S_n(T_2)}{S_n(T_1)Q_n(T_2)R_n(T_1,T_2)}\left(g_2,\varphi_n\right)\right\}\varphi_n\\
&=&\sum_{n=1}^{\infty}\frac{Q_n(T_1)S_n(T_2)}{S_n(T_1)Q_n(T_2)}\left\{\left(a^{(k-1)},\varphi_n\right)-\frac{1}{R_n(T_1,T_2)}\left[Q_n(T_2)\left(g_1,\varphi_n\right)-Q_n(T_1)\left(g_2,\varphi_n\right)\right]\right\}\varphi_n\\
&=&\sum_{n=1}^{\infty}\frac{Q_n(T_1)S_n(T_2)}{S_n(T_1)Q_n(T_2)}\left\{\left(a^{(k-1)}-a,\varphi_n\right)\right\}\varphi_n=\cdots=\sum_{n=1}^{\infty}\left[\frac{Q_n(T_1)S_n(T_2)}{S_n(T_1)Q_n(T_2)}\right]^k\left(a^{(0)}-a,\varphi_n\right)\varphi_n.
\end{eqnarray*}
By Lemma \ref{SQQS_bound}, we know that
\begin{equation*}
\left\|a^{(k)}-a\right\|^2\leq d^{2k}\sum_{n=1}^{\infty}\left(a^{(0)}-a,\varphi_n\right)^2=d^{2k}\left\|a^{(0)}-a\right\|^2.
\end{equation*}
Similarly, we have
\begin{equation*}
f^{(k)}-f=\sum_{n=1}^{\infty}\left[\frac{Q_n(T_1)S_n(T_2)}{S_n(T_1)Q_n(T_2)}\right]^k\left(f^{(0)}-f,\varphi_n\right)\varphi_n
\end{equation*}
and
\begin{equation*}
\left\|f^{(k)}-f\right\|^2\leq d^{2k}\sum_{n=1}^{\infty}\left(f^{(0)}-f,\varphi_n\right)^2=d^{2k}\left\|f^{(0)}-f\right\|^2.
\end{equation*}
\end{proof}
It follows from (\ref{akk}) and (\ref{fkk}) that the iterated solutions $a^{(k)}$ and $f^{(k)}$ can be rewritten as
\begin{eqnarray}\label{rec_ak}
a^{(k)}&=&\sum_{n=1}^{\infty}\left[\frac{Q_n(T_1)S_n(T_2)}{S_n(T_1)Q_n(T_2)}(a^{(k-1)},\varphi_n)+\frac{1}{S_n(T_1)}(g_1,\varphi_n)-\frac{Q_n(T_1)}{S_n(T_1)Q_n(T_2)}(g_2,\varphi_n)\right]\varphi_n\\\nonumber
&=&\sum_{n=1}^{\infty}\frac{Q_n(T_1)S_n(T_2)}{S_n(T_1)Q_n(T_2)}\left[\frac{Q_n(T_1)S_n(T_2)}{S_n(T_1)Q_n(T_2)}(a^{(k-2)},\varphi_n)+\frac{1}{S_n(T_1)}(g_1,\varphi_n)-\frac{Q_n(T_1)}{S_n(T_1)Q_n(T_2)}(g_2,\varphi_n)\right]\varphi_n\\\nonumber
&&+\sum_{n=1}^{\infty}\left[\frac{1}{S_n(T_1)}(g_1,\varphi_n)-\frac{Q_n(T_1)}{S_n(T_1)Q_n(T_2)}(g_2,\varphi_n)\right]\varphi_n\\\nonumber
&=&\sum_{n=1}^{\infty}\left\{\left[\frac{Q_n(T_1)S_n(T_2)}{S_n(T_1)Q_n(T_2)}\right]^2\left(a^{(k-2)},\varphi_n\right)\right.\\\nonumber
&&\left.+\left[1+\frac{Q_n(T_1)S_n(T_2)}{S_n(T_1)Q_n(T_2)}\right]\left[\frac{1}{S_n(T_1)}(g_1,\varphi_n)-\frac{Q_n(T_1)}{S_n(T_1)Q_n(T_2)}(g_2,\varphi_n)\right]\right\}\varphi_n\\\nonumber
&=&\sum_{n=1}^{\infty}\left\{\left[\frac{Q_n(T_1)S_n(T_2)}{S_n(T_1)Q_n(T_2)}\right]^k\left(a^{(0)},\varphi_n\right)+M_{k-1}\left[\frac{1}{S_n(T_1)}(g_1,\varphi_n)-\frac{Q_n(T_1)}{S_n(T_1)Q_n(T_2)}(g_2,\varphi_n)\right]\right\}\varphi_n
\end{eqnarray}
and
\begin{equation}\label{rec_fk}
f^{(k)}=\sum_{n=1}^{\infty}\left\{\left[\frac{Q_n(T_1)S_n(T_2)}{S_n(T_1)Q_n(T_2)}\right]^k\left(f^{(0)},\varphi_n\right)+M_{k-1}\left[\frac{1}{Q_n(T_2)}(g_2,\varphi_n)-\frac{S_n(T_2)}{S_n(T_1)Q_n(T_2)}(g_1,\varphi_n)\right]\right\}\varphi_n,
\end{equation}
where $M_{k}:=\sum_{j=0}^{k}\left(\frac{Q_n(T_1)S_n(T_2)}{S_n(T_1)Q_n(T_2)}\right)^j$ and it satisfies the following estimate:
\begin{equation*}
M_{k}=\left[1-\left(\frac{Q_n(T_1)S_n(T_2)}{S_n(T_1)Q_n(T_2)}\right)^{k+1}\right]\Big/\left(1-\frac{Q_n(T_1)S_n(T_2)}{S_n(T_1)Q_n(T_2)}\right)
\leq\frac{1}{1-d}.
\end{equation*}
Then we have the following results for iterated solutions$\left(a^{(k)},f^{(k)}\right)$.
\begin{lem}\label{reg_afk}
Let $a,f,f^{(0)}\in D((-\Delta)^p)$ for $p\in[0,1]$. Then we have
\begin{equation}\label{afkest}
\|a^{(k)}\|_{D((-\Delta)^p)}+\|f^{(k)}\|_{D((-\Delta)^p)}\leq C\left(d^k+1\right)\left(\|a\|_{D((-\Delta)^p)}+\|f\|_{D((-\Delta)^p)}\right)+Cd^k\|f^{(0)}\|_{D((-\Delta)^p)}.
\end{equation}
\end{lem}
\begin{proof}
It can be learned from (\ref{rec_ak}) that
\begin{eqnarray*}
\|a^{(k)}\|_{D((-\Delta)^p)}^2&\leq& 2\sum_{n=1}^{\infty}\left[\frac{Q_n(T_1)S_n(T_2)}{S_n(T_1)Q_n(T_2)}\right]^{2k}\lambda_n^{2p}\left(a^{(0)},\varphi_n\right)^2\\
&&+4M_{k-1}^2\sum_{n=1}^{\infty}\lambda_n^{2p}\left[\frac{1}{S_n^2(T_1)}\left(g_1,\varphi_n\right)^2+\frac{Q_n^2(T_1)}{S_n^2(T_1)Q_n^2(T_2)}\left(g_2,\varphi_n\right)^2\right]^{2}.
\end{eqnarray*}
By (\ref{a0}), we can obtain from H{\"o}lder's inequality that
\begin{equation*}
\|a^{(0)}\|_{D((-\Delta)^p)}^2\leq 2\sum_{n=1}^{\infty}\frac{\lambda_n^{2p}}{S_n^2(T_1)}\left[(g_1,\varphi_n)^2+Q_n^2(T_1)(f^{(0)},\varphi_n)^2\right].
\end{equation*}
Then we have 
\begin{eqnarray*}
\|a^{(k)}\|_{D((-\Delta)^p)}^2&\leq& \left[2d^{2k}+\frac{4}{(1-d)^2}\right]\sum_{n=1}^{\infty}\frac{\lambda_n^{2p}}{S_n^2(T_1)}(g_1,\varphi_n)^2\\
&&+2d^{2k}\sum_{n=1}^{\infty}\frac{Q_n^2(T_1)}{S_n^2(T_1)}\lambda_n^{2p}(f^{(0)},\varphi_n)^2+\frac{4}{(1-d)^2}\sum_{n=1}^{\infty}\frac{\lambda_n^{2p}}{Q_n^2(T_2)}(g_2,\varphi_n)^2\\
&\leq&C\left(d^{2k}+1\right)\sum_{n=1}^{\infty}\lambda_n^{2p+2}(g_1,\varphi_n)^2+Cd^{2k}\sum_{n=1}^{\infty}\lambda_n^{2p}(f^{(0)},\varphi_n)^2+C\sum_{n=1}^{\infty}\lambda_n^{2p+2}(g_2,\varphi_n)^2\\
&\leq&Cd^{2k}\left(\|a\|_{D((-\Delta)^p)}^2+\|f\|_{D((-\Delta)^p)}^2+\|f^{(0)}\|_{D((-\Delta)^p)}^2\right)+C\left(\|a\|_{D((-\Delta)^p)}^2+\|f\|_{D((-\Delta)^p)}^2\right).
\end{eqnarray*}
Similarly, we can deduce that
\begin{eqnarray*}
\|f^{(k)}\|_{D((-\Delta)^p)}^2&\leq& 2\sum_{n=1}^{\infty}\left[\frac{Q_n(T_1)S_n(T_2)}{S_n(T_1)Q_n(T_2)}\right]^{2k}\lambda_n^{2p}\left(f^{(0)},\varphi_n\right)^2\\
&&+4M_{k-1}^2\sum_{n=1}^{\infty}\lambda_n^{2p}\left[\frac{1}{Q_n^2(T_2)}\left(g_2,\varphi_n\right)^2+\frac{S_n^2(T_2)}{S_n^2(T_1)Q_n^2(T_2)}\left(g_1,\varphi_n\right)^2\right]^{2}\\
&\leq&2d^{2k}\|f^{(0)}\|_{D((-\Delta)^p)}^2+C\left(\|g_1\|_{D((-\Delta)^{p+1})}^2+\|g_2\|_{D((-\Delta)^{p+1})}^2\right)\\
&\leq&2d^{2k}\|f^{(0)}\|_{D((-\Delta)^p)}^2+C\left(\|a\|_{D((-\Delta)^{p})}^2+\|f\|_{D((-\Delta)^{p})}^2\right).
\end{eqnarray*}
Finally, we obtain that 
\begin{equation*}
\|a^{(k)}\|_{D((-\Delta)^p)}+\|f^{(k)}\|_{D((-\Delta)^p)}\leq C\left(d^k+1\right)\left(\|a\|_{D((-\Delta)^p)}+\|f\|_{D((-\Delta)^p)}\right)+Cd^k\|f^{(0)}\|_{D((-\Delta)^p)}.
\end{equation*}
\end{proof}

It is known that problems (\ref{SP1}) and (\ref{SP2}) are ill-posed. To obtain stable approximate solutions, we apply the quasi-boundary value method for both problems (\ref{SP1}) and (\ref{SP2}). That is, by setting $f_{\mu}^{(0)}=f^{(0)}$, we consider solving the following regularized problem:
\begin{align}\label{RSP1}
\left\{
\begin{aligned}
&\partial_{t}^{\alpha}v^{(k)}(t)-\Delta v^{(k)}(t)=f_{\mu}^{(k)} & &in\quad\Omega\times(0,T),\\
&v^{(k)}(t)=0& &on\quad\partial\Omega\times[0,T],\\
&v^{(k)}(0)=a_{\mu}^{(k)},~\partial_tv^{(k)}(0)=0& &in\quad\overline{\Omega},\\
&v^{(k)}(T_1)-\mu v^{(k)}(0)=g_1& &in\quad\overline{\Omega}.
\end{aligned}
\right.
\end{align}
The solution to (\ref{RSP1}) can be obtained as
\begin{eqnarray}\label{amuk}
a_{\mu}^{(k)}&=&(\mathbb{S}(T_1)-\mu\mathbb{I})^{-1}\left(g_1-\mathbb{Q}(T_1)f_{\mu}^{(k)}\right),\quad k=0,1,2,\cdots\\\nonumber
&=&\sum_{n=1}^{\infty}\frac{1}{S_n(T_1)-\mu}\left[(g_1,\varphi_n)-Q_n(T_1)(f_{\mu}^{(k)},\varphi_n)\right]\varphi_n.
\end{eqnarray}
With obtained $a_{\mu}^{(k)}$ in (\ref{amuk}), we consider solving the following regularized problem:
\begin{align}\label{RSP2}
\left\{
\begin{aligned}
&\partial_{t}^{\alpha}w^{(k)}(t)-\Delta w^{(k)}(t)=f_{\mu}^{(k)} & &in\quad\Omega\times(0,T),\\
&w^{(k)}(t)=0& &on\quad\partial\Omega\times[0,T],\\
&w^{(k)}(0)=a_{\mu}^{(k-1)},~\partial_tv^{(k)}(0)=0& &in\quad\overline{\Omega},\\
&w^{(k)}(T_2)+\mu f_{\mu}^{(k)}=g_2& &in\quad\overline{\Omega}.
\end{aligned}
\right.
\end{align}
The solution to (\ref{RSP2}) can be obtained as
\begin{eqnarray}\label{fmuk}
f_{\mu}^{(k)}&=&(\mathbb{Q}(T_2)+\mu\mathbb{I})^{-1}\left(g_2-\mathbb{S}(T_2)a_{\mu}^{(k-1)}\right),\quad k=1,2,\cdots\\\nonumber
&=&\sum_{n=1}^{\infty}\frac{1}{Q_n(T_2)+\mu}\left[(g_2,\varphi_n)-S_n(T_2)(a_{\mu}^{(k-1)},\varphi_n)\right]\varphi_n.
\end{eqnarray}
For $k>0$, by substituting (\ref{fmuk}) into (\ref{amuk}), we obtain 
\begin{eqnarray}\label{amukk}
a_{\mu}^{(k)}&=&(\mathbb{S}(T_1)-\mu\mathbb{I})^{-1}\left[g_1-\mathbb{Q}(T_1)(\mathbb{Q}(T_2)+\mu\mathbb{I})^{-1}\left(g_2-\mathbb{S}(T_2)a_{\mu}^{(k-1)}\right)\right]\\\nonumber
&=&\sum_{n=1}^{\infty}\left[\frac{Q_n(T_1)S_n(T_2)}{(S_n(T_1)-\mu)(Q_n(T_2)+\mu)}(a_{\mu}^{(k-1)},\varphi_n)+\frac{1}{S_n(T_1)-\mu}(g_1,\varphi_n)\right.\\\nonumber
&&\left.-\frac{Q_n(T_1)}{(S_n(T_1)-\mu)(Q_n(T_2)+\mu)}(g_2,\varphi_n)\right]\varphi_n.
\end{eqnarray}
Conversely, by substituting (\ref{amuk}) into (\ref{fmuk}), we obtain 
\begin{eqnarray}\label{fmukk}
f_{\mu}^{(k)}&=&(\mathbb{Q}(T_2)+\mu\mathbb{I})^{-1}\left[g_2-\mathbb{S}(T_2)(\mathbb{S}(T_1)-\mu\mathbb{I})^{-1}\left(g_1-\mathbb{Q}(T_1)f_{\mu}^{(k-1)}\right)\right]\\\nonumber
&=&\sum_{n=1}^{\infty}\left[\frac{Q_n(T_1)S_n(T_2)}{(S_n(T_1)-\mu)(Q_n(T_2)+\mu)}(f_{\mu,\delta}^{(k-1)},\varphi_n)+\frac{1}{Q_n(T_2)+\mu}(g_2,\varphi_n)\right.\\\nonumber
&&\left.-\frac{S_n(T_2)}{(S_n(T_1)-\mu)(Q_n(T_2)+\mu)}(g_1,\varphi_n)\right]\varphi_n.
\end{eqnarray}
The regularized solutions $(a_{\mu}^{k},f_{\mu}^{(k)})$ satisfy the following estimates.
\begin{lem}\label{reg_regusolu}
Let $a,f,f^{(0)}\in D((-\Delta)^p)$ for $p\in[0,1]$. Then for $q\in[p,p+1]$, we have
\begin{equation*}
\left\|a_{\mu}^{(k)}\right\|_{D((-\Delta)^q)}+\left\|f_{\mu}^{(k)}\right\|_{D((-\Delta)^q)}\leq C\mu^{p-q}\left(d^k+1\right)\left(\|a\|_{D((-\Delta)^p)}+\|f\|_{D((-\Delta)^p)}\right)+C\mu^{p-q}d^k\|f^{(0)}\|_{D((-\Delta)^p)}.
\end{equation*}
\end{lem}
\begin{proof}
By (\ref{amuk}), we know that
\begin{equation*}
a_{\mu}^{(0)}=\sum_{n=1}^{\infty}\frac{1}{S_n(T_1)-\mu}\left[(g_1,\varphi_n)-Q_n(T_1)(f^{(0)},\varphi_n)\right]\varphi_n,
\end{equation*}
thus we have
\begin{eqnarray*}
\|a_{\mu}^{(0)}\|_{D((-\Delta)^q)}^2&\leq&C\sum_{n=1}^{\infty}\left[\frac{\lambda_n^{2+2q}}{(\underline{C_1}+\mu\lambda_n)^2}(g_1,\varphi_n)^2+\frac{\lambda_n^{2q}}{(\underline{C_1}+\mu\lambda_n)^2}(f^{(0)},\varphi_n)^2\right]\\
&\leq&C\sum_{n=1}^{\infty}\left[\frac{\lambda_n^{2(q-p)}}{(\underline{C_1}+\mu\lambda_n)^2}\lambda_n^{2(p+1)}(g_1,\varphi_n)^2+\frac{\lambda_n^{2(q-p)}}{(\underline{C_1}+\mu\lambda_n)^2}\lambda_n^{2p}(f^{(0)},\varphi_n)^2\right]\\
&\leq&C\sum_{n=1}^{\infty}\frac{(\mu\lambda_n)^{2(q-p)}}{(\underline{C_1}+\mu\lambda_n)^{2(q-p)}}\cdot\frac{\mu^{2(p-q)}}{(\underline{C_1}+\mu\lambda_n)^{2+2(p-q)}}\left[\lambda_n^{2(p+1)}(g_1,\varphi_n)^2+\lambda_n^{2p}(f^{(0)},\varphi_n)^2\right]\\
&\leq&C\mu^{2(p-q)}\sum_{n=1}^{\infty}\left(\|g_1\|_{D((-\Delta)^{p+1})}^2+\|f^{(0)}\|_{D((-\Delta)^{p})}^2\right)\\
&\leq&C\mu^{2(p-q)}\left(\|a\|_{D((-\Delta)^p)}^2+\|f\|_{D((-\Delta)^p)}^2+\|f^{(0)}\|_{D((-\Delta)^p)}^2\right).
\end{eqnarray*}
It can be learned from Lemma \ref{SQQS_bound} that 
\begin{equation}\label{QSmu}
0<\frac{Q_n(T_1)S_n(T_2)}{(S_n(T_1)-\mu)(Q_n(T_2)+\mu)}<\frac{Q_n(T_1)S_n(T_2)}{S_n(T_1)Q_n(T_2)}\leq d<1.
\end{equation}
Consequently, we have
\begin{eqnarray*}
\left\|a_{\mu}^{(1)}\right\|_{D((-\Delta)^q)}&\leq& d\left\|a_{\mu}^{0}\right\|_{D((-\Delta)^q)}+C\mu^{p-q}\left(\|a\|_{D((-\Delta)^p)}+\|f\|_{D((-\Delta)^p)}\right)\\
&\leq&C\mu^{p-q}(d+1)\left(\|a\|_{D((-\Delta)^p)}+\|f\|_{D((-\Delta)^p)}\right)+C\mu^{p-q}d\left\|f^{0}\right\|_{D((-\Delta)^p)}.
\end{eqnarray*}
Next we assume that for $j=2,\cdots,k$, there holds that
\begin{equation*}
\|a_{\mu}^{(j)}\|_{D((-\Delta)^q)}\leq C\mu^{p-q}\left(d^j+1\right)\left(\|a\|_{D((-\Delta)^p)}+\|f\|_{D((-\Delta)^p)}\right)+C\mu^{p-q}d^j\|f^{(0)}\|_{D((-\Delta)^p)}.
\end{equation*}
Then for $j=k+1$, by repeating above arguments, we have
\begin{eqnarray*}
\|a_{\mu}^{(k+1)}\|_{D((-\Delta)^q)}&\leq&d\left\|a_{\mu}^{(k)}\right\|_{D((-\Delta)^q)}+C\mu^{p-q}\left(\|a\|_{D((-\Delta)^p)}+\|f\|_{D((-\Delta)^p)}\right)\\
&\leq&C\mu^{p-q}\left(d^{k+1}+1\right)\left(\|a\|_{D((-\Delta)^p)}+\|f\|_{D((-\Delta)^p)}\right)+C\mu^{p-q}d^{k+1}\|f^{(0)}\|_{D((-\Delta)^p)}\\
&&+C\mu^{p-q}\left(\|a\|_{D((-\Delta)^p)}+\|f\|_{D((-\Delta)^p)}\right)\\
&\leq&C\mu^{p-q}\left(d^{k+1}+1\right)\left(\|a\|_{D((-\Delta)^p)}+\|f\|_{D((-\Delta)^p)}\right)+C\mu^{p-q}d^{k+1}\|f^{(0)}\|_{D((-\Delta)^p)}.
\end{eqnarray*}
Similarly, we can obtain that
\begin{equation*}
\|f_{\mu}^{(k)}\|_{D((-\Delta)^q)}\leq C\mu^{p-q}(d^k+1)\left(\|a\|_{D((-\Delta)^p)}+\|f\|_{D((-\Delta)^p)}\right)+C\mu^{p-q}d^k\|f^{(0)}\|_{D((-\Delta)^p)},
\end{equation*}
for all $k$. Then the proof is accomplished.
\end{proof}
The following result can be obtained directly from above lemma.
\begin{lem}\label{reg_iter_solu}
For $a,f,f^{(0)}\in D((-\Delta)^p)$ for $p\in[0,1]$. Then for $q\in[p,p+1]$, the solutions to problems (\ref{RSP1}) and (\ref{RSP2}) satisfy the following estimate:
\begin{equation}\label{reg_iter_solu1}
\left\|v^{(k)}(t)\right\|_{D((-\Delta)^q)}+\left\|w^{(k)}(t)\right\|_{D((-\Delta)^q)}\leq C\mu^{p-q}\left(d^{k+1}+1\right)\left(\|a\|_{D((-\Delta)^p)}+\|f\|_{D((-\Delta)^p)}\right)+C\mu^{p-q}d^k\|f^{(0)}\|_{D((-\Delta)^p)}.
\end{equation}
\end{lem}
\begin{proof}
It can be easily deduced that the solutions to problems (\ref{RSP1}) and (\ref{RSP2}) are
\begin{eqnarray*}
v^{(k)}(t)&=&\mathbb{S}(t)a_{\mu}^{(k)}+\mathbb{Q}(t)f_{\mu}^{(k)},\\
w^{(k)}(t)&=&\mathbb{S}(t)a_{\mu}^{(k-1)}+\mathbb{Q}(t)f_{\mu}^{(k)}.
\end{eqnarray*}
Then (\ref{reg_iter_solu1}) can be obtained directly from Lemmas \ref{lem2.1} and \ref{reg_regusolu}.
\end{proof}

\begin{lem}\label{RSPlem1}
Let $a,f,f^{(0)}\in D((-\Delta)^p)$ for $p\in(0,1]$. $(a^{(k)},f^{(k)})$ are given by (\ref{ak}) and (\ref{fk}). $(a_{\mu}^{(k)},f_{\mu}^{(k)})$ are defined by (\ref{amuk}) and (\ref{fmuk}). Then we have
\begin{equation}\label{amukak}
\left\|a_{\mu}^{(k)}-a^{(k)}\right\|+\left\|f_{\mu}^{(k)}-f^{(k)}\right\|\leq C\mu^{p}\left(d^k+1\right)\left(\|a\|_{D((-\Delta)^p)}+\|f\|_{D((-\Delta)^p)}\right)+C\mu^pd^k\|f^{(0)}\|_{D((-\Delta)^p)}.
\end{equation}
\end{lem}
\begin{proof}
By a recursive approach, it follows from (\ref{amukk}) and (\ref{akk}) that
\begin{eqnarray*}
a_{\mu}^{(k)}-a^{(k)}&=&\sum_{n=1}^{\infty}\frac{Q_n(T_1)S_n(T_2)}{(S_n(T_1)-\mu)(Q_n(T_2)+\mu)}(a_{\mu}^{(k-1)}-a^{(k-1)},\varphi_n)\varphi_n\\
&&+\left[\frac{Q_n(T_1)S_n(T_2)}{(S_n(T_1)-\mu)(Q_n(T_2)+\mu)}-\frac{Q_n(T_1)S_n(T_2)}{S_n(T_1)Q_n(T_2)}\right](a^{(k-1)},\varphi_n)\varphi_n\\
&&+\left[\frac{1}{S_n(T_1)-\mu}-\frac{1}{S_n(T_1)}\right](g_1,\varphi_n)\varphi_n\\
&&-\left[\frac{Q_n(T_1)}{(S_n(T_1)-\mu)(Q_n(T_2)+\mu)}-\frac{Q_n(T_1)}{S_n(T_1)Q_n(T_2)}\right](g_2,\varphi_n)\varphi_n\\
&=:&\sum_{n=1}^{\infty}\frac{Q_n(T_1)S_n(T_2)}{(S_n(T_1)-\mu)(Q_n(T_2)+\mu)}(a_{\mu}^{(k-1)}-a^{(k-1)},\varphi_n)\varphi_n+J_1+J_2+J_3.
\end{eqnarray*}
Next, we divide it into four steps to accomplish the proof. \\
\emph{Step 1:} Estimation of $J_1$. It can be easily deduced that
\begin{eqnarray*}
\|J_1\|^2&=&\sum_{n=1}^{\infty}\left[\frac{Q_n(T_1)S_n(T_2)}{S_n(T_1)Q_n(T_2)}\right]^2\left[\frac{\mu(Q_n(T_2)-S_n(T_1))+\mu^2}{(S_n(T_1)-\mu)(Q_n(T_2)+\mu)}\right]^2(a^{(k-1)},\varphi_n)^2\\
&\leq&Cd^2\sum_{n=1}^{\infty}\left[\frac{\mu\lambda_n+\mu^2\lambda_n^2}{(\underline{C_1}+\mu\lambda_n)(\underline{C_2}+\mu\lambda_n)}\right]^2(a^{(k-1)},\varphi_n)^2\\
&\leq&Cd^2\sum_{n=1}^{\infty}\left[\frac{\mu\lambda_n^{1-p}+\mu^2\lambda_n^{2-p}}{(\underline{C_1}+\mu\lambda_n)(\underline{C_2}+\mu\lambda_n)}\right]^2\lambda_n^{2p}(a^{(k-1)},\varphi_n)^2\\
&\leq&Cd^2\mu^{2p}\|a^{(k-1)}\|_{D((-\Delta)^p)}^2.
\end{eqnarray*}
By (\ref{afkest}) in Lemma \ref{reg_afk}, we have
\begin{equation*}
\|J_1\|\leq C\mu^p\left(d^{k}+1\right)\left(\|a\|_{D((-\Delta)^p)}+\|f\|_{D((-\Delta)^p)}\right)+C\mu^pd^{k}\|f^{(0)}\|_{D((-\Delta)^p)}.
\end{equation*}
\emph{Step 2:} Estimation of $J_2$.\\
\begin{eqnarray*}
\|J_2\|^2&=&\sum_{n=1}^{\infty}\left[\frac{\mu}{S_n(T_1)(S_n(T_1)-\mu)}\right]^2(g_1,\varphi_n)^2\leq C\sum_{n=1}^{\infty}\left[\frac{\mu\lambda_n^{1-p}}{\underline{C_1}(\underline{C_1}+\mu\lambda_n)}\right]^2\lambda_n^{2p+2}(g_1,\varphi_n)^2\\
&\leq&C\mu^{2p}\sum_{n=1}^{\infty}\lambda_n^{2p+2}(g_1,\varphi_n)^2\leq C\mu^{2p}\left(\|a\|_{D((-\Delta)^p)}^2+\|f\|_{D((-\Delta)^p)}^2\right).
\end{eqnarray*}
\emph{Step 3:} Estimation of $J_3$.\\
\begin{eqnarray*}
\|J_3\|^2&=&\sum_{n=1}^{\infty}\left[\frac{Q_n(T_1)S_n(T_2)}{S_n(T_1)Q_n(T_2)}\right]^{2}\frac{1}{S_n(T_2)^2}\left[\frac{\mu(Q_n(T_2)-S_n(T_1))+\mu^2}{(S_n(T_1)-\mu)(Q_n(T_2)+\mu)}\right]^2(g_2,\varphi_n)^2\\
&\leq&Cd^2\sum_{n=1}^{\infty}\left[\frac{\mu\lambda_n+\mu^2\lambda_n^2}{(\underline{C_1}+\mu\lambda_n)(\underline{C_2}+\mu\lambda_n)}\right]^2\lambda_n^2(g_2,\varphi_n)^2\\
&\leq&Cd^2\mu^{2p}\sum_{n=1}^{\infty}\lambda_n^{2p+2}(g_2,\varphi_n)^2\leq Cd^2\mu^{2p}\left(\|a\|_{D((-\Delta)^p)}^2+\|f\|_{D((-\Delta)^p)}^2\right).
\end{eqnarray*}
\emph{Step 4:} Estimation of $a_{\mu}^{(k)}-a^{(k)}$. We accomplish the proof by  the method of induction. By (\ref{amuk}) and (\ref{a0}), we have
\begin{eqnarray*}
\|a_{\mu}^{(0)}-a^{(0)}\|^2&=&\sum_{n=1}^{\infty}\left[\frac{1}{S_n(T_1)-\mu}-\frac{1}{S_n(T_1)}\right]^2\left[(g_1,\varphi_n)-Q_n(T_1)(f^{(0)},\varphi_n)\right]^2\\
&\leq&2\sum_{n=1}^{\infty}\left[\frac{\mu}{S_n(T_1)(S_n(T_1)-\mu)}\right]^2\left[(g_1,\varphi_n)^2+Q_n(T_1)^2(f^{(0)},\varphi_n)^2\right]\\
&\leq&2\sum_{n=1}^{\infty}\left[\frac{\mu\lambda_n^{1-p}}{\underline{C_1}(\underline{C_1}+\mu\lambda_n)}\right]^2\lambda_n^{2p+2}\left[(g_1,\varphi_n)^2+Q_n(T_1)^2(f^{(0)},\varphi_n)^2\right]\\
&\leq&C\mu^{2p}\left(\|a\|_{D((-\Delta)^p)}^2+\|f\|_{D((-\Delta)^p)}^2+\|f^{(0)}\|_{D((-\Delta)^p)}^2\right).
\end{eqnarray*}
Moreover, we have
\begin{eqnarray*}
a_{\mu}^{(1)}-a^{(1)}&=&\sum_{n=1}^{\infty}\frac{Q_n(T_1)S_n(T_2)}{(S_n(T_1)-\mu)(Q_n(T_2)+\mu)}(a_{\mu}^{(0)}-a^{(0)},\varphi_n)\varphi_n\\
&&+\left[\frac{Q_n(T_1)S_n(T_2)}{(S_n(T_1)-\mu)(Q_n(T_2)+\mu)}-\frac{Q_n(T_1)S_n(T_2)}{S_n(T_1)Q_n(T_2)}\right](a^{(0)},\varphi_n)\varphi_n\\
&&+\left[\frac{1}{S_n(T_1)-\mu}-\frac{1}{S_n(T_1)}\right](g_1,\varphi_n)\varphi_n\\
&&-\left[\frac{Q_n(T_1)}{(S_n(T_1)-\mu)(Q_n(T_2)+\mu)}-\frac{Q_n(T_1)}{S_n(T_1)Q_n(T_2)}\right](g_2,\varphi_n)\varphi_n.
\end{eqnarray*}
By using (\ref{QSmu}) and repeating above arguments, we can easily obtain that
\begin{eqnarray*}
\|a_{\mu}^{(1)}-a^{(1)}\|&\leq&Cd\|a_{\mu}^{(0)}-a^{(0)}\|+C\mu^pd\left(\|a\|_{D((-\Delta)^p)}+\|f\|_{D((-\Delta)^p)}+\|f^{(0)}\|_{D((-\Delta)^p)}\right)\\
&&+C\mu^p\left(\|a\|_{D((-\Delta)^p)}+\|f\|_{D((-\Delta)^p)}\right)
+C\mu^pd\left(\|a\|_{D((-\Delta)^p)}+\|f\|_{D((-\Delta)^p)}\right)\\
&\leq& C\mu^{p}\left(d+1\right)\left(\|a\|_{D((-\Delta)^p)}+\|f\|_{D((-\Delta)^p)}\right)+C\mu^pd\|f^{(0)}\|_{D((-\Delta)^p)}.
\end{eqnarray*}
Next we assume that for $j=2,\cdots,k$, we have the following estimate:
\begin{equation}\label{aj_induc}
\|a_{\mu}^{(j)}-a^{(j)}\|\leq C\mu^{p}\left(d^j+1\right)\left(\|a\|_{D((-\Delta)^p)}+\|f\|_{D((-\Delta)^p)}\right)+C\mu^pd^j\|f^{(0)}\|_{D((-\Delta)^p)}.
\end{equation}
Then it holds that
\begin{eqnarray*}
a_{\mu}^{(k+1)}-a^{(k+1)}&=&\sum_{n=1}^{\infty}\frac{Q_n(T_1)S_n(T_2)}{(S_n(T_1)-\mu)(Q_n(T_2)+\mu)}(a_{\mu}^{(k)}-a^{(k)},\varphi_n)\varphi_n\\
&&+\left[\frac{Q_n(T_1)S_n(T_2)}{(S_n(T_1)-\mu)(Q_n(T_2)+\mu)}-\frac{Q_n(T_1)S_n(T_2)}{S_n(T_1)Q_n(T_2)}\right](a^{(k)},\varphi_n)\varphi_n\\
&&+\left[\frac{1}{S_n(T_1)-\mu}-\frac{1}{S_n(T_1)}\right](g_1,\varphi_n)\varphi_n\\
&&-\left[\frac{Q_n(T_1)}{(S_n(T_1)-\mu)(Q_n(T_2)+\mu)}-\frac{Q_n(T_1)}{S_n(T_1)Q_n(T_2)}\right](g_2,\varphi_n)\varphi_n.
\end{eqnarray*}
Repeating above arguments, we have
\begin{eqnarray*}
\|a_{\mu}^{(k+1)}-a^{(k+1)}\|&\leq&Cd\|a_{\mu}^{(k)}-a^{(k)}\|+Cd\left(d^k+1\right)\mu^p\left(\|a\|_{D((-\Delta)^p)}+\|f\|_{D((-\Delta)^p)}\right)+Cd^{k+1}\mu^p\|f^{(0)}\|_{D((-\Delta)^p)}\\
&&+C\mu^p\left(\|a\|_{D((-\Delta)^p)}+\|f\|_{D((-\Delta)^p)}\right)
+Cd\mu^p\left(\|a\|_{D((-\Delta)^p)}+\|f\|_{D((-\Delta)^p)}\right).
\end{eqnarray*}
By using the inequality (\ref{aj_induc}), we have
\begin{eqnarray*}
\|a_{\mu}^{(k+1)}-a^{(k+1)}\|&\leq&C\mu^{p}d\left[\left(d^k+1\right)\left(\|a\|_{D((-\Delta)^p)}+\|f\|_{D((-\Delta)^p)}\right)+d^k\|f^{(0)}\|_{D((-\Delta)^p)}\right]\\
&&+Cd\left(d^k+1\right)\mu^p\left(\|a\|_{D((-\Delta)^p)}+\|f\|_{D((-\Delta)^p)}\right)+Cd^{k+1}\mu^p\|f^{(0)}\|_{D((-\Delta)^p)}\\
&&+C\mu^p\left(\|a\|_{D((-\Delta)^p)}+\|f\|_{D((-\Delta)^p)}\right)
+Cd\mu^p\left(\|a\|_{D((-\Delta)^p)}+\|f\|_{D((-\Delta)^p)}\right)\\
&\leq&C\mu^{p}\left(d^{k+1}+1\right)\left(\|a\|_{D((-\Delta)^p)}+\|f\|_{D((-\Delta)^p)}\right)+C\mu^pd^{k+1}\|f^{(0)}\|_{D((-\Delta)^p)},
\end{eqnarray*}
which means that estimate (\ref{aj_induc}) also holds for $j=k+1$. Thus we can learn that (\ref{aj_induc}) holds for all $k\in\mathbb{N}$.

Similarly, we can obtain the following estimate for $f_{\mu}^{(k)}-f^{(k)}$:
\begin{equation*}
\|f_{\mu}^{(k)}-f^{(k)}\|\leq C\mu^{p}\left(d^{k}+1\right)\left(\|a\|_{D((-\Delta)^p)}+\|f\|_{D((-\Delta)^p)}\right)+C\mu^pd^{k}\|f^{(0)}\|_{D((-\Delta)^p)}.
\end{equation*}
Finally, we can obtain the estimate (\ref{amukak}).
\end{proof}

\begin{lem}\label{noise_err}
Let $a,f,f^{(0)}\in D((-\Delta)^p)$ for $p\in(0,1]$. $(a_{\mu,\delta}^{(k)},f_{\mu,\delta}^{(k)})$ are defined in (\ref{amukk}) and (\ref{fmukk}) corresponding to the noisy measurements $(g_1^{\delta},g_2^{\delta})$ that satisfies condition (\ref{noise}). Then we have
\begin{equation}\label{noise_err1}
\left\|a_{\mu,\delta}^{(k)}-a_{\mu}^{(k)}\right\|+\left\|f_{\mu,\delta}^{(k)}-f_{\mu}^{(k)}\right\|\leq C\left(d^k+1\right)\frac{\delta}{\mu}.
\end{equation}
\end{lem}
\begin{proof}
By (\ref{amukk}), we know that 
\begin{eqnarray*}
a_{\mu,\delta}^{(k)}-a_{\mu}^{(k)}&=&\sum_{n=1}^{\infty}\frac{Q_n(T_1)S_n(T_2)}{(S_n(T_1)-\mu)(Q_n(T_2)+\mu)}\left(a_{\mu,\delta}^{(k-1)}-a_{\mu}^{(k-1)},\varphi_n\right)\varphi_n\\
&&+\sum_{n=1}^{\infty}\left[\frac{1}{S_n(T_1)-\mu}\left(g_1^{\delta}-g_1,\varphi_n\right)-\frac{Q_n(T_1)}{(S_n(T_1)-\mu)(Q_n(T_2)+\mu)}\left(g_2^{\delta}-g_2,\varphi_n\right)\right]\varphi_n\\
&=&\sum_{n=1}^{\infty}\left[\frac{Q_n(T_1)S_n(T_2)}{(S_n(T_1)-\mu)(Q_n(T_2)+\mu)}\right]^k\left(a_{\mu,\delta}^{(0)}-a_{\mu}^{(0)},\varphi_n\right)\varphi_n\\
&&+\sum_{n=1}^{\infty}N_{k-1}\left[\frac{1}{S_n(T_1)-\mu}\left(g_1^{\delta}-g_1,\varphi_n\right)-\frac{Q_n(T_1)}{(S_n(T_1)-\mu)(Q_n(T_2)+\mu)}\left(g_2^{\delta}-g_2,\varphi_n\right)\right]\varphi_n,
\end{eqnarray*}
where $N_k=\sum_{j=0}^{k}\left(\frac{Q_n(T_1)S_n(T_1)}{(S_n(T_1)-\mu)(Q_n(T_2)+\mu)}\right)^j$, and $N_k$ satisfies
\begin{eqnarray*}
N_k&=&\left[1-\left(\frac{Q_n(T_1)S_n(T_1)}{(S_n(T_1)-\mu)(Q_n(T_2)+\mu)}\right)^{k+1}\right]\Big/\left[1-\left(\frac{Q_n(T_1)S_n(T_1)}{(S_n(T_1)-\mu)(Q_n(T_2)+\mu)}\right)\right]\\
&\leq&1\Big/\left(1-\frac{Q_n(T_1)S_n(T_1)}{S_n(T_1)Q_n(T_2)}\right)\leq\frac{1}{1-d}.
\end{eqnarray*}
Since $f_{\mu,\delta}^{(0)}=f_{\mu}^{(0)}=f^{(0)}$,
\begin{eqnarray*}
a_{\mu,\delta}^{(0)}-a_{\mu}^{(0)}&=&\sum_{n=1}^{\infty}\frac{1}{S_n(T_1)-\mu}\left[(g_1^{\delta},\varphi_n)-Q_n(T_1)(f_{\mu,\delta}^{(0)},\varphi_n)\right]\varphi_n\\
&&-\sum_{n=1}^{\infty}\frac{1}{S_n(T_1)-\mu}\left[(g_1,\varphi_n)-Q_n(T_1)(f_{\mu}^{(0)},\varphi_n)\right]\varphi_n\\
&=&\sum_{n=1}^{\infty}\frac{1}{S_n(T_1)-\mu}\left(g_1^{\delta}-g_1,\varphi_n\right)\varphi_n,
\end{eqnarray*}
then it holds that
\begin{equation*}
\left\|a_{\mu,\delta}^{(0)}-a_{\mu}^{(0)}\right\|^2=\sum_{n=1}^{\infty}\frac{1}{(S_n(T_1)-\mu)^2}\left(g_1^{\delta}-g_1,\varphi_n\right)^2\leq C\frac{\delta^2}{\mu^2}.
\end{equation*}
By the fact that $N_k<\frac{1}{1-d}$, we have
\begin{eqnarray*}
\left\|a_{\mu,\delta}^{(k)}-a_{\mu}^{(k)}\right\|^2&\leq&2\sum_{n=1}^{\infty}\left[\frac{Q_n(T_1)S_n(T_2)}{(S_n(T_1)-\mu)(Q_n(T_2)+\mu)}\right]^{2k}\left(a_{\mu,\delta}^{(0)}-a_{\mu}^{(0)},\varphi_n\right)^2\\
&&+4N_{k-1}^2\sum_{n=1}^{\infty}\frac{1}{(S_n(T_1)-\mu)^2}\left[\left(g_1^{\delta}-g_1,\varphi_n\right)^2+\left(\frac{Q_n(T_1)}{Q_n(T_2)+\mu}\right)^2\left(g_2^{\delta}-g_2,\varphi_n\right)^2\right]\\
&\leq&C\left(d^{2k}+1\right)\frac{\delta^2}{\mu^2}.
\end{eqnarray*}
Similarly, we can also deduce the following estimate for $f_{\mu,\delta}^{(0)}-f_{\mu}^{(0)}$:
\begin{equation*}
\left\|f_{\mu,\delta}^{(k)}-f_{\mu}^{(k)}\right\|\leq C\left(d^k+1\right)\frac{\delta}{\mu}.
\end{equation*}
Finally, the proof is accomplished.
\end{proof}
Combining Lemmas \ref{SPlem2}, \ref{RSPlem1} and \ref{noise_err} together, we finally arrive at the main result of this section.
\begin{thm}\label{main1}
Let $a,f,f^{(0)}\in D((-\Delta)^p)$ for $p\in(0,1]$. $(a,f)$ are the exact solutions defined in (\ref{exactaf}). $(a_{\mu,\delta}^{(k)},f_{\mu,\delta}^{(k)})$ are defined in (\ref{amukk}) and (\ref{fmukk}) corresponding to the noisy measurements $(g_1^{\delta},g_2^{\delta})$ that satisfies condition (\ref{noise}). Then we have 
\begin{eqnarray}\label{main_res}
&&\left\|a_{\mu,\delta}^{(k)}-a\right\|+\left\|f_{\mu,\delta}^{(k)}-f\right\|\leq C\left(d^k+1\right)\left[\frac{\delta}{\mu}+\mu^p\left(\left\|a\right\|_{D((-\Delta)^p)}+\left\|f\right\|_{D((-\Delta)^p)}\right)\right]\\\nonumber
&&+d^k\left(\left\|a^{(0)}-a\right\|+\left\|f^{(0)}-f\right\|+\mu^p\left\|f^{(0)}\right\|_{D((-\Delta)^p)}\right),
\end{eqnarray}
where $a^{(0)}$ is given in (\ref{a0}).
\end{thm}

\begin{rem}\label{rem1}
It can be learned that if $k\rightarrow\infty$, the estimate (\ref{main_res}) reduces to the following form:
\begin{equation*}
\left\|a_{\mu,\delta}-a\right\|+\left\|f_{\mu,\delta}-f\right\|\leq C\left[\frac{\delta}{\mu}+\mu^p\left(\left\|a\right\|_{D((-\Delta)^p)}+\left\|f\right\|_{D((-\Delta)^p)}\right)\right],
\end{equation*}
where $\lim_{k\rightarrow\infty}a_{\mu,\delta}^{(k)}=a_{\mu,\delta}$ and $\lim_{k\rightarrow\infty}f_{\mu,\delta}^{(k)}=f_{\mu,\delta}$. We learn that the initial guess $f^{(0)}$ has slightly influence on the convergence result as $k$ increases. Under such situation, if we follow the \emph{a priori} choice rule of the regularization parameter, by taking $\mu=O(\delta^{\frac{1}{p+1}})$, the optimal convergence order in (\ref{main_res}) becomes $O(\delta^{\frac{p}{p+1}})$.
\end{rem}
\begin{rem}\label{rem2}
If $p=0$ in Theorem \ref{main1}, it can be viewed as a special case of $a,f,f^{(0)}\in L^2(\Omega)$. In such situation, we can learn from (\ref{main_res}) that the reconstructed solutions $(a_{\mu,\delta}^{(k)},f_{\mu,\delta}^{(k)})$ cannot converge to the exact solutions in $L^2$-sense. However, related convergence results can also be obtained in some weaker-norm senses such as Sobolev spaces $H^{m}(\Omega)$ with negative order $m$. One may refer \cite{Zhang+Zhou-2022} for detailed arguments and we do not consider the case here.
\end{rem}

\section{Spatially semidiscrete schemes and error estimates}\label{sec4}

In this section, we provide spatially semidiscrete schemes based on finite element method for solving both the direct problem and inverse problem. 

For $h\in(0,1)$, we denote $\{\mathcal{T}_h\}$ as a family of shape regular and quasi-uniform partitions of the region $\Omega$. The parameter $h$ is the maximum diameter of the elements in $\mathcal{T}_h$. Then we define the following finite elements space $X_h$ by
\begin{equation}\label{finitespace}
X_h:=\{\chi\in H_0^1(\Omega):\chi~\text{is a linear function over} ~\tau~\forall\tau\in\mathcal{T}_h\}.
\end{equation}
Next, following the standard description as in \cite{Jin+Lazarov+Zhou-2013}, we define two projection operators, i.e., $L^2$-projection $P_h: L^2(\Omega)\rightarrow X_h$ and Ritz projection $R_h: H_0^1(\Omega)\rightarrow X_h$ as follows:
\begin{eqnarray*}
&&(P_h\psi,\chi)=(\psi,\chi),\quad\forall\psi\in L^2(\Omega), \forall\chi\in X_h,\\
&&(R_h\psi,\chi)=(\psi,\chi),\quad\forall\psi\in H_0^1(\Omega), \forall\chi\in X_h.
\end{eqnarray*}
Then it can be learned from \cite{Thomee-2006} that both the projections satisfy the following properties:
\begin{equation}\label{appPhRh}
\left\|P_h\psi-\psi\right\|+\left\|R_h\psi-\psi\right\|\leq Ch^2\|\psi\|_{H^2(\Omega)},\quad\forall\psi\in D(-\Delta).
\end{equation}

\subsection{Error estimate for standard Galerkin method}\label{sec4.1} 

The semidiscrete scheme for solving direct problem (\ref{P}) based on standard Galerkin finite elements method can be expressed as follows:
\begin{eqnarray}\label{variaform}
&&(\partial_{t}^{\alpha}u_h(t),\chi)+(\nabla u_h(t),\nabla\chi)=(F,\chi),\quad\forall \chi\in X_h,\quad t\in(0,T],\\\nonumber
&&u_h(0)=P_h a,\quad\partial_t u_h(0)=P_hb.
\end{eqnarray}
To find the solution $u_h(t)\in X_h$ of (\ref{variaform}), we introduce the discrete Laplace operator $\Delta_h: X_h\rightarrow X_h$ by following definition:
\begin{equation}\label{defdislap}
-(\Delta_h\psi,\chi)=(\nabla\psi,\nabla\chi),\quad\forall\psi, \chi\in X_h.
\end{equation}
Then the semidiscrete scheme (\ref{variaform}) can be rewritten as 
\begin{align}\label{semiP}
\left\{
\begin{aligned}
&\partial_{t}^{\alpha}u_h(t)-\Delta_hu_h(t)=P_hF,& &t\in(0,T],\\
&u_h(0)=P_ha,~\partial_t u_h(0)=P_hb.& &
\end{aligned}
\right. 
\end{align}
We denote $\{\lambda_n^h\}_{n=1}^{K}$ and $\{\varphi_n^h\}_{n=1}^{K}$ as the eigenvalues and eigenfunctions of $-\Delta_h$, respectively. For general source term, the solution $u_h(t)$ to problem (\ref{semiP}) can be written as 
\begin{equation*}
u_h(t)=\mathbb{S}_h(t)P_ha+\mathbb{T}_h(t)P_hb+\int_0^t\mathbb{G}_h(t-\tau)P_hF(\tau)d\tau,
\end{equation*}
where the solution operators $\mathbb{S}_h(t),\mathbb{T}_h(t)$ and $\mathbb{G}_h(t)$ are defined by
\begin{eqnarray*}
&&\mathbb{S}_h(t)v=\sum_{n=1}^{K}S_n^h(t)(v,\varphi_n^h)\varphi_n^h:=\sum_{n=1}^{K}E_{\alpha,1}(-\lambda_n^h t^{\alpha})(v,\varphi_n^h)\varphi_n^h,\\
&&\mathbb{T}_h(t)v=\sum_{n=1}^{K}T_n^h(t)(v,\varphi_n^h)\varphi_n^h:=\sum_{n=1}^{K}tE_{\alpha,2}(-\lambda_n^h t^{\alpha})(v,\varphi_n^h)\varphi_n^h,\\
&&\mathbb{G}_h(t)v=\sum_{n=1}^{K}G_n^h(t)(v,\varphi_n^h)\varphi_n^h:=\sum_{n=1}^{K}t^{\alpha-1}E_{\alpha,\alpha}(-\lambda_n^h t^{\alpha})(v,\varphi_n^h)\varphi_n^h.
\end{eqnarray*}

When the source term is assumed to be time-independent, the semidiscrete scheme for solving problem (\ref{IP}) becomes
\begin{align}\label{semiIP}
\left\{
\begin{aligned}
&\partial_{t}^{\alpha}u_h(t)-\Delta_hu_h(t)=P_hf,& &t\in(0,T],\\
&u_h(0)=P_ha,~\partial_t u_h(0)=0,& &
\end{aligned}
\right. 
\end{align}
the solution $u_h(t)$ to problem (\ref{semiIP}) can be reduced to:
\begin{equation}\label{semi_solu}
u_h(t)=\mathbb{S}_h(t)P_ha+\mathbb{Q}_h(t)P_hf,
\end{equation}
where the solution operator $\mathbb{Q}_h(t)$ can be defined by 
\begin{equation*}
\mathbb{Q}_h(t)v=\sum_{n=1}^{K}Q_n^h(t)(v,\varphi_n^h)\varphi_n^h:=\sum_{n=1}^{K}t^{\alpha}E_{\alpha,\alpha+1}(-\lambda_n t^{\alpha})(v,\varphi_n^h)\varphi_n^h.
\end{equation*}

By the property of monotonically increasing of $\left\{\bar{\lambda}_n^h\right\}_{n=1}^{K}$, it can be learned that Lemmas \ref{SQ_bound}, \ref{SQQS_bound} also hold for operators $\mathbb{S}_h$ and $\mathbb{Q}_h$. Then the regularity of solution $u_h(t)$ defined by (\ref{semi_solu})  can be obtained by similar approach of Proposition \ref{direct2}. Thus we present the following result without proof.
\begin{prop}\label{semi_direct}
Let $u_h(t)$ be defined in (\ref{semi_solu}). If $a,~f\in D((-\Delta)^p)$ for $p\in[0,1]$, then $u(t)$ is the solution to problem (\ref{semiIP}), and satisfies the following estimate:
\begin{equation}\label{solu_regaf}
\|\partial_t^{(m)}u_h(t)\|_{D((-\Delta)^q)}\leq Ct^{-m-\alpha(q-p)}\left(\|a\|_{D((-\Delta)^p)}+t^{\alpha}\|f\|_{D((-\Delta)^p)}\right)
\end{equation}
for any integer $m\geq0$ and $q\in[p,p+1]$. 
\end{prop}

The error estimate between the solution $u(t)$ of problem (\ref{IP}) and the solution $u_h(t)$ of problem (\ref{semiIP}) are given in the following lemma.
\begin{lem}\label{error_u_uh}
Let $a,f\in D((-\Delta)^p)$, $u(t)$ and $u_h(t)$ be the solutions to problems (\ref{IP}) and (\ref{semiIP}), respectively. Then we have 
\begin{equation}\label{eq_err_u_uh}
\left\|u(t)-u_h(t)\right\|\leq Ch^2t^{-\alpha(1-p)}\left(\|a\|_{D((-\Delta)^p)}+t^{\alpha}\|f\|_{D((-\Delta)^p)}\right).
\end{equation}
\end{lem}
\noindent The proof is analogue to Theorem 3.2 in \cite{Jin+Lazarov+Zhou-2016}, thus we omit it here.

Now we propose the standard Galerkin method for solving the inverse problem based on the alternating approach provided in Section \ref{sec3}. For the given $f^{(0)}$, let $f_h^{(0)}=P_hf^{(0)}$, we consider solving the following regularized semidiscrete problem:
\begin{align}\label{SSP1}
\left\{
\begin{aligned}
&\partial_{t}^{\alpha}v_h^{(k)}(t)-\Delta_hv_h^{(k)}(t)=f_{\mu,h}^{(k)},& &t\in(0,T),\\
&v_h^{(k)}(0)=a_{\mu,h}^{(k)},~\partial_t v_h^{(k)}(0)=0,& &\\
&v_h^{(k)}(T_1)-\mu v_{h}^{(k)}(0)=P_hg_1.& &
\end{aligned}
\right.
\end{align}
The solution to (\ref{SSP1}) can be obtained as
\begin{eqnarray}\label{sak}
a_{\mu,h}^{(k)}&=&(\mathbb{S}_h(T_1)-\mu\mathbb{I})^{-1}\left[P_hg_1-\mathbb{Q}_h(T_1)f_{\mu,h}^{(k)}\right]\\\nonumber
&=&\sum_{n=1}^{K}\frac{1}{S_n^{h}(T_1)-\mu}(P_hg_1,\varphi_n^h)\varphi_n^h-\sum_{n=1}^{K}\frac{Q_n^h(T_1)}{S_n^{h}(T_1)-\mu}(f_{\mu,h}^{(k)},\varphi_n^h)\varphi_n^h.
\end{eqnarray}
With obtained $a_{\mu,h}^{(k)}$ in (\ref{sak}), we consider solving the following regularized semidiscrete problem:
\begin{align}\label{SSP2}
\left\{
\begin{aligned}
&\partial_{t}^{\alpha}w_h^{(k)}(t)-\Delta_hw_h^{(k)}(t)=f_{\mu,h}^{(k)},& &t\in(0,T),\\
&w_h^{(k)}(0)=a_{\mu,h}^{(k-1)},~\partial_t w_h^{(k)}(0)=0,& &\\
&w_h^{(k)}(T_2)+\mu f_{\mu,h}^{(k)}=P_hg_2.& &
\end{aligned}
\right.
\end{align}
The solution to (\ref{SSP2}) can be obtained as
\begin{eqnarray}\label{sfk}
f_{\mu,h}^{(k)}&=&(\mathbb{Q}_h(T_2)+\mu\mathbb{I})^{-1}\left[P_hg_2-\mathbb{S}_h(T_2)a_{\mu,h}^{(k-1)}\right]\\\nonumber
&=&\sum_{n=1}^{K}\frac{1}{Q_n^{h}(T_2)+\mu}(P_hg_2,\varphi_n^h)\varphi_n^h-\sum_{n=1}^{K}\frac{S_n^h(T_2)}{Q_n^{h}(T_2)+\mu}(a_{\mu,h}^{(k-1)},\varphi_n^h)\varphi_n^h.
\end{eqnarray}

Analogue to Lemma \ref{reg_regusolu},  we have the following regularity results for the obtained solutions $\left(a_{\mu,h}^{(k)},f_{\mu,h}^{(k)}\right)$.

\begin{lem}\label{reg_semiregusolu}
Let $a,f,f^{(0)}\in D((-\Delta)^p)$ for $p\in[0,1]$. Then for $q\in[p,p+1]$, we have
\begin{equation*}
\left\|a_{\mu,h}^{(k)}\right\|_{D((-\Delta)^q)}+\left\|f_{\mu,h}^{(k)}\right\|_{D((-\Delta)^q)}\leq C\mu^{p-q}\left(d^{k}+1\right)\left(\|a\|_{D((-\Delta)^p)}+\|f\|_{D((-\Delta)^p)}\right)+C\mu^{p-q}d^{k}\|f^{(0)}\|_{D((-\Delta)^p)}.
\end{equation*}
\end{lem}
Following Lemma 4.4 in \cite{Zhang+Zhou-2022}, we have the following result which is crucial for subsequent analysis.
\begin{lem}\label{estsource}
Let $v^{(k)}$ be the solution to problem (\ref{RSP1}). For $p\in[0,1]$, we have
\begin{equation}\label{res_estsource}
\left\|\int_0^t\mathbb{G}(t-\tau)\Delta_h(R_h-P_h)v^{k}(\tau)d\tau\right\|\leq C h^2t^{-\alpha(1-p)}\left(\|a_{\mu}^{(k)}\|_{D((-\Delta)^p)}+t^{\alpha}\|f_{\mu}^{(k)}\|_{D((-\Delta)^p)}\right).
\end{equation}
\end{lem}
\begin{proof}
We assume $y(t)$ is  the solution to the following problem
\begin{align}\label{SSSP2}
\left\{
\begin{aligned}
&\partial_{t}^{\alpha}y(t)-\Delta_hy(t)=P_hf_{\mu}^{(k)},& &t\in(0,T),\\
&y(0)=P_hv^{(k)}(0),~\partial_t y(0)=0.& &
\end{aligned}
\right.
\end{align}
Then we have
\begin{equation*}
y(t)-v^{(k)}(t)=\left(y(t)-P_hv^{(k)}(t)\right)+\left(P_hv^{(k)}(t)-v^{(k)}(t)\right)=:z_1(t)+z_2(t).
\end{equation*}
Following the regularity estimate (\ref{solu_regaf}) defined in Proposition \ref{direct2} and the approximate property  (\ref{appPhRh}) of the $L^2$-projection of $P_h$, we can obtain that 
\begin{equation*}
\|z_2(t)\|\leq Ch^2\|v^{(k)}(t)\|_{D(-\Delta)}\leq Ch^2t^{-\alpha(1-p)}\left(\|a_{\mu}^{(k)}\|_{D((-\Delta)^p)}+t^{\alpha}\|f_{\mu}^{(k)}\|_{D((-\Delta)^p)}\right).
\end{equation*}
Then we consider the estimate of $z_1(t)$. By using the property $P_h\Delta=\Delta_hR_h$ in finite element method theory, $z_1(t)$ satisfies the following problem
\begin{align*}
\left\{
\begin{aligned}
&\partial_{t}^{\alpha}z_1(t)-\Delta_hz_1(t)=\Delta_h\left(P_h-R_h\right)v^{(k)}(t),& &t\in(0,T),\\
&z_1(0)=0,~\partial_t z_1(0)=0.& &
\end{aligned}
\right.
\end{align*}
The solution $z_1(t)$ can be represented as 
\begin{equation*}
z_1(t)=\int_0^t\mathbb{G}(t-\tau)\Delta_h\left(P_h-R_h\right)v^{(k)}(\tau)d\tau.
\end{equation*}
Moreover, by Lemma \ref{error_u_uh}, it holds that
\begin{equation*}
\|y(t)-v^{(k)}(t)\|\leq Ch^2t^{-\alpha(1-p)}\left(\|a_{\mu}^{(k)}\|_{D((-\Delta)^p)}+t^{\alpha}\|f_{\mu}^{(k)}\|_{D((-\Delta)^p)}\right).
\end{equation*}
Then we know that
\begin{eqnarray*}
\|z_1(t)\|&=&\|(y(t)-v^{(k)}(t))-z_2(t)\|\leq\|y(t)-v^{(k)}(t)\|+\|z_2(t)\|\\
&\leq&Ch^2t^{-\alpha(1-p)}\left(\|a_{\mu}^{(k)}\|_{D((-\Delta)^p)}+t^{\alpha}\|f_{\mu}^{(k)}\|_{D((-\Delta)^p)}\right).
\end{eqnarray*}
Finally, the proof is accomplished.
\end{proof}

The following Lemma presents the error estimates between the solutions $(v^{(k)}(t),w^{(k)}(t))$ to problems (\ref{SP1}) and $(\ref{SP2})$ and semidiscrete solutions $(v_h^{(k)}(t),w_h^{(k)}(t))$ to problems (\ref{SSP1}) and $(\ref{SSP2})$. 
\begin{lem}\label{semi_error}
We assume that $a,f,f^{(0)}\in D((-\Delta)^p)$ for $p\in(0,1]$. Let $(v^{(k)}(t),w^{(k)}(t))$ be the solutions to problems (\ref{SP1}) and $(\ref{SP2})$ and $(v_h^{(k)}(t),w_h^{(k)}(t))$ be the solutions to problems (\ref{SSP1}) and $(\ref{SSP2})$, respectively. Then we have 
\begin{eqnarray}\label{semi_solu_err}
&&\left\|v^{(k)}(t)-v_h^{(k)}(t)\right\|+\left\|w^{(k)}(t)-w_h^{(k)}(t)\right\|\\\nonumber
&\leq& Ch^2t^{-\alpha(1-p)}\left(d^{k+1}+1\right)\left(\|a\|_{D((-\Delta)^p)}+\|f\|_{D((-\Delta)^p)}\right)+Ch^2t^{-\alpha(1-p)}d^k\|f^{(0)}\|_{D((-\Delta)^p)}.
\end{eqnarray}
\end{lem}
\begin{proof}
By (\ref{eq_err_u_uh}) in Lemma \ref{error_u_uh}, we know that
\begin{equation*}
\|v^{(k)}(t)-v_h^{(k)}(t)\|\leq Ch^2t^{-\alpha(1-p)}\left(\|a_{\mu}^{(k)}\|_{D((-\Delta)^p)}+\|f_{\mu}^{(k)}\|_{D((-\Delta)^p)}\right).
\end{equation*}
Combined with Lemma \ref{reg_regusolu} with $q=p$, we can obtain (\ref{semi_solu_err}).
\end{proof}

Next, we derive the error estimates between $(a_{\mu}^{(k)},f_{\mu}^{(k)})$ and $(a_{\mu,h}^{(k)},f_{\mu,h}^{(k)})$.

\begin{lem}\label{semilem1}
Let $a,f,f^{(0)}\in D((-\Delta)^p)$ for $p\in(0,1]$. $(a_{\mu}^{(k)},f_{\mu}^{(k)})$ are defined by (\ref{amuk}) and (\ref{fmuk}). $(a_{\mu,h}^{(k)},f_{\mu,h}^{(k)})$ are given by (\ref{sak}) and (\ref{sfk}). Then we have
\begin{equation}\label{semilem1res}
\left\|a_{\mu,h}^{(k)}-a_{\mu}^{(k)}\right\|+\left\|f_{\mu,h}^{(k)}-f_{\mu}^{(k)}\right\|\leq  \frac{C}{\mu}h^2\left(d^k+1\right)\left(\|a\|_{D((-\Delta)^p)}+\|f\|_{D((-\Delta)^p)}\right)+\frac{C}{\mu}h^2d^k\|f^{(0)}\|_{D((-\Delta)^p)}.
\end{equation}
\end{lem}
\begin{proof}
For simplicity, we only deduce the error estimate of $a_{\mu,h}^{(k)}-a_{\mu}^{(k)}$. The estimate for $f_{\mu,h}^{(k)}-f_{\mu}^{(k)}$ can be obtained by similar argument. 

We split the term $a_{\mu,h}^{(k)}-a_{\mu}^{(k)}$ into two parts, i.e.,
\begin{equation*}
a_{\mu,h}^{(k)}-a_{\mu}^{(k)}=\left(a_{\mu,h}^{(k)}-P_ha_{\mu}^{(k)}\right)+\left(P_ha_{\mu}^{(k)}-a_{\mu}^{(k)}\right).
\end{equation*} 
Let $q=1$ in Lemma \ref{reg_regusolu}, by using the property of Ritz projection in (\ref{appPhRh}), we can obtain that
\begin{eqnarray}\label{semilem11}
&&\|P_ha_{\mu}^{(k)}-a_{\mu}^{(k)}\|\leq Ch^2\|a_{\mu}^{(k)}\|_{D(-\Delta)}\\\nonumber
&\leq& Ch^2\mu^{p-1}\left(d^k+1\right)\left(\|a\|_{D((-\Delta)^p)}+\|f\|_{D((-\Delta)^p)}\right)+Ch^2\mu^{p-1}d^k\|f^{(0)}\|_{D((-\Delta)^p)}.
\end{eqnarray}
%and
%\begin{equation*}
%\left\|P_hf_{\mu}^{(k)}-f_{\mu}^{(k)}\right\|\leq Ch^2\|a_{\mu}^{(k)}\|_{D(-\Delta)}\leq Ch^2\mu^{p-1}(d^k+1)\left(\|a\|_{D((-\Delta)^p)}+\|f\|_{D((-\Delta)^p)}\right)+Ch^2\mu^{p-1}d^k\|f^{(0)}\|_{D((-\Delta)^p)}.
%\end{equation*}
Let $(v^{(k)}(t),w^{(k)}(t))$ be the solutions to problems (\ref{RSP1}) and (\ref{RSP2}), $(v_h^{(k)}(t),w_h^{(k)}(t))$ be the solutions to problem (\ref{SSP1}) and (\ref{SSP2}), respectively. We denote $\theta_1^{(k)}(t)=P_hv^{(k)}(t)-v_h^{(k)}(t)$ and $\theta_2^{(k)}(t)=P_hw^{(k)}(t)-w_h^{(k)}(t)$. Then it can be learned that $\theta_1^{(k)}(t)$ satisfies the following problem:
\begin{align}\label{theta1k}
\left\{
\begin{aligned}
&\partial_{t}^{\alpha}\theta_1^{(k)}(t)-\Delta_h\theta_1^{(k)}(t)=P_hf_{\mu}^{(k)}-f_{\mu,h}^{(k)}+\Delta_h\left(R_h-P_h\right)v^{(k)}(t),& &t\in(0,T),\\
&\theta_1^{(k)}(0)=P_ha_{\mu}^{(k)}-a_{\mu,h}^{(k)},~\partial_t \theta_1^{(k)}(0)=0,&&\\
&\theta_1^{(k)}(T_1)-\mu\theta_1^{(k)}(0)=0,&&
\end{aligned}
\right.
\end{align}
and $\theta_2^{(k)}(t)$ satisfies
\begin{align}\label{theta2k}
\left\{
\begin{aligned}
&\partial_{t}^{\alpha}\theta_2^{(k)}(t)-\Delta_h\theta_2^{(k)}(t)=P_hf_{\mu}^{(k)}-f_{\mu,h}^{(k)}+\Delta_h\left(R_h-P_h\right)w^{(k)}(t),& &t\in(0,T),\\
&\theta_2^{(k)}(0)=P_ha_{\mu}^{(k-1)}-a_{\mu,h}^{(k-1)},~\partial_t \theta_2^{(k)}(0)=0,&&\\
&\theta_2^{(k)}(T_2)+\mu\left(P_hf_{\mu}^{(k)}-f_{\mu,h}^{(k)}\right)=0.&&
\end{aligned}
\right.
\end{align}
Thus the solutions to (\ref{theta1k}) and (\ref{theta2k}) can be written as 
\begin{equation*}
\theta_1^{(k)}(t)=\mathbb{S}_h(t)\left(P_ha_{\mu}^{(k)}-a_{\mu,h}^{(k)}\right)+\mathbb{Q}_h(t)\left(P_hf_{\mu}^{(k)}-f_{\mu,h}^{(k)}\right)+\int_0^t\mathbb{G}_h(t-\tau)\Delta_h\left(R_h-P_h\right)v^{(k)}(\tau)d\tau,
\end{equation*}
and 
\begin{equation*}
\theta_2^{(k)}(t)=\mathbb{S}_h(t)\left(P_ha_{\mu}^{(k-1)}-a_{\mu,h}^{(k-1)}\right)+\mathbb{Q}_h(t)\left(P_hf_{\mu}^{(k)}-f_{\mu,h}^{(k)}\right)+\int_0^t\mathbb{G}_h(t-\tau)\Delta_h\left(R_h-P_h\right)w^{(k)}(\tau)d\tau.
\end{equation*}
By the third conditions in (\ref{theta1k}) and (\ref{theta2k}), we have
\begin{equation}\label{Pamuh}
P_ha_{\mu}^{(k)}-a_{\mu,h}^{(k)}=-(\mathbb{S}_h(T_1)-\mu\mathbb{I})^{-1}\left[\mathbb{Q}_h(T_1)\left(P_hf_{\mu}^{(k)}-f_{\mu,h}^{(k)}\right)+\int_0^{T_1}\mathbb{G}_h(T_1-\tau)\Delta_h\left(R_h-P_h\right)v^{(k)}(\tau)d\tau\right]
\end{equation}
and
\begin{equation*}
P_hf_{\mu}^{(k)}-f_{\mu,h}^{(k)}=-(\mathbb{Q}_h(T_2)+\mu\mathbb{I})^{-1}\left[\mathbb{S}_h(T_2)\left(P_ha_{\mu}^{(k-1)}-a_{\mu,h}^{(k-1)}\right)+\int_0^{T_2}\mathbb{G}_h(T_2-\tau)\Delta_h\left(R_h-P_h\right)w^{(k)}(\tau)d\tau\right],
\end{equation*}
which implies that
\begin{eqnarray*}
P_ha_{\mu}^{(k)}-a_{\mu,h}^{(k)}&=&(\mathbb{S}_h(T_1)-\mu\mathbb{I})^{-1}\mathbb{Q}_h(T_1)(\mathbb{Q}_h(T_2)+\mu\mathbb{I})^{-1}\mathbb{S}_h(T_2)\left(P_ha_{\mu}^{(k-1)}-a_{\mu,h}^{(k-1)}\right)\\
&&+(\mathbb{S}_h(T_1)-\mu\mathbb{I})^{-1}\mathbb{Q}_h(T_1)(\mathbb{Q}_h(T_2)+\mu\mathbb{I})^{-1}\int_0^{T_2}\mathbb{G}_h(T_2-\tau)\Delta_h\left(R_h-P_h\right)w^{(k)}(\tau)d\tau\\
&&-(\mathbb{S}_h(T_1)-\mu\mathbb{I})^{-1}\int_0^{T_1}\mathbb{G}_h(T_1-\tau)\Delta_h\left(R_h-P_h\right)v^{(k)}(\tau)d\tau\\
&=:&(\mathbb{S}_h(T_1)-\mu\mathbb{I})^{-1}\mathbb{Q}_h(T_1)(\mathbb{Q}_h(T_2)+\mu\mathbb{I})^{-1}\mathbb{S}_h(T_2)\left(P_ha_{\mu}^{(k-1)}-a_{\mu,h}^{(k-1)}\right)+J_1^h+J_2^h.
\end{eqnarray*}
Next, we divide it into three steps to accomplish the proof. \\
\emph{Step 1:} Estimation of $J_1^h$. For $\forall v\in L^2(\Omega)$, it can be learned from (\ref{SQ_S}) and (\ref{SQ_Q}) in Lemma \ref{SQ_bound} that
\begin{eqnarray*}
\left\|(\mathbb{S}_h(T_1)-\mu\mathbb{I})^{-1}\mathbb{Q}_h(T_1)(\mathbb{Q}_h(T_2)+\mu\mathbb{I})^{-1}v\right\|&=&\left[\sum_{n=1}^{K}\left(\frac{Q_n^h(T_1)}{(S_n^h(T_1)-\mu)(Q_n^h(T_2)+\mu)}\right)^2(v,\varphi_n^h)^2\right]^{\frac12}\\
&\leq&\frac{1}{\mu}\left[\sum_{n=1}^{K}\left(\frac{\overline{C_1}}{(\underline{C_2}+\mu\lambda_n^h)}\right)^2(v,\varphi_n^h)^2\right]^{\frac12}\leq\frac{C}{\mu}\|v\|.
\end{eqnarray*}
Then we apply the estimate (\ref{res_estsource}) in Lemma \ref{estsource}, it holds that
\begin{eqnarray*}
\|J_1^h\|&=&\left\|(\mathbb{S}_h(T_1)-\mu\mathbb{I})^{-1}\mathbb{Q}_h(T_1)(\mathbb{Q}_h(T_2)+\mu\mathbb{I})^{-1}\int_0^{T_2}\mathbb{G}_h(T_2-\tau)\Delta_h\left(R_h-P_h\right)w^{(k)}(\tau)d\tau\right\|\\
&\leq&\frac{C}{\mu}\left\|\int_0^{T_2}\mathbb{G}_h(T_2-\tau)\Delta_h\left(R_h-P_h\right)w^{(k)}(\tau)d\tau\right\|\\
&\leq&\frac{C}{\mu}h^2T_2^{-\alpha(1-p)}\left(\|a_{\mu}^{(k-1)}\|_{D((-\Delta)^p)}+T_2^{\alpha}\|f_{\mu}^{(k)}\|_{D((-\Delta)^p)}\right).
\end{eqnarray*}
Letting $q=p$ in Lemma \ref{reg_regusolu}, we have
\begin{equation*}
\|J_1^h\|\leq\frac{C}{\mu}h^2\left(d^k+1\right)\left(\|a\|_{D((-\Delta)^p)}+\|f\|_{D((-\Delta)^p)}\right)+\frac{C}{\mu}h^2d^k\|f^{(0)}\|_{D((-\Delta)^p)}.
\end{equation*}
\emph{Step 2:} Estimation of $J_2^h$. Similar to the procedure of deducing the estimate of $J_1^h$, it can be obtained that
\begin{eqnarray*}
\|J_2^h\|&=&\left\|(\mathbb{S}_h(T_1)-\mu\mathbb{I})^{-1}\int_0^{T_1}\mathbb{G}_h(T_1-\tau)\Delta_h\left(R_h-P_h\right)v^{(k)}(\tau)d\tau\right\|\\
&\leq&\frac{C}{\mu}\left\|\int_0^{T_1}\mathbb{G}_h(T_1-\tau)\Delta_h\left(R_h-P_h\right)v^{(k)}(\tau)d\tau\right\|\\
&\leq&\frac{C}{\mu}h^2T_1^{-\alpha(1-p)}\left(\|a_{\mu}^{(k)}\|_{D((-\Delta)^p)}+t^{\alpha}\|f_{\mu}^{(k)}\|_{D((-\Delta)^p)}\right).
\end{eqnarray*}
Letting $q=p$ in Lemma \ref{reg_regusolu}, we have
\begin{equation*}
\|J_2^h\|\leq\frac{C}{\mu}h^2\left(d^k+1\right)\left(\|a\|_{D((-\Delta)^p)}+\|f\|_{D((-\Delta)^p)}\right)+\frac{C}{\mu}h^2d^k\|f^{(0)}\|_{D((-\Delta)^p)}.
\end{equation*}
\emph{Step 3:} Estimation of $P_ha_{\mu}^{(k)}-a_{\mu,h}^{(k)}$. We deduce this estimation by induction. Since $f_{\mu,h}^{(0)}=P_hf_{\mu}^{(0)}=P_hf^{(0)}$, it follows from (\ref{Pamuh}) that
\begin{eqnarray*}
\left\|P_ha_{\mu}^{(0)}-a_{\mu,h}^{(0)}\right\|&=&\left\|(\mathbb{S}_h(T_1)-\mu\mathbb{I})^{-1}\int_0^{T_1}\mathbb{G}_h(T_1-\tau)\Delta_h\left(R_h-P_h\right)v^{(0)}(\tau)d\tau\right\|\\
&\leq&\frac{C}{\mu}h^2\left(\|a_{\mu}^{(0)}\|_{D((-\Delta)^p)}+\|f^{(0)}\|_{D((-\Delta)^p)}\right).
\end{eqnarray*}
Letting $q=p$ in Lemma \ref{reg_regusolu}, we obtain
\begin{equation*}
\left\|P_ha_{\mu}^{(0)}-a_{\mu,h}^{(0)}\right\|\leq\frac{C}{\mu}h^2\left(\|a\|_{D((-\Delta)^p)}+\|f\|_{D((-\Delta)^p)}+\|f^{(0)}\|_{D((-\Delta)^p)}\right).
\end{equation*}
Moreover, we have
\begin{eqnarray*}
\left\|P_ha_{\mu}^{(1)}-a_{\mu,h}^{(1)}\right\|&\leq& d\left\|P_ha_{\mu}^{(0)}-a_{\mu,h}^{(0)}\right\|+\|J_1^h\|+\|J_2^h\|\\
&\leq&\frac{C}{\mu}h^2\left(d+1\right)\left(\|a\|_{D((-\Delta)^p)}+\|f\|_{D((-\Delta)^p)}\right)+\frac{C}{\mu}h^2d\|f^{(0)}\|_{D((-\Delta)^p)}.
\end{eqnarray*}
Next we assume that for $j=2,\cdots,k$, we have the following estimate
\begin{equation}\label{Paj_induc}
\left\|P_ha_{\mu}^{(j)}-a_{\mu,h}^{(j)}\right\|\leq\frac{C}{\mu}h^2\left(d^j+1\right)\left(\|a\|_{D((-\Delta)^p)}+\|f\|_{D((-\Delta)^p)}\right)+\frac{C}{\mu}h^2d^j\|f^{(0)}\|_{D((-\Delta)^p)}.
\end{equation}
Then for $j=k+1$, by repeating above arguments, we have
\begin{eqnarray*}
\left\|P_ha_{\mu}^{(k+1)}-a_{\mu,h}^{(k+1)}\right\|&\leq&d\left\|P_ha_{\mu}^{(k)}-a_{\mu,h}^{(k)}\right\|+\|J_1^h\|+\|J_2^h\|\\
&\leq&\frac{C}{\mu}h^2\left(d^{k+1}+1\right)\left(\|a\|_{D((-\Delta)^p)}+\|f\|_{D((-\Delta)^p)}\right)+\frac{C}{\mu}h^2d^{k+1}\|f^{(0)}\|_{D((-\Delta)^p)}.
\end{eqnarray*}
Let $p=0$ in (\ref{semilem11}), we have
\begin{eqnarray*}
\left\|a_{\mu}^{(k)}-a_{\mu,h}^{(k)}\right\|&\leq&\left\|a_{\mu}^{(k)}-P_ha_{\mu}^{(k)}\right\|+\left\|P_ha_{\mu}^{(k)}-a_{\mu,h}^{(k)}\right\|\\
&\leq&\frac{C}{\mu}h^2\left(d^k+1\right)\left(\|a\|_{D((-\Delta)^p)}+\|f\|_{D((-\Delta)^p)}\right)+\frac{C}{\mu}h^2d^k\|f^{(0)}\|_{D((-\Delta)^p)}\\
&&+\frac{C}{\mu}h^2\left(d^k+1\right)\left(\|a\|+\|f\|\right)+\frac{C}{\mu}h^2d^k\|f^{(0)}\|\\
&\leq&\frac{C}{\mu}h^2\left(d^k+1\right)\left(\|a\|_{D((-\Delta)^p)}+\|f\|_{D((-\Delta)^p)}\right)+\frac{C}{\mu}h^2d^k\|f^{(0)}\|_{D((-\Delta)^p)}.
\end{eqnarray*}
Similarly, we can deduce the corresponding estimate of $\left\|f_{\mu}^{(k)}-f_{\mu,h}^{(k)}\right\|$. Therefore, the proof is accomplished.
\end{proof}
\begin{lem}\label{semi_noise_err}
Let $a,f,f^{(0)}\in D((-\Delta)^p)$ for $p\in(0,1]$. $(a_{\mu,h,\delta}^{(k)},f_{\mu,h,\delta}^{(k)})$ are defined in (\ref{sak}) and (\ref{sfk}) corresponding to the noisy measurements $(g_1^{\delta},g_2^{\delta})$ that satisfies condition (\ref{noise}). Then we have
\begin{equation}\label{semi_noise_err1}
\left\|a_{\mu,h,\delta}^{(k)}-a_{\mu,h}^{(k)}\right\|+\left\|f_{\mu,h,\delta}^{(k)}-f_{\mu,h}^{(k)}\right\|\leq C\left(d^k+1\right)\frac{\delta}{\mu}.
\end{equation}
\end{lem}
\begin{proof}
By (\ref{sak}) and (\ref{sfk}), we have
\begin{eqnarray*}
a_{\mu,h,\delta}^{(k)}-a_{\mu,h}^{(k)}&=&\sum_{n=1}^{K}\frac{Q_n^h(T_1)S_n^h(T_2)}{(S_n^h(T_1)-\mu)(Q_n^h(T_2)+\mu)}\left(a_{\mu,h,\delta}^{(k-1)}-a_{\mu,h}^{(k-1)},\varphi_n^h\right)\varphi_n^h\\
&&+\sum_{n=1}^{K}\left[\frac{1}{S_n^h(T_1)-\mu}(P_hg_1^{\delta}-P_hg_1,\varphi_n^h)-\frac{Q_n^h(T_1)}{(S_n^h(T_1)-\mu)(Q_n^h(T_2)+\mu)}(P_hg_2^{\delta}-P_hg_2,\varphi_n^h)\right]\varphi_n^h\\
&=&\sum_{n=1}^{K}\left[\frac{Q_n^h(T_1)S_n^h(T_2)}{(S_n^h(T_1)-\mu)(Q_n^h(T_2)+\mu)}\right]^k\left(a_{\mu,h,\delta}^{(0)}-a_{\mu,h}^{(0)},\varphi_n^h\right)\varphi_n^h\\
&&+N_{k-1}^h\sum_{n=1}^{K}\left[\frac{1}{S_n^h(T_1)-\mu}(P_hg_1^{\delta}-P_hg_1,\varphi_n^h)-\frac{Q_n^h(T_1)}{(S_n^h(T_1)-\mu)(Q_n^h(T_2)+\mu)}(P_hg_2^{\delta}-P_hg_2,\varphi_n^h)\right]\varphi_n^h,
\end{eqnarray*}
where $N_k^h=\sum_{j=0}^{k}\left(\frac{Q_n^h(T_1)S_n^h(T_1)}{(S_n^h(T_1)-\mu)(Q_n^h(T_2)+\mu)}\right)^j$, and $N_k^h$ satisfies
\begin{eqnarray*}
N_k^h&=&\left[1-\left(\frac{Q_n^h(T_1)S_n^h(T_1)}{(S_n^h(T_1)-\mu)(Q_n^h(T_2)+\mu)}\right)^{k+1}\right]\Big/\left[1-\left(\frac{Q_n^h(T_1)S_n^h(T_1)}{(S_n^h(T_1)-\mu)(Q_n^h(T_2)+\mu)}\right)\right]\\
&\leq&1\Big/\left(1-\frac{Q_n^h(T_1)S_n^h(T_1)}{S_n^h(T_1)Q_n^h(T_2)}\right)\leq\frac{1}{1-d}.
\end{eqnarray*}

Since $f_{\mu,h,\delta}^{(0)}=f_{\mu,h}^{0}=P_hf_{\mu}^{(0)}=P_hf^{(0)}$, we have
\begin{equation*}
\left\|a_{\mu,h,\delta}^{(0)}-a_{\mu,h}^{(0)}\right\|^2=\sum_{n=1}^{K}\frac{1}{(S_n^h(T_1)-\mu)^2}\left(P_h(g_1^{\delta}-g_1),\varphi_n^h\right)^2\leq C\frac{\delta^2}{\mu^2}.
\end{equation*}
By the fact that $N_{k}^h\leq\frac{1}{1-d}$, we have
\begin{eqnarray*}
\left\|a_{\mu,h,\delta}^{(k)}-a_{\mu,h}^{(k)}\right\|^2&\leq&2\sum_{n=1}^{K}\left[\frac{Q_n^h(T_1)S_n^h(T_2)}{(S_n^h(T_1)-\mu)(Q_n^h(T_2)+\mu)}\right]^{2k}\left(a_{\mu,h,\delta}^{(0)}-a_{\mu,h}^{(0)},\varphi_n^h\right)^2\\
&&+2\sum_{n=1}^K\left(\frac{N_{k-1}^h}{S_n^h(T_1)-\mu}\right)^2\left[\left(P_hg_1^{\delta}-P_hg_1,\varphi_n^h\right)^2+\left(\frac{Q_n^h(T_1)}{Q_n^h(T_2)+\mu}\right)^2\left(P_hg_2^{\delta}-P_hg_2,\varphi_n^h\right)^2\right]\\
&\leq&C\left(d^{2k}+1\right)\frac{\delta^2}{\mu^2}.
\end{eqnarray*}
Similarly, we can obtain the corresponding estimate of $f_{\mu,h,\delta}^{(0)}-f_{\mu,h}^{(0)}$ as follows:
\begin{equation*}
\left\|f_{\mu,h,\delta}^{(k)}-f_{\mu,h}^{(k)}\right\|\leq C\left(d^k+1\right)\frac{\delta}{\mu}.
\end{equation*}
Consequently, we accomplish the proof of (\ref{semi_noise_err1}).
\end{proof}

Combining Lemmas \ref{SPlem2}, \ref{RSPlem1}, \ref{semilem1} and \ref{semi_noise_err} together, we come to the following result.
\begin{thm}\label{main2}
Let $a,f,f^{(0)}\in D((-\Delta)^p)$ for $p\in(0,1]$.  $(a,f)$ are the exact solutions defined in (\ref{exactaf}). $(a_{\mu,h,\delta}^{(k)},f_{\mu,h,\delta}^{(k)})$ are defined in (\ref{sak}) and (\ref{sfk}) corresponding to the noisy measurements $(g_1^{\delta},g_2^{\delta})$ that satisfies condition (\ref{noise}). Then we have 
\begin{eqnarray}\label{main_res1}
&&\left\|a_{\mu,h,\delta}^{(k)}-a\right\|+\left\|f_{\mu,h,\delta}^{(k)}-f\right\|\\\nonumber
&\leq& C\left(d^k+1\right)\left[\frac{\delta}{\mu}+\left(\mu^p+\frac{h^2}{\mu}\right)\left(\left\|a\right\|_{D((-\Delta)^p)}+\left\|f\right\|_{D((-\Delta)^p)}\right)\right]\\\nonumber
&&+d^k\left[\left\|a-a^{(0)}\right\|+\left\|f-f^{(0)}\right\|+C\left(\mu^p+\frac{h^2}{\mu}\right)\left\|f^{(0)}\right\|_{D((-\Delta)^p)}\right],
\end{eqnarray}
where $a^{(0)}$ is given in (\ref{a0}).
\end{thm}

\begin{rem}\label{rem3}
It can be learned that if $k\rightarrow\infty$, the estimate (\ref{main_res1}) reduces to the following form:
\begin{equation*}
\left\|a_{\mu,h,\delta}-a\right\|+\left\|f_{\mu,h,\delta}-f\right\|\leq C\left[\frac{\delta}{\mu}+\left(\mu^p+\frac{h^2}{\mu}\right)\left(\left\|a\right\|_{D((-\Delta)^p)}+\left\|f\right\|_{D((-\Delta)^p)}\right)\right],
\end{equation*}
where $a_{\mu,h,\delta}=\lim_{k\rightarrow\infty}a_{\mu,h,\delta}^{(k)}$ and $f_{\mu,h,\delta}=\lim_{k\rightarrow\infty}f_{\mu,h,\delta}^{(k)}$. Under such situation, if we follow the \emph{a priori} choice rule of the regularization parameter, by taking $\mu=O(\delta^{\frac{1}{p+1}})$ and $h=O(\delta^{\frac12})$, the optimal convergence order in (\ref{main_res}) becomes $O(\delta^{\frac{p}{p+1}})$.
\end{rem}

\subsection{Error estimate for lumped mass method}\label{sec4.2} 
In this part, we define an approximate $L^2(\Omega)$-inner product in $X_h$ by
\begin{equation*}
\left(\omega,\chi\right)_h=\sum_{\tau\in\mathcal{T}_h}\frac{|\tau|}{d+1}\sum_{j=1}^{d+1}\omega(z_j^{\tau})\chi(z_j^{\tau}),
\end{equation*}
where $z_j^{\tau} (j=1,\cdots,d+1)$ denote the vertices of the $d$-complex $\tau\in\mathcal{T}_h$. Moreover, we introduce a projection operator $\overline{P}_h: L^2(\Omega)\rightarrow X_h$ as follows:
\begin{equation*}
\left(\bar{P}_hf, \chi\right)_h=(f,\chi),\quad\forall\chi\in X_h.
\end{equation*}
Then the lumped mass method for solving direct problem (\ref{P}) is to find solution $\bar{u}_h(t)$ that satisfies:
\begin{eqnarray}\label{variaformlum}
&&\left(\partial_{t}^{\alpha}\bar{u}_h(t),\chi\right)_h+(\nabla\bar{u}_h(t),\nabla\chi)=(F,\chi),\quad\forall \chi\in X_h,\quad t\in(0,T],\\\nonumber
&&\bar{u}_h(0)=\bar{P}_h a,\quad\partial_t \bar{u}_h(0)=\bar{P}_h b.
\end{eqnarray}
Therefore, we introduce another discrete Lapacian $\bar{\Delta}_h:X_h\rightarrow X_h$ by 
\begin{equation}\label{dis_Laplum}
-\left(\bar{\Delta}_h\psi,\chi\right)_h=(\nabla\psi,\nabla\chi).
\end{equation}
Then the semidiscrete scheme (\ref{variaformlum}) can be rewritten as 
\begin{align}\label{semiPlum}
\left\{
\begin{aligned}
&\partial_{t}^{\alpha}\bar{u}_h(t)-\bar{\Delta}_h\bar{u}_h(t)=\bar{P}_hF,& &t\in(0,T],\\
&\bar{u}_h(0)=\bar{P}_ha,~\partial_t \bar{u}_h(0)=\bar{P}_hb.& &
\end{aligned}
\right. 
\end{align}
Similar to standard Galerkin method, the solution to problem (\ref{semiPlum}) can be represented as
\begin{equation*}
\bar{u}_h(t):=\bar{\mathbb{S}}_h(t)\bar{P}_ha+\bar{\mathbb{T}}_h(t)\bar{P}_hb+\int_0^t\bar{\mathbb{G}}_h(t-\tau)\bar{P}_hF(\tau)d\tau,
\end{equation*}
where the solution operators $\bar{\mathbb{S}}_h(t), \bar{\mathbb{T}}_h(t)$ and $\bar{\mathbb{G}}_h(t)$ are defined by
\begin{eqnarray*}
&&\bar{\mathbb{S}}_h(t)v=\sum_{n=1}^{K}\bar{S}_n^h(t)(v,\bar{\varphi}_n^h)\bar{\varphi}_n^h:=\sum_{n=1}^{K}E_{\alpha,1}(-\bar{\lambda}_n^h t^{\alpha})(v,\bar{\varphi}_n^h)\bar{\varphi}_n^h,\\
&&\bar{\mathbb{T}}_h(t)v=\sum_{n=1}^{K}\bar{T}_n^h(t)(v,\bar{\varphi}_n^h)\bar{\varphi}_n^h:=\sum_{n=1}^{K}tE_{\alpha,2}(-\bar{\lambda}_n^h t^{\alpha})(v,\bar{\varphi}_n^h)\bar{\varphi}_n^h,\\
&&\bar{\mathbb{G}}_h(t)v=\sum_{n=1}^{K}\bar{G}_n^h(t)(v,\bar{\varphi}_n^h)\bar{\varphi}_n^h:=\sum_{n=1}^{K}t^{\alpha-1}E_{\alpha,\alpha}(-\bar{\lambda}_n^h t^{\alpha})(v,\bar{\varphi}_n^h)\bar{\varphi}_n^h,
\end{eqnarray*}
and $\left\{\bar{\lambda}_n^h,\bar{\varphi}_n^h\right\}_{n=1}^{K}$ denotes the eigensystem of $-\bar{\Delta}_h$. In the case that $F(t)\equiv f$ is time-independent and $b\equiv0$, the solution can be reduced to
\begin{equation}\label{semilumsolu}
\bar{u}_h(t):=\bar{\mathbb{S}}_h(t)\bar{P}_ha+\bar{\mathbb{Q}}_h(t)\bar{P}_hf,
\end{equation}
where the solution operator $\bar{\mathbb{Q}}_h(t)$ is defined by
\begin{equation*}
\bar{\mathbb{Q}}_h(t)v=\sum_{n=1}^{K}\bar{Q}_n^h(t)(v,\bar{\varphi}_n^h)\bar{\varphi}_n^h:=\sum_{n=1}^{K}t^{\alpha}E_{\alpha,\alpha+1}(-\bar{\lambda}_n^h t^{\alpha})(v,\bar{\varphi}_n^h)\bar{\varphi}_n^h.
\end{equation*}
By the property of monotonically increasing of $\left\{\bar{\lambda}_n^h\right\}_{n=1}^{K}$, it can be learned that Lemmas \ref{SQ_bound}, \ref{SQQS_bound} also hold for operators $\bar{\mathbb{S}}_h$ and $\bar{\mathbb{Q}}_h$. 

To measure the error between the inner products $(\cdot,\cdot)_h$ and $(\cdot,\cdot)$, we define an operator $\mathcal{E}_h:X_h\rightarrow X_h$ by
\begin{equation}\label{inner_error}
\left(\nabla\mathcal{E}_h\psi,\nabla\chi\right)=\epsilon_h(\chi,\psi):=(\chi,\psi)_h-(\chi,\psi),\quad\forall\chi,\psi\in X_h.
\end{equation}
It can be learned from \cite{Chatzipantelidis+Lazarov+Thomée-2012} that the following estimate holds for $\mathcal{E}_h$.
\begin{lem}\cite{Chatzipantelidis+Lazarov+Thomée-2012}
Let $\bar{\Delta}_h$ and $\mathcal{E}_h$ be defined by $(\ref{dis_Laplum})$ and $(\ref{inner_error})$, respectively. Then we have
\begin{equation*}
\left\|\nabla\mathcal{E}_h\chi\right\|+h\left\|\bar{\Delta}_h\mathcal{E}_h\chi\right\|\leq Ch^{p+1}\left\|\nabla^p\chi\right\|,\quad\forall\chi\in X_h,\quad p=0,1.
\end{equation*}
\end{lem}

However, in order to obtain an optimal error estimate analogue to standard Galerkin method, we always assume in this section that the region $\Omega$ can be divided into symmetric triangulations. Under such assumption, a stronger condition \cite{Chatzipantelidis+Lazarov+Thomée-2012}
\begin{equation}\label{h2bound}
\|\mathcal{E}_h\psi\|\leq ch^2\|\psi\|,\quad\psi\in X_h,
\end{equation}
is satisfied which is crucial for subsequent analysis. 

Next we apply the lumped mass method for solving the inverse problem based on the alternating approach provided in Section \ref{sec3}. For the given $f^{(0)}$, let $\bar{f}_h^{(0)}=\bar{P}_hf^{(0)}$, we consider solving the following regularized semidiscrete problem:
\begin{align}\label{SSP1lum}
\left\{
\begin{aligned}
&\partial_{t}^{\alpha}\bar{v}_h^{(k)}(t)-\bar{\Delta}_h\bar{v}_h^{(k)}(t)=\bar{f}_{\mu,h}^{(k)},& &t\in(0,T),\\
&\bar{v}_h^{(k)}(0)=\bar{a}_{\mu,h}^{(k)},~\partial_t \bar{v}_h^{(k)}(0)=0,& &\\
&\bar{v}_h^{(k)}(T_1)-\mu\bar{v}_h^{(k)}(0)=\bar{P}_hg_1.& &
\end{aligned}
\right.
\end{align}
The solution to (\ref{SSP1lum}) can be expressed by
\begin{eqnarray}\label{saklum}
\bar{a}_{\mu,h}^{(k)}&=&(\bar{\mathbb{S}}_h(T_1)-\mu\mathbb{I})^{-1}\left[\bar{P}_hg_1-\bar{\mathbb{Q}}_h(T_1)\bar{f}_{\mu,h}^{(k)}\right]\\\nonumber
&=&\sum_{n=1}^{K}\frac{1}{\bar{S}_n^{h}(T_1)-\mu}(\bar{P}_hg_1,\bar{\varphi}_n^h)\bar{\varphi}_n^h-\sum_{n=1}^{K}\frac{\bar{Q}_n^h(T_1)}{\bar{S}_n^{h}(T_1)-\mu}(\bar{f}_{\mu,h}^{(k)},\bar{\varphi}_n^h)\bar{\varphi}_n^h.
\end{eqnarray}
With obtained $\bar{a}_{\mu,h}^{(k)}$ in (\ref{saklum}), we consider solving the following regularized semidiscrete problem:
\begin{align}\label{SSP2lum}
\left\{
\begin{aligned}
&\partial_{t}^{\alpha}\bar{w}_h^{(k)}(t)-\bar{\Delta}_h\bar{w}_h^{(k)}(t)=\bar{f}_{\mu,h}^{(k)},& &t\in(0,T),\\
&\bar{w}_h^{(k)}(0)=\bar{a}_{\mu,h}^{(k-1)},~\partial_t \bar{w}_h^{(k)}(0)=0,& &\\
&\bar{w}_h^{(k)}(T_2)+\mu\bar{f}_{\mu,h}^{(k)}=\bar{P}_hg_2.
\end{aligned}
\right.
\end{align}
The solution to (\ref{SSP2lum}) can be obtained as
\begin{eqnarray}\label{sfklum}
\bar{f}_{\mu,h}^{(k)}&=&(\bar{\mathbb{Q}}_h(T_2)+\mu\mathbb{I})^{-1}\left[\bar{P}_hg_2-\bar{\mathbb{S}}_h(T_2)\bar{a}_{\mu,h}^{(k-1)}\right]\\\nonumber
&=&\sum_{n=1}^{K}\frac{1}{\bar{Q}_n^{h}(T_2)+\mu}(\bar{P}_hg_2,\bar{\varphi}_n^h)\bar{\varphi}_n^h-\sum_{n=1}^{K}\frac{\bar{S}_n^h(T_2)}{\bar{Q}_n^{h}(T_2)+\mu}(\bar{a}_{\mu,h}^{(k-1)},\bar{\varphi}_n^h)\bar{\varphi}_n^h.
\end{eqnarray}
Similar to Lemma \ref{estsource}, we present the following key  estimate for the main results.
\begin{lem}\label{estsourcelum}
Let $u_h(t)$ be the solution to problem (\ref{semiIP}). For $a,f\in D((-\Delta)^p)$  with $p\in(0,1]$, we have
\begin{equation}\label{res_estsourcelum}
\left\|\int_0^t\bar{\mathbb{G}}_h(t-\tau)\bar{\Delta}_h\mathcal{E}_h\partial_{\tau}^{\alpha}u_h(\tau)d\tau\right\|\leq C h^2t^{-\alpha(1-p)}\left(\|a\|_{D((-\Delta)^p)}+t^{\alpha}\|f\|_{D((-\Delta)^p)}\right).
\end{equation}
\end{lem}
\begin{proof}
Assume that $y(t)$ is the solution to following problem
\begin{align*}
\left\{
\begin{aligned}
&\partial_{t}^{\alpha}y(t)-\bar{\Delta}_hy(t)=P_hf,& &t\in(0,T),\\
&y(0)=P_ha,\quad\partial_t y(0)=0.
\end{aligned}
\right.
\end{align*}
Let $\theta(t)=u_h(t)-y(t)$, then it satisfies 
\begin{align}
\left\{
\begin{aligned}
&\partial_{t}^{\alpha}\theta(t)-\bar{\Delta}_h\theta(t)=\bar{\Delta}_h\mathcal{E}_h\partial_t^{\alpha}u_h(t),& &t\in(0,T),\\
&\theta(0)=0,\quad\partial_t \theta(0)=0.
\end{aligned}
\right.
\end{align}
The solution can be expressed as
\begin{equation*}
\theta(t)=\int_0^t\bar{\mathbb{G}}_h(t-\tau)\bar{\Delta}_h\mathcal{E}_h\partial_{\tau}^{\alpha}u_h(\tau)d\tau.
\end{equation*}
By the definition of Caputo derivative, we can rewrite $\theta(t)$ as
\begin{eqnarray*}
\theta(t)&=&\int_0^t\bar{\mathbb{G}}_h(t-\tau)\bar{\Delta}_h\mathcal{E}_h\partial_{\tau}^{\alpha}u_h(\tau)d\tau=\int_0^t\bar{\mathbb{T}}_h(t-\tau)\bar{\Delta}_h\mathcal{E}_h\partial_{\tau}^{2}u_h(\tau)d\tau\\
&=&\left[\int_0^{\frac t2}+\int_{\frac t2}^t\right]\bar{\mathbb{T}}_h(t-\tau)\bar{\Delta}_h\mathcal{E}_h\partial_{\tau}^{2}u_h(\tau)d\tau=:J_1(t)+J_2(t).
\end{eqnarray*}
We consider $J_1(t)$ firstly. By integration by parts, it holds that
\begin{eqnarray*}
J_1(t)&=&\int_0^{\frac t2}\bar{\mathbb{T}}_h(t-\tau)\bar{\Delta}_h\mathcal{E}_h\partial_{\tau}^{2}u_h(\tau)d\tau\\
&=&\left[\bar{\mathbb{T}}_h(t-\tau)\bar{\Delta}_h\mathcal{E}_h\partial_{\tau}u_h(\tau)\right]_0^{\frac t2}+\int_0^{\frac t2}\bar{\mathbb{S}}_h(t-\tau)\bar{\Delta}_h\mathcal{E}_h\partial_{\tau}u_h(\tau)d\tau\\
&=&\bar{\mathbb{T}}_h(\frac t2)\bar{\Delta}_h\mathcal{E}_h\partial_tu_h(\frac t2)+\int_0^{\frac t2}\bar{\mathbb{S}}_h(t-\tau)\bar{\Delta}_h\mathcal{E}_h\partial_{\tau}u_h(\tau)d\tau.
\end{eqnarray*}
By the fact that $\bar{\mathbb{T}}_h(t)\bar{\Delta}_hv=\partial_t^{\alpha}\bar{\mathbb{T}}_h(t)v$, it follows from (\ref{h2bound}), Lemma \ref{lem2.1} and Lemma \ref{semi_direct} with  $m=1,q=0$ that 
\begin{eqnarray*}
\left\|\bar{\mathbb{T}}_h(\frac t2)\bar{\Delta}_h\mathcal{E}_h\partial_tu_h(\frac t2)\right\|&\leq& Ct^{1-\alpha}\left\|\mathcal{E}_h\partial_tu_h(\frac t2)\right\|\leq Ct^{1-\alpha}h^2\left\|\partial_tu_h(\frac t2)\right\|\\
&\leq&Ch^2t^{\alpha(p-1)}\left(\left\|a\right\|_{D\left((-\Delta)^p\right)}+t^{\alpha}\left\|f\right\|_{D\left((-\Delta)^p\right)}\right).
\end{eqnarray*}
Moreover, since $\bar{\mathbb{S}}_h(t)\bar{\Delta}_hv=\partial_t^{\alpha}\bar{\mathbb{S}}_h(t)v$, let $m=1,q=0$ in Lemma \ref{semi_direct}, we have
\begin{eqnarray*}
&&\left\|\int_0^{\frac t2}\bar{\mathbb{S}}_h(t-\tau)\bar{\Delta}_h\mathcal{E}_h\partial_{\tau}u_h(\tau)d\tau\right\|\leq\int_0^{\frac t2}\left\|\bar{\mathbb{S}}_h(t-\tau)\bar{\Delta}_h\mathcal{E}_h\partial_{\tau}u_h(\tau)\right\|d\tau\\
&\leq&Ch^2\int_0^{\frac t2}(t-\tau)^{-\alpha}\tau^{\alpha p-1}d\tau\left(\left\|a\right\|_{D\left((-\Delta)^p\right)}+t^{\alpha}\left\|f\right\|_{D\left((-\Delta)^p\right)}\right)\\
&\leq&Ch^2t^{\alpha(p-1)}\left(\left\|a\right\|_{D\left((-\Delta)^p\right)}+t^{\alpha}\left\|f\right\|_{D\left((-\Delta)^p\right)}\right).
\end{eqnarray*}
Thus we have 
\begin{equation*}
\left\|J_1(t)\right\|\leq Ch^2t^{\alpha(p-1)}\left(\left\|a\right\|_{D\left((-\Delta)^p\right)}+t^{\alpha}\left\|f\right\|_{D\left((-\Delta)^p\right)}\right).
\end{equation*}
Next, we consider $J_2(t)$. We also use the fact that $\bar{\mathbb{T}}_h(t)\bar{\Delta}_hv=\partial_t^{\alpha}\bar{\mathbb{T}}_h(t)v$, let $m=2,q=0$ in Lemma \ref{semi_direct}, we have
\begin{eqnarray*}
\left\|J_2(t)\right\|&\leq&\int_{\frac t2}^t\left\|\bar{\mathbb{T}}_h(t-\tau)\bar{\Delta}_h\mathcal{E}_h\partial_{\tau}^2u_h(\tau)\right\|d\tau\\
&\leq&C\int_{\frac t2}^t(t-\tau)^{1-\alpha}\left\|\mathcal{E}_h\partial_{\tau}^2u_h(\tau)\right\|d\tau\\
&\leq&Ch^2\int_{\frac t2}^t(t-\tau)^{1-\alpha}\tau^{\alpha p-2}d\tau\left(\left\|a\right\|_{D\left((-\Delta)^p\right)}+t^{\alpha}\left\|f\right\|_{D\left((-\Delta)^p\right)}\right)\\
&\leq&Ch^2t^{\alpha(p-1)}\left(\left\|a\right\|_{D\left((-\Delta)^p\right)}+t^{\alpha}\left\|f\right\|_{D\left((-\Delta)^p\right)}\right).
\end{eqnarray*}
Combined with the estimate of $J_1(t)$, we obtain $(\ref{res_estsourcelum})$.
\end{proof}
\begin{lem}\label{semilumlem1}
Let $a,f,f^{(0)}\in D((-\Delta)^p)$ for $p\in(0,1]$. $(\bar{a}_{\mu,h}^{(k)},\bar{f}_{\mu,h}^{(k)})$ are defined by (\ref{saklum}) and (\ref{sfklum}). $(a_{\mu,h}^{(k)},f_{\mu,h}^{(k)})$ are given by (\ref{sak}) and (\ref{sfk}). Then we have
\begin{equation}\label{semilumlem1res}
\left\|\bar{a}_{\mu,h}^{(k)}-a_{\mu,h}^{(k)}\right\|+\left\|\bar{f}_{\mu,h}^{(k)}-f_{\mu,h}^{(k)}\right\|\leq Ch^2\mu^{-1}\left(d^k+1\right)\left(\left\|a\right\|_{D((-\Delta)^p)}+\left\|f\right\|_{D((-\Delta)^p)}\right)+Ch^2\mu^{-1}d^k\left\|f^{(0)}\right\|_{D((-\Delta)^p)}.
\end{equation}
\end{lem}
\begin{proof}
For simplicity, we only consider the estimate of $\bar{a}_{\mu,h}^{(k)}-a_{\mu,h}^{(k)}$. Let $\bar{v}_{h}^{(k)}(t)$ and $v_{h}^{(k)}(t)$ be solutions to problems (\ref{SSP1lum}) and (\ref{SSP1}), respectively. Denote $\theta(t):=\bar{v}_{h}^{(k)}(t)-v_{h}^{(k)}(t)$, it follows from the variational forms (\ref{variaform}) and (\ref{variaformlum}) that
\begin{align}\label{lem411}
\left\{
\begin{aligned}
&\partial_{t}^{\alpha}\theta(t)-\bar{\Delta}_h\theta(t)=\left(\bar{f}_{\mu,h}^{(k)}-f_{\mu,h}^{(k)}\right)+\bar{\Delta}_h\mathcal{E}_h\partial_t^{\alpha}v_h^{(k)}(t)-\bar{\Delta}_h\mathcal{E}_hf_{\mu,h}^{(k)},& &t\in(0,T),\\
&\theta(0)=\bar{a}_{\mu,h}^{(k)}-a_{\mu,h}^{(k)},\quad\partial_t \theta(0)=0,\\&\theta(T_1)-\mu\left(\bar{a}_{\mu,h}^{(k)}-a_{\mu,h}^{(k)}\right)=\bar{P}_hg_1-P_hg_1.
\end{aligned}
\right.
\end{align}
By the third condition in problem (\ref{lem411}), we have
\begin{eqnarray*}
\bar{a}_{\mu,h}^{(k)}-a_{\mu,h}^{(k)}&=&\left(\bar{\mathbb{S}}_h(T_1)-\mu\mathbb{I}\right)^{-1}\left(\bar{P}_hg_1-P_hg_1\right)\\
&&-\left(\bar{\mathbb{S}}_h(T_1)-\mu\mathbb{I}\right)^{-1}\bar{\mathbb{Q}}_h(T_1)\left[\left(\bar{f}_{\mu,h}^{(k)}-f_{\mu,h}^{(k)}\right)-\bar{\Delta}_h\mathcal{E}_hf_{\mu,h}^{(k)}\right]\\
&&-\left(\bar{\mathbb{S}}_h(T_1)-\mu\mathbb{I}\right)^{-1}\int_0^{T_1}\bar{\mathbb{G}}_h(T_1-\tau)\bar{\Delta}_h\mathcal{E}\partial_{\tau}^{\alpha}v_h^{(k)}(\tau)d\tau.
\end{eqnarray*}
Similarly, we have
\begin{eqnarray*}
\bar{f}_{\mu,h}^{(k)}-f_{\mu,h}^{(k)}&=&\left(\bar{\mathbb{Q}}_h(T_2)+\mu\mathbb{I}\right)^{-1}\left(\bar{P}_hg_2-P_hg_2\right)
\\
&&-\left(\bar{\mathbb{Q}}_h(T_2)+\mu\mathbb{I}\right)^{-1}\left[\bar{\mathbb{S}}_h(T_2)\left(\bar{a}_{\mu,h}^{(k-1)}-a_{\mu,h}^{(k-1)}\right)-\bar{\mathbb{Q}}_h(T_2)\bar{\Delta}_h\mathcal{E}_hf_{\mu,h}^{(k)}\right]\\
&&-\left(\bar{\mathbb{Q}}_h(T_2)+\mu\mathbb{I}\right)^{-1}\int_0^{T_2}\bar{\mathbb{G}}_h(T_2-\tau)\bar{\Delta}_h\mathcal{E}_h\partial_{\tau}^{\alpha}w_h^{(k)}(\tau)d\tau.
\end{eqnarray*}
Substituting $\bar{f}_{\mu,h}^{(k)}-f_{\mu,h}^{(k)}$ into $\bar{a}_{\mu,h}^{(k)}-a_{\mu,h}^{(k)}$, we have
\begin{eqnarray*}
\bar{a}_{\mu,h}^{(k)}-a_{\mu,h}^{(k)}&=&\left(\bar{\mathbb{S}}_h(T_1)-\mu\mathbb{I}\right)^{-1}\left(\bar{P}_hg_1-P_hg_1\right)\\
&&-\left(\bar{\mathbb{S}}_h(T_1)-\mu\mathbb{I}\right)^{-1}\bar{\mathbb{Q}}_h(T_1)\left(\bar{\mathbb{Q}}_h(T_2)+\mu\mathbb{I}\right)^{-1}\left(\bar{P}_hg_2-P_hg_2\right)\\
&&+\left(\bar{\mathbb{S}}_h(T_1)-\mu\mathbb{I}\right)^{-1}\bar{\mathbb{Q}}_h(T_1)\left(\bar{\mathbb{Q}}_h(T_2)+\mu\mathbb{I}\right)^{-1}\bar{\mathbb{S}}_h(T_2)\left(\bar{a}_{\mu,h}^{(k-1)}-a_{\mu,h}^{(k-1)}\right)\\
&&-\left(\bar{\mathbb{S}}_h(T_1)-\mu\mathbb{I}\right)^{-1}\bar{\mathbb{Q}}_h(T_1)\left(\bar{\mathbb{Q}}_h(T_2)+\mu\mathbb{I}\right)^{-1}\bar{\mathbb{Q}}_h(T_2)\bar{\Delta}_h\mathcal{E}_hf_{\mu,h}^{(k)}\\
&&+\left(\bar{\mathbb{S}}_h(T_1)-\mu\mathbb{I}\right)^{-1}\bar{\mathbb{Q}}_h(T_1)\bar{\Delta}_h\mathcal{E}_hf_{\mu,h}^{(k)}\\
&&+\left(\bar{\mathbb{S}}_h(T_1)-\mu\mathbb{I}\right)^{-1}\bar{\mathbb{Q}}_h(T_1)\left(\bar{\mathbb{Q}}_h(T_2)+\mu\mathbb{I}\right)^{-1}\int_0^{T_2}\bar{\mathbb{G}}_h(T_2-\tau)\bar{\Delta}_h\mathcal{E}_h\partial_{\tau}^{\alpha}w_h^{(k)}(\tau)d\tau\\
&&-\left(\bar{\mathbb{S}}_h(T_1)-\mu\mathbb{I}\right)^{-1}\int_0^{T_1}\bar{\mathbb{G}}_h(T_1-\tau)\bar{\Delta}_h\mathcal{E}\partial_{\tau}^{\alpha}v_h^{(k)}(\tau)d\tau\\
&=&\left(\bar{\mathbb{S}}_h(T_1)-\mu\mathbb{I}\right)^{-1}\bar{\mathbb{Q}}_h(T_1)\left(\bar{\mathbb{Q}}_h(T_2)+\mu\mathbb{I}\right)^{-1}\bar{\mathbb{S}}_h(T_1)\left(\bar{a}_{\mu,h}^{(k-1)}-a_{\mu,h}^{(k-1)}\right)+\sum_{j=1}^{3}J_i^h,
\end{eqnarray*}
where
\begin{eqnarray*}
J_1^h&:=&\left(\bar{\mathbb{S}}_h(T_1)-\mu\mathbb{I}\right)^{-1}\left[\left(\bar{P}_hg_1-P_hg_1\right)-\bar{\mathbb{Q}}_h(T_1)\left(\bar{\mathbb{Q}}_h(T_2)+\mu\mathbb{I}\right)^{-1}\left(\bar{P}_hg_2-P_hg_2\right)\right],\\
J_2^h&:=&\left(\bar{\mathbb{S}}_h(T_1)-\mu\mathbb{I}\right)^{-1}\bar{\mathbb{Q}}_h(T_1)\left[\bar{\Delta}_h\mathcal{E}_hf_{\mu,h}^{(k)}-\left(\bar{\mathbb{Q}}_h(T_2)+\mu\mathbb{I}\right)^{-1}\bar{\mathbb{Q}}_h(T_2)\bar{\Delta}_h\mathcal{E}_hf_{\mu,h}^{(k)}\right],\\
J_3^h&:=&\left(\bar{\mathbb{S}}_h(T_1)-\mu\mathbb{I}\right)^{-1}\bar{\mathbb{Q}}_h(T_1)\left(\bar{\mathbb{Q}}_h(T_2)+\mu\mathbb{I}\right)^{-1}\int_0^{T_2}\bar{\mathbb{G}}_h(T_2-\tau)\bar{\Delta}_h\mathcal{E}_h\partial_{\tau}^{\alpha}w_h^{(k)}(\tau)d\tau\\
&&-\left(\bar{\mathbb{S}}_h(T_1)-\mu\mathbb{I}\right)^{-1}\int_0^{T_1}\bar{\mathbb{G}}_h(T_1-\tau)\bar{\Delta}_h\mathcal{E}\partial_{\tau}^{\alpha}v_h^{(k)}(\tau)d\tau.
\end{eqnarray*}
Let $k=0$, we have
\begin{eqnarray*}
\left\|\bar{a}_{\mu,h}^{(0)}-a_{\mu,h}^{(0)}\right\|&\leq&\left\|\left(\bar{\mathbb{S}}_h(T_1)-\mu\mathbb{I}\right)^{-1}\left(\bar{P}_hg_1-P_hg_1\right)\right\|\\
&&+\left\|\left(\bar{\mathbb{S}}_h(T_1)-\mu\mathbb{I}\right)^{-1}\bar{\mathbb{Q}}_h(T_1)\left(\bar{f}_{\mu,h}^{(0)}-f_{\mu,h}^{(0)}\right)\right\|\\
&&+\left\|\left(\bar{\mathbb{S}}_h(T_1)-\mu\mathbb{I}\right)^{-1}\bar{\mathbb{Q}}_h(T_1)\bar{\Delta}_h\mathcal{E}_hf_{\mu,h}^{(0)}\right\|\\
&&+\left\|\left(\bar{\mathbb{S}}_h(T_1)-\mu\mathbb{I}\right)^{-1}\int_0^{T_1}\bar{\mathbb{G}}_h(T_1-\tau)\bar{\Delta}_h\mathcal{E}\partial_{\tau}^{\alpha}v_h^{(0)}(\tau)d\tau\right\|.
\end{eqnarray*}
By the fact that $\left(\bar{P}_hg_1,\chi\right)_h=\left(P_hg_1,\chi\right)=\left(g_1,\chi\right)$, we have
\begin{eqnarray*}
0&=&\left(\bar{P}_hg_1,\chi\right)_h-\left(P_hg_1,\chi\right)\\
&=&\left(\bar{P}_hg_1,\chi\right)_h-\left(\bar{P}_hg_1,\chi\right)+\left(\bar{P}_hg_1,\chi\right)-\left(P_hg_1,\chi\right)\\
&=&\left(\nabla\mathcal{E}_h\bar{P}_hg_1,\nabla\chi\right)+\left(\bar{P}_hg_1-P_hg_1,\chi\right)\\
&=&-\left(\Delta_h\mathcal{E}_h\bar{P}_hg_1,\chi\right)+\left(\bar{P}_hg_1-P_hg_1,\chi\right),
\end{eqnarray*}
which implies that $\bar{P}_{h}g_1-P_hg_1=\Delta_h\mathcal{E}_h\bar{P}_hg_1$. Then we have the following results
\begin{eqnarray*}
&&\left\|\left(\bar{\mathbb{S}}_h(T_1)-\mu\mathbb{I}\right)^{-1}\Delta_h\mathcal{E}_h\bar{P}_hg_1\right\|
\leq C\mu^{-1}\left\|\Delta_h\mathcal{E}_h\bar{P}_hg_1\right\|\\
&\leq&Ch^2\mu^{-1}\left(\left\|a\right\|+\left\|f\right\|\right)\leq Ch^2\mu^{-1}\left(\left\|a\right\|_{D((-\Delta)^p)}+\left\|f\right\|_{D((-\Delta)^p)}\right).
\end{eqnarray*}
It can be learned from above deduction that $\bar{f}_{\mu,h}^{(0)}-f_{\mu,h}^{(0)}=\Delta_h\mathcal{E}_h\bar{P}_hf^{(0)}$, thus we have
\begin{eqnarray*}
&&\left\|\left(\bar{\mathbb{S}}_h(T_1)-\mu\mathbb{I}\right)^{-1}\bar{\mathbb{Q}}_h(T_1)\left(\bar{f}_{\mu,h}^{(0)}-f_{\mu,h}^{(0)}\right)\right\|+\left\|\left(\bar{\mathbb{S}}_h(T_1)-\mu\mathbb{I}\right)^{-1}\bar{\mathbb{Q}}_h(T_1)\bar{\Delta}_h\mathcal{E}_hf_{\mu,h}^{(0)}\right\|\\
&\leq& C\mu^{-1}\left(\left\|\mathcal{E}_h\bar{f}^{(0)}\right\|+\left\|\mathcal{E}_hf_{\mu,h}^{(0)}\right\|\right)\leq Ch^2\mu^{-1}\left\|f^{(0)}\right\|\leq Ch^2\mu^{-1}\left\|f^{(0)}\right\|_{D((-\Delta)^p)}
\end{eqnarray*}
and
\begin{eqnarray*}
&&\left\|\left(\bar{\mathbb{S}}_h(T_1)-\mu\mathbb{I}\right)^{-1}\int_0^{T_1}\bar{\mathbb{G}}_h(T_1-\tau)\bar{\Delta}_h\mathcal{E}\partial_{\tau}^{\alpha}v_h^{(0)}(\tau)d\tau\right\|\\
&\leq&Ch^2\mu^{-1}\left(\left\|a_{\mu,h}^{(0)}\right\|_{D((-\Delta)^p)}+\left\|f_{\mu,h}^{(0)}\right\|_{D((-\Delta)^p)}\right).
\end{eqnarray*}
By Lemme \ref{reg_semiregusolu}, we have
\begin{equation*}
\left\|\bar{a}_{\mu,h}^{(0)}-a_{\mu,h}^{(0)}\right\|\leq Ch^2\mu^{-1}\left(\left\|a\right\|_{D((-\Delta)^p)}+\left\|f\right\|_{D((-\Delta)^p)}+\left\|f^{(0)}\right\|_{D((-\Delta)^p)}\right).
\end{equation*}
Moreover, it holds that
\begin{eqnarray*}
\left\|J_1^h\right\|&\leq& C\mu^{-1}\left(\left\|\bar{P}_hg_1-P_hg_1\right\|+\left\|\bar{P}_hg_1-P_hg_1\right\|\right)\\
&\leq&C\mu^{-1}\left(\left\|\Delta_h\mathcal{E}_h\bar{P}_hg_1\right\|+\left\|\Delta_h\mathcal{E}_h\bar{P}_hg_2\right\|\right)\\
&\leq&Ch^2\mu^{-1}\left(\left\|a\right\|_{D((-\Delta)^p)}+\left\|f\right\|_{D((-\Delta)^p)}\right),
\end{eqnarray*}
\begin{eqnarray*}
\left\|J_2^h\right\|&=&\left\|\mu\left(\bar{\mathbb{S}}_h(T_1)-\mu\mathbb{I}\right)^{-1}\bar{\mathbb{Q}}_h(T_1)\left(\bar{\mathbb{Q}}_h(T_2)+\mu\mathbb{I}\right)^{-1}\bar{\Delta}_h\mathcal{E}_hf_{\mu,h}^{(k)}\right\|\leq C\mu^{-1}\left\|\bar{\Delta}_h\mathcal{E}_hf_{\mu,h}^{(k)}\right\|\\
&\leq&\frac{C}{\mu}h^2\left(d^k+1\right)\left(\left\|a\right\|_{D((-\Delta)^p)}+\left\|f\right\|_{D((-\Delta)^p)}\right)+\frac{C}{\mu}h^2d^k\left\|f^{(0)}\right\|_{D((-\Delta)^p)}
\end{eqnarray*}
and
\begin{eqnarray*}
\left\|J_3^h\right\|&\leq&C\mu^{-1}\left(\left\|a_{\mu,h}^{(k)}\right\|_{D((-\Delta)^p)}+\left\|a_{\mu,h}^{(k-1)}\right\|_{D((-\Delta)^p)}+\left\|f_{\mu,h}^{(k)}\right\|_{D((-\Delta)^p)}\right)\\
&\leq&\frac{C}{\mu}h^2\left(d^k+1\right)\left(\left\|a\right\|_{D((-\Delta)^p)}+\left\|f\right\|_{D((-\Delta)^p)}\right)+\frac{C}{\mu}h^2d^k\left\|f^{(0)}\right\|_{D((-\Delta)^p)}.
\end{eqnarray*}
The estimate of $J_3^h$ follows from Lemma \ref{reg_semiregusolu} with $p=q$. Therefore, we have
\begin{equation*}
\left\|\bar{a}_{\mu,h}^{(1)}-a_{\mu,h}^{(1)}\right\|\leq Cd\left\|\bar{a}_{\mu,h}^{(0)}-a_{\mu,h}^{(0)}\right\|+\frac{C}{\mu}h^2\left(d+1\right)\left(\left\|a\right\|_{D((-\Delta)^p)}+\left\|f\right\|_{D((-\Delta)^p)}\right)+\frac{C}{\mu}h^2d\left\|f^{(0)}\right\|_{D((-\Delta)^p)}.
\end{equation*}
Applying the method of induction as in Lemma \ref{semilem1}, we can obtain that 
\begin{equation*}
\left\|\bar{a}_{\mu,h}^{(k)}-a_{\mu,h}^{(k)}\right\|\leq \frac{C}{\mu}h^2\left(d^k+1\right)\left(\left\|a\right\|_{D((-\Delta)^p)}+\left\|f\right\|_{D((-\Delta)^p)}\right)+\frac{C}{\mu}h^2d^k\left\|f^{(0)}\right\|_{D((-\Delta)^p)}.
\end{equation*}
Similarly, we can obtain the corresponding estimate of $\left\|\bar{f}_{\mu,h}^{(k)}-f_{\mu,h}^{(k)}\right\|$. Finally, we accomplish the proof.
\end{proof}
Repeating the arguments in Lemma \ref{semi_noise_err}, we have the following estimate.
\begin{lem}\label{semi_noise_errlum}
Let $a,f,f^{(0)}\in D((-\Delta)^p)$ for $p\in(0,1]$. $(\bar{a}_{\mu,h,\delta}^{(k)},\bar{f}_{\mu,h,\delta}^{(k)})$ are defined in (\ref{saklum}) and (\ref{sfklum}) corresponding to the noisy measurements $(g_1^{\delta},g_2^{\delta})$ that satisfies condition (\ref{noise}). Then we have
\begin{equation}\label{semi_noise_err1lum}
\left\|\bar{a}_{\mu,h,\delta}^{(k)}-\bar{a}_{\mu,h}^{(k)}\right\|+\left\|\bar{f}_{\mu,h,\delta}^{(k)}-\bar{f}_{\mu,h}^{(k)}\right\|\leq C\left(d^k+1\right)\frac{\delta}{\mu}.
\end{equation}
\end{lem}
Combining Lemma \ref{SPlem2}, Lemma \ref{RSPlem1}, Lemma \ref{semilem1}, Lemma \ref{semilumlem1} and Lemma \ref{semi_noise_errlum}, we obtain the main convergence result for the lump mass method.
\begin{thm}\label{main3}
Let $a,f,f^{(0)}\in D((-\Delta)^p)$ for $p\in(0,1]$.  $(a,f)$ are the exact solutions defined in (\ref{exactaf}). $(\bar{a}_{\mu,h,\delta}^{(k)},\bar{f}_{\mu,h,\delta}^{(k)})$ are defined in (\ref{saklum}) and (\ref{sfklum}) corresponding to the noisy measurements $(g_1^{\delta},g_2^{\delta})$ that satisfies condition (\ref{noise}). Then we have 
\begin{eqnarray}\label{main_res2}
&&\left\|a_{\mu,h,\delta}^{(k)}-a\right\|+\left\|f_{\mu,h,\delta}^{(k)}-f\right\|\\\nonumber
&\leq& C\left(d^k+1\right)\left[\frac{\delta}{\mu}+\left(\mu^p+\frac{h^2}{\mu}\right)\left(\left\|a\right\|_{D((-\Delta)^p)}+\left\|f\right\|_{D((-\Delta)^p)}\right)\right]\\\nonumber
&&+d^k\left[\left\|a-a^{(0)}\right\|+\left\|f-f^{(0)}\right\|+C\left(\mu^p+\frac{h^2}{\mu}\right)\left\|f^{(0)}\right\|_{D((-\Delta)^p)}\right],
\end{eqnarray}
where $a^{(0)}$ is given in (\ref{a0}). 
\end{thm}

\begin{rem}\label{rem4}
The estimate (\ref{main_res2}) is derived under the assumption that the region $\Omega$ is divided into symmetric triangulations. For nonsymmetric partitions, the condition (\ref{h2bound}) does not hold. Thus the convergence order in (\ref{main_res2}) is invalid. Other techniques should be taken for related analysis which is beyond the content of this paper.
\end{rem}

\section{Numerical implementations}\label{sec5}

In this section, we present numerical examples to demonstrate the effectiveness and applicability of the proposed method for reconstructing the initial value and source term simultaneously. The numerical examples aim at illustrating the behavior of the proposed method under various conditions. Specifically, we investigate the influence  of noise level, regularization parameter and discretized parameter on the solution accuracy and computational efficiency. Through these numerical experiments, we  provide insights into the practical utility and robustness of the proposed method. 
The additional measurements are generated by adding random perturbation, i.e.
\begin{eqnarray*}
&&g_1^{\delta}=u(T_1)+\epsilon\max_{x\in\overline{\Omega}}u(T_1)\left(2\cdot \text{randn}(\text{size}(g_1))-1\right),\\
&&g_2^{\delta}=u(T_2)+\epsilon\max_{x\in\overline{\Omega}}u(T_2)\left(2\cdot \text{randn}(\text{size}(g_2))-1\right),
\end{eqnarray*}
where $\delta=\epsilon\cdot\max\left(\|g_1\|,\|g_2\|\right)$ denotes the noise level. 
In order to measure the accuracy of the reconstructed solutions, we use the following relative errors:
\begin{equation*}
re_a:=\frac{\|a_{\mu,h,\delta}^{(k)}-a\|}{\|a\|},\quad re_f:=\frac{\|f_{\mu,h,\delta}^{(k)}-f\|}{\|f\|},
\end{equation*}
where $(a,f)$ denote the exact solutions. 

\subsection{One-dimensional case}\label{sec5.1}
We take $\Omega=(0,1)$, and divide it into $K$ subintervals equally with mesh size $h=\frac{1}{K}$. It can be learned from \cite{Jin+Lazarov+Zhou-2013} that the eigensystem of $-\bar{\Delta}_h$ based on lumped mass method is
\begin{equation*}
\bar{\lambda}_n^h=\frac{4}{h^2}\sin^2\frac{n\pi}{2K},\quad\bar{\varphi}_n^h=\sqrt{2}\sin(n\pi x_n), \quad n=1,\cdots,K-1,
\end{equation*}
while the eigensystem of $-\Delta_h$ based on standard Galerkin method is
\begin{equation*}
\lambda_n^h=\frac{\bar{\lambda}_n^h}{1-\frac{h^2}{6}\bar{\lambda}_n^h},\quad\varphi_n^h=\sqrt{2}\sin(n\pi x_n), \quad n=1,\cdots,K-1.
\end{equation*}
Thus the reconstructed solutions of the initial value and source term can be computed through (\ref{sak}), (\ref{sfk}), (\ref{saklum}) and (\ref{sfklum}) directly. 

\noindent\textbf{Example 1.} Let $\alpha=1.6$, $d=0.4$ and $T_2=1.5$. We take  $a(x)=\sin(2\pi x)$ and $f(x)=\sin(\pi x)$.\\
With the given $a(x)$ and $f(x)$, the exact solution to direct problem (\ref{IP}) can be given by
\begin{equation}\label{ex1eq1}
u(x,t)=E_{\alpha,1}(-4\pi t^{\alpha})\sin(2\pi x)+t^{\alpha}E_{\alpha,\alpha+1}(-\pi t^{\alpha})\sin(\pi x).
\end{equation}
Thus the exact measurement can be obtained by letting $t=T_1, T_2$ in (\ref{ex1eq1}). 

\begin{figure}[htbp]
  \centering
  \subfigure{
    \includegraphics[width=0.35\textwidth,height=1.5in]{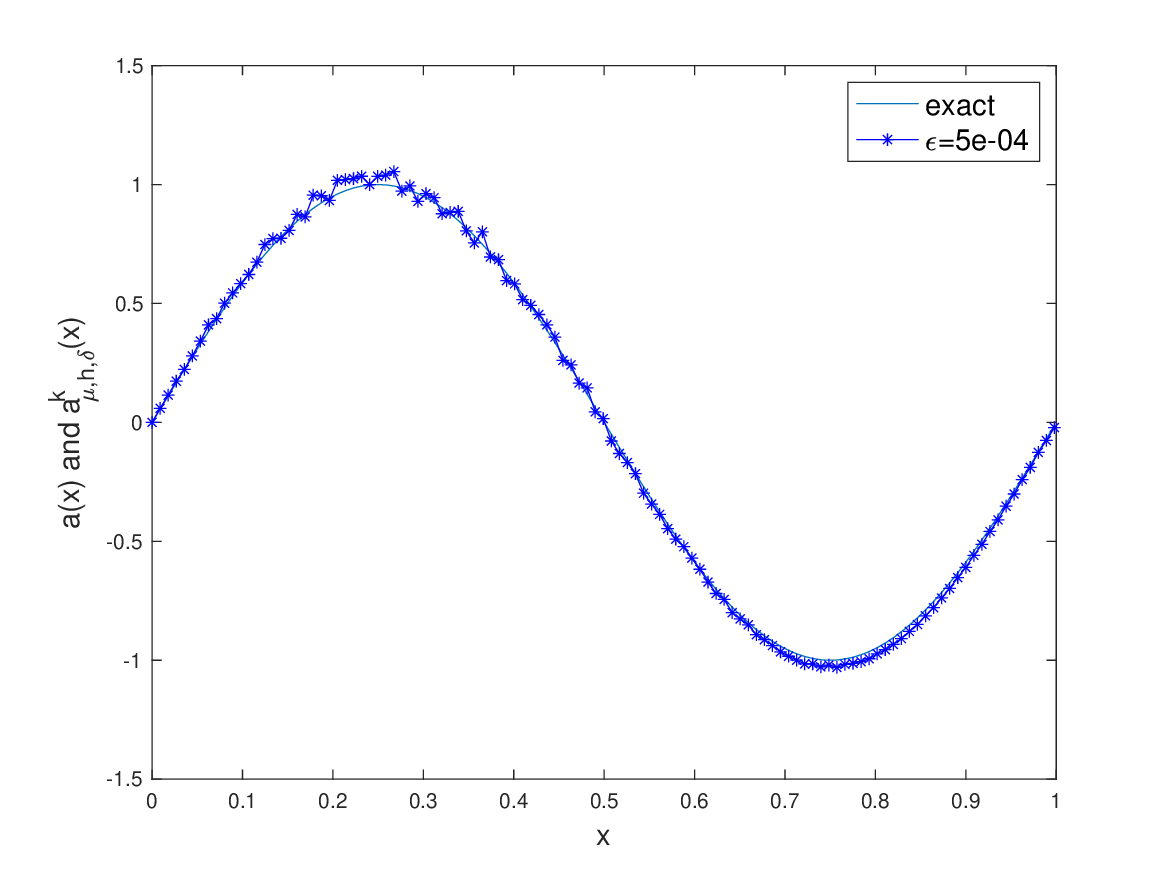}}
  \subfigure{
    \includegraphics[width=0.35\textwidth,height=1.5in]{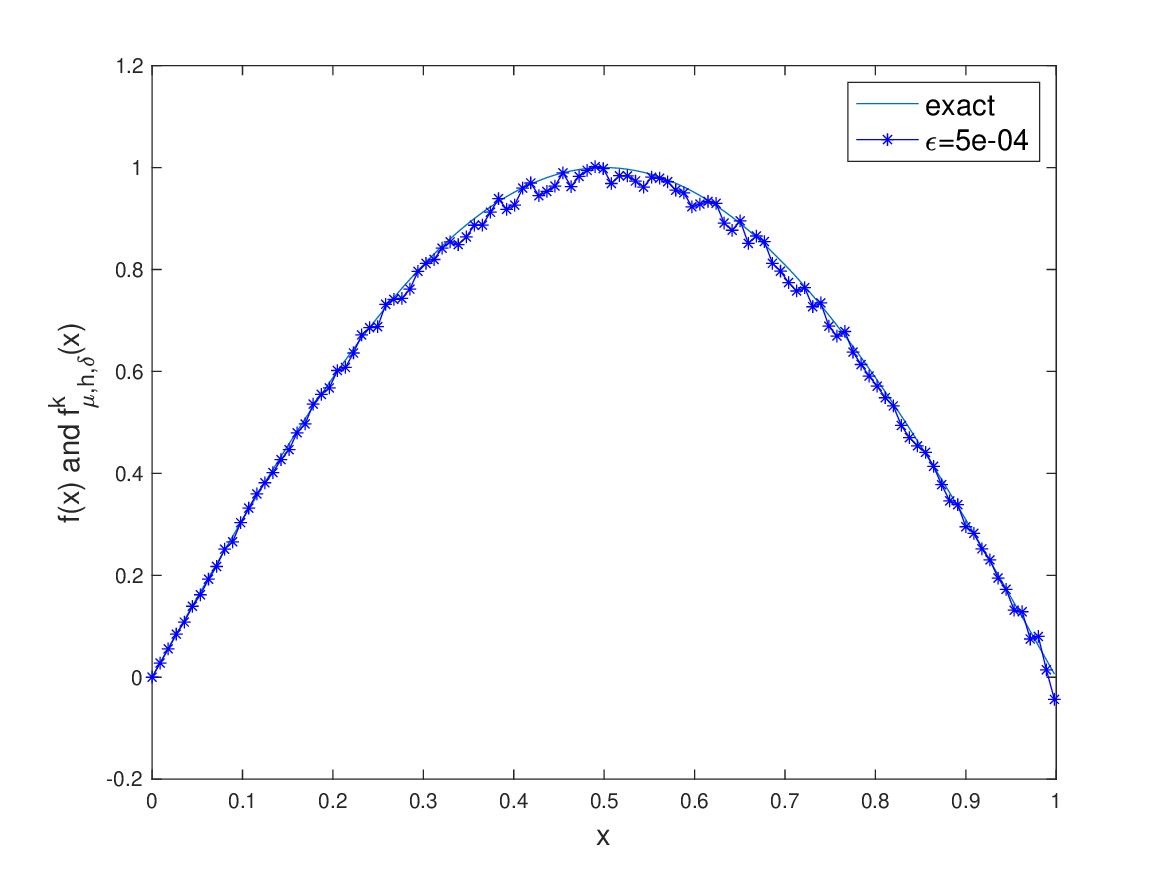}}\\
      \subfigure{
    \includegraphics[width=0.35\textwidth,height=1.5in]{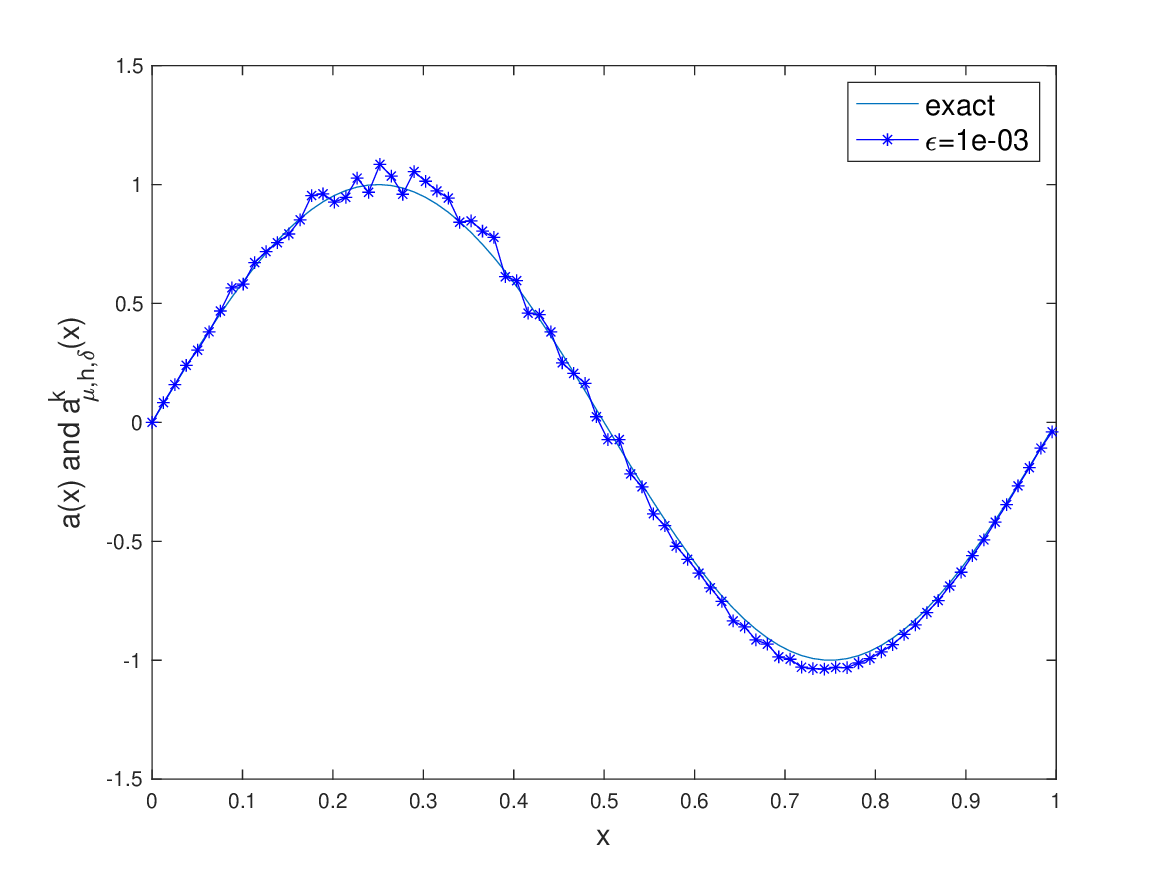}}
      \subfigure{
    \includegraphics[width=0.35\textwidth,height=1.5in]{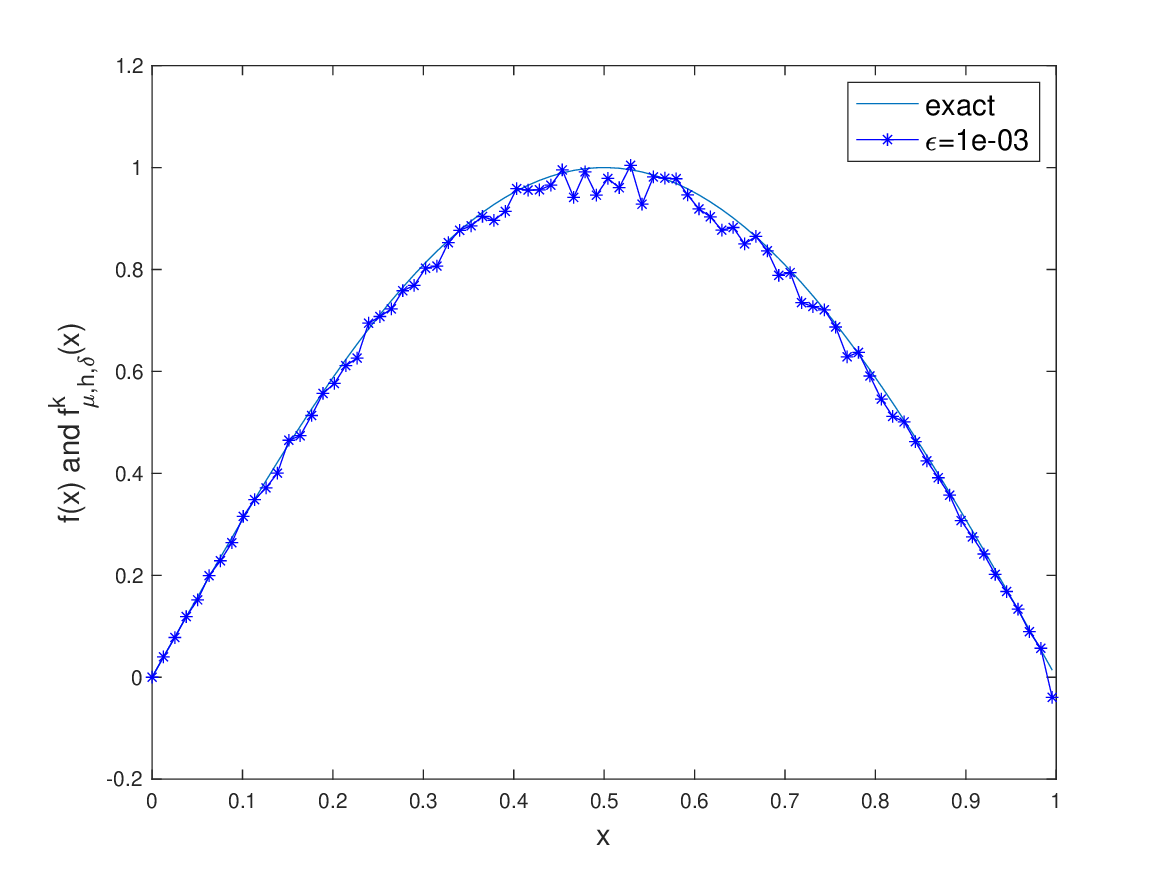}}\\
  \subfigure{
    \includegraphics[width=0.35\textwidth,height=1.5in]{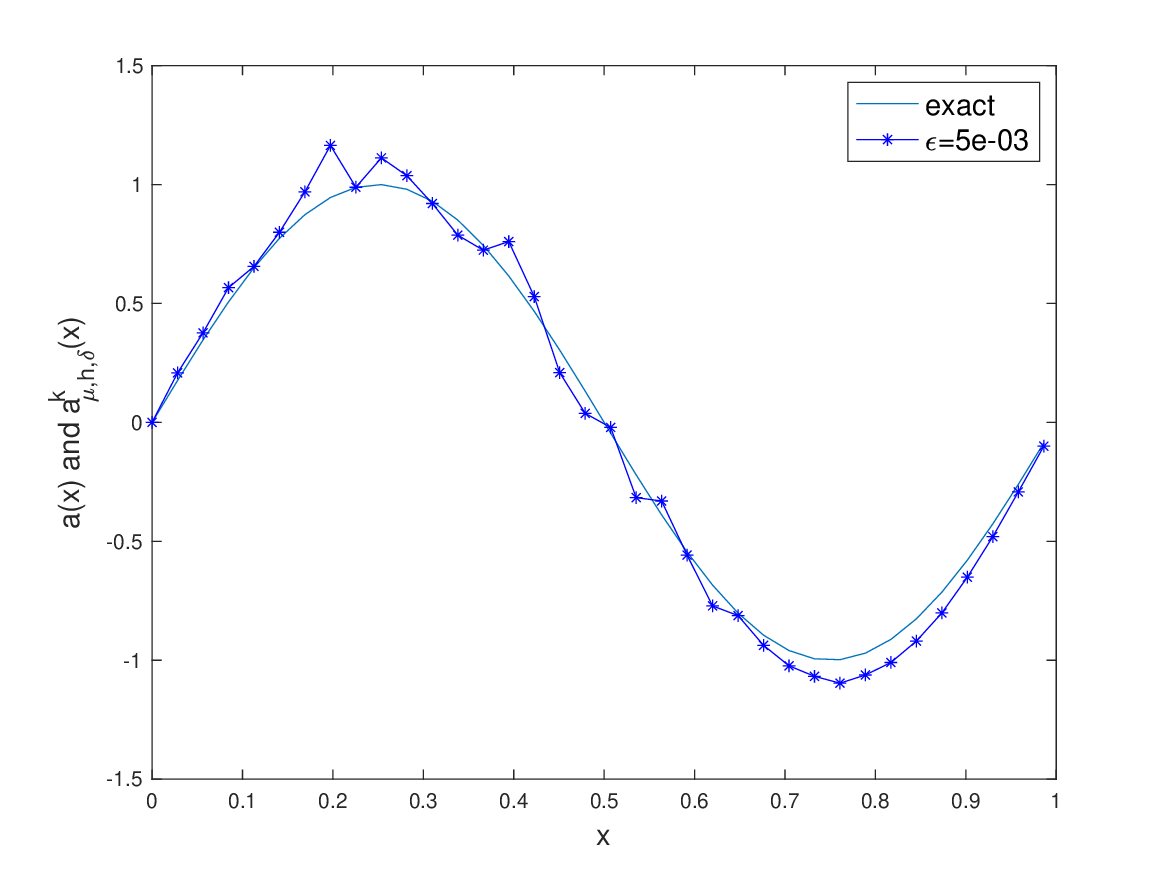}}
  \subfigure{
    \includegraphics[width=0.35\textwidth,height=1.5in]{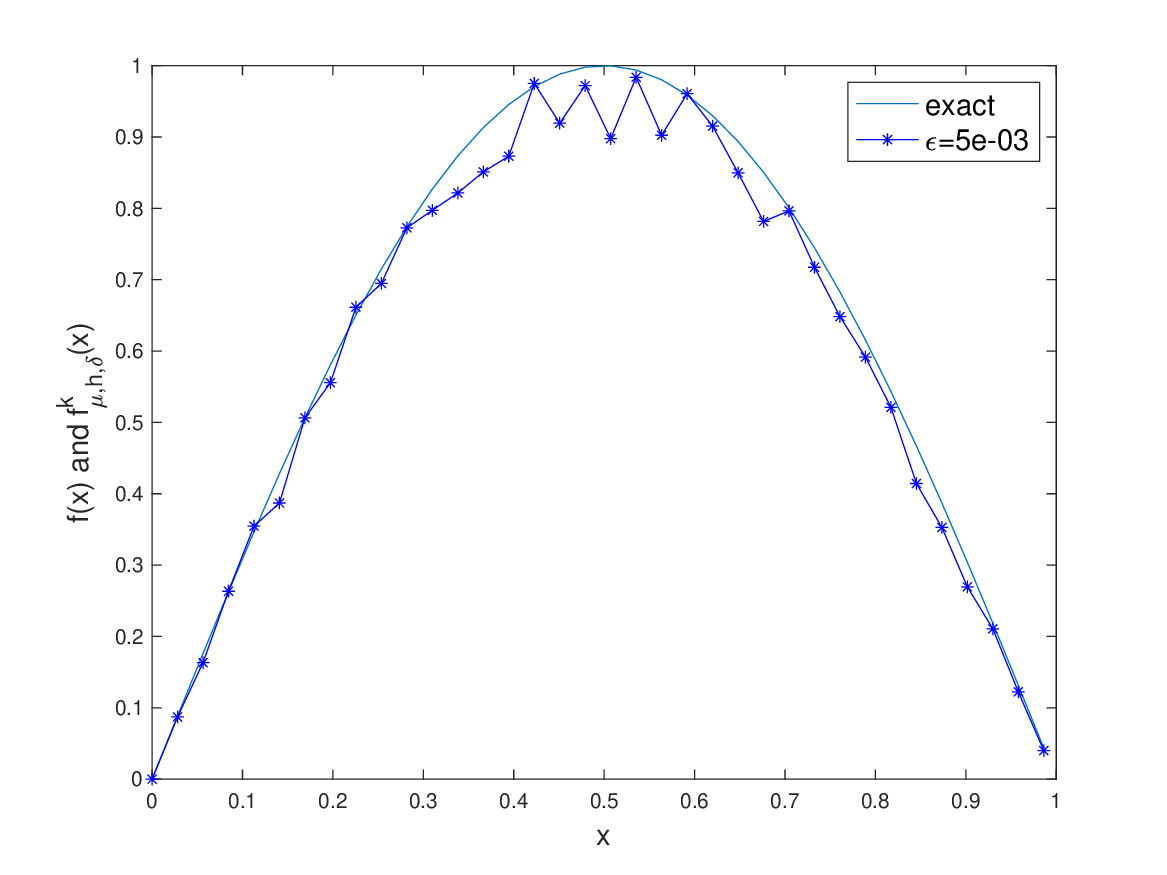}}\\
      \subfigure{
    \includegraphics[width=0.35\textwidth,height=1.5in]{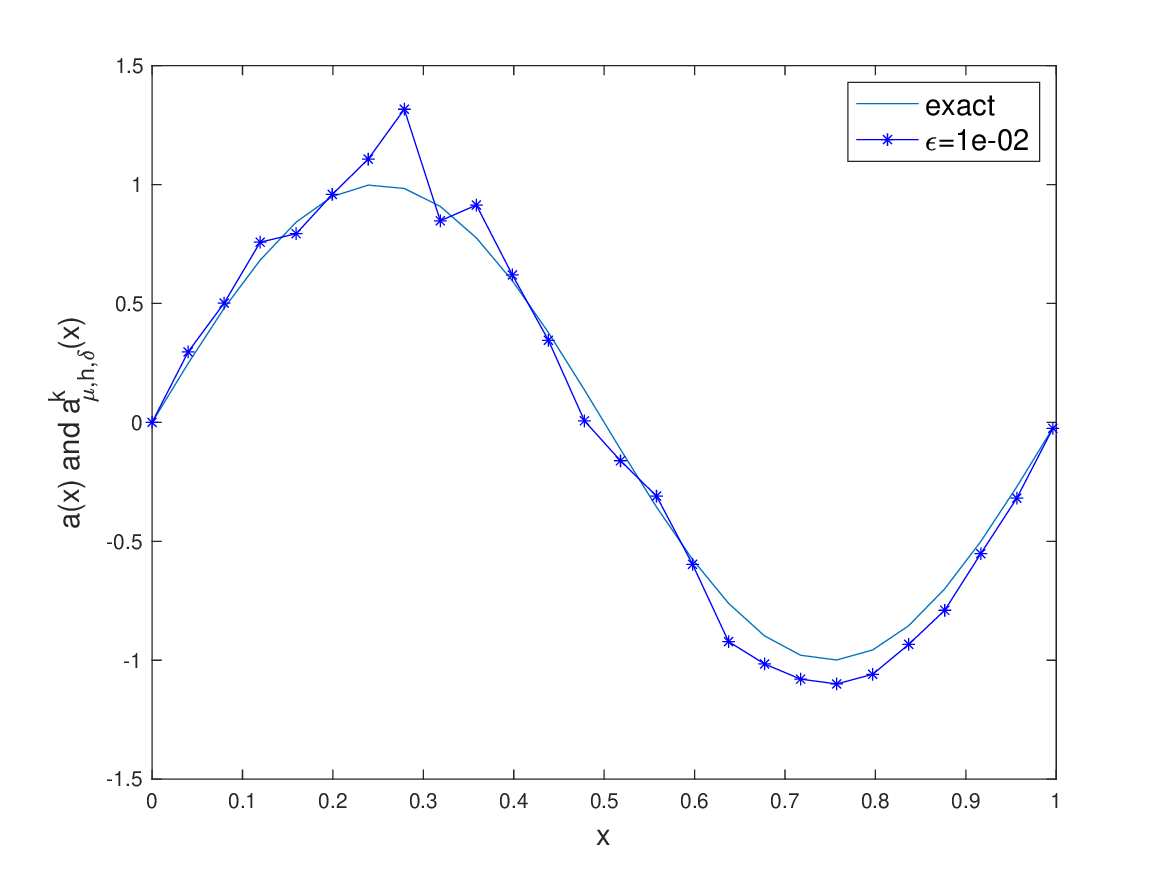}}
      \subfigure{
    \includegraphics[width=0.35\textwidth,height=1.5in]{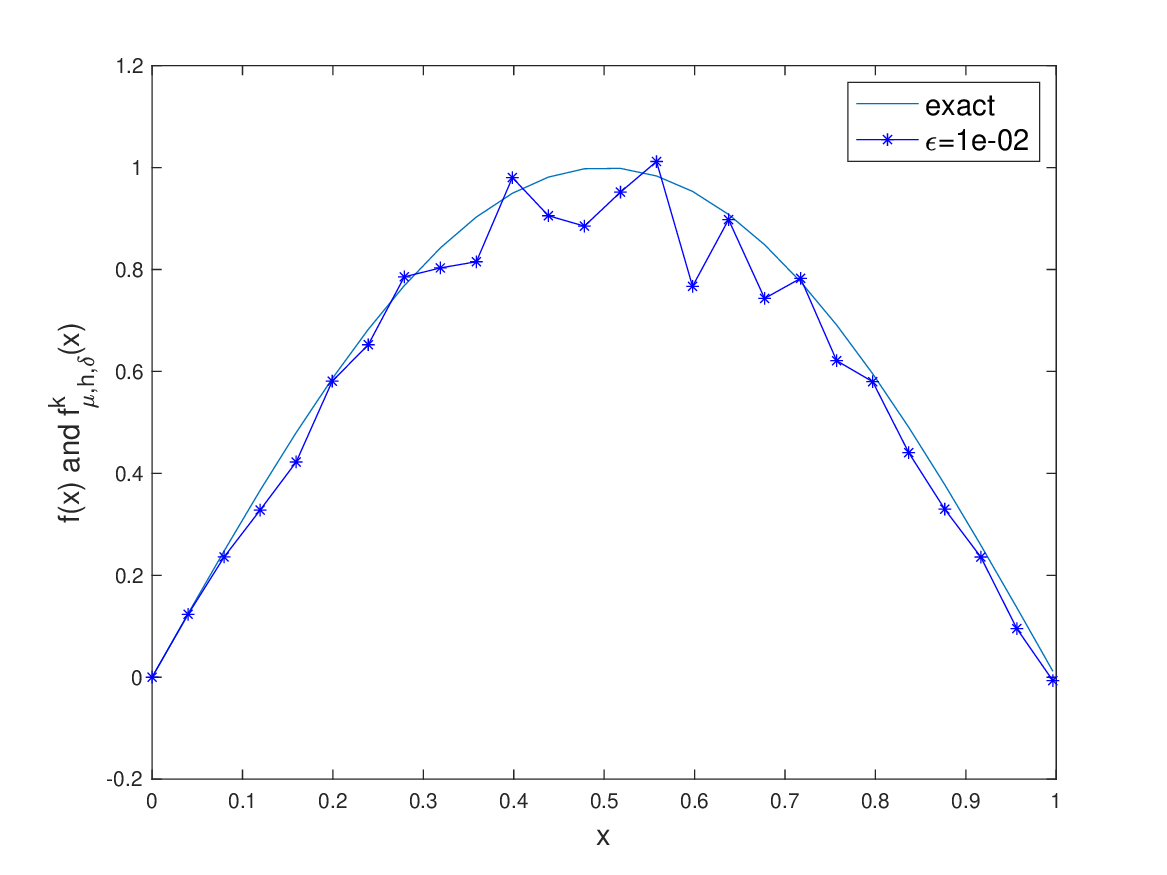}}
  \caption{Approximate solutions $a_{\mu,h,\delta}^{(k)}$ (left column) and $f_{\mu,h,\delta}^{(k)}$(right column) of Example 1 with various $\epsilon$.}\label{fig1}
\end{figure}

Since the functions $a(x), f(x)\in D((-\Delta))$, we set $h=\delta^{\frac12}$ and $\mu=0.2\delta^{\frac23}$. As mentioned in Remark \ref{rem3}, we anticipate the convergence order to be optimal at $O(\delta^{\frac12})$. Utilizing formulae (\ref{sak}) and (\ref{sfk}), we compute the approximate solutions of $(a,f)$. To verify the convergence results, we carry out a numerical experiment across various noise levels of $\epsilon=0.0005,0.001,0.005,0.01$, respectively. Figure \ref{fig1} demonstrating excellent agreement between the reconstructed solutions and the exact solutions. As $\epsilon$ increases, the grids become sparser due to the growing of mesh size $h$. The relative errors and convergence orders are given in Table \ref{tab1}. It can be seen that as $\epsilon$ increases, the relative errors also tend to be large. And the numerical findings presented in Table \ref{tab1} align closely with the theoretical results outlined in Theorem \ref{main2} for reconstructing smooth initial values and source terms.

\begin{table}[h]
\centering
\begin{tabular}{ccccc}
\hline
$\epsilon$ & $re_a$ & Order(a) & $re_f$ & Order(f) \\
\hline
$0.0005$ & 0.0331 & 0.4645 & 0.0211 & 0.7708 \\
\hline
$0.001$ & 0.0488 & 0.4786 & 0.0290 & 0.5036 \\
\hline
$0.005$ & 0.01093 & 0.3614 & 0.0546 & 0.4631 \\
\hline
$0.01$ & 0.1466 & 0.3430 & 0.0876 & 0.4175 \\
\hline
%$0.0006$ & 0.0390 & $\delta=9.5279e-05$ & 0.0213 & $\delta=9.5279e-05$  \\
%\hline
\end{tabular}
\caption{The relative errors and convergence orders based on standard Galerkin method of Example 1.}\label{tab1}
\end{table}

\noindent\textbf{Example 2.} Let $\alpha=1.5$, $d=0.3$, and $T_2=1.5$. We take\\
\begin{minipage}{0.5\textwidth}
\begin{equation*}
a(x)=\left\{
\begin{aligned}
&3x,& &0<x\leq0.3,\\
&0.9,& &0.3<x\leq0.7,\\
&3(1-x),& &0.7<x\leq1,
\end{aligned}
\right.
\end{equation*}
\end{minipage}%
\begin{minipage}{0.5\textwidth}
\begin{equation*}
f(x)=\left\{
\begin{aligned}
&x,& &0<x\leq0.5,\\
&1-x,& &0.5<x\leq1.
\end{aligned}
\right.
\end{equation*}
\end{minipage}

\begin{figure}[htbp]
  \centering
  \subfigure{
    \includegraphics[width=0.35\textwidth,height=1.5in]{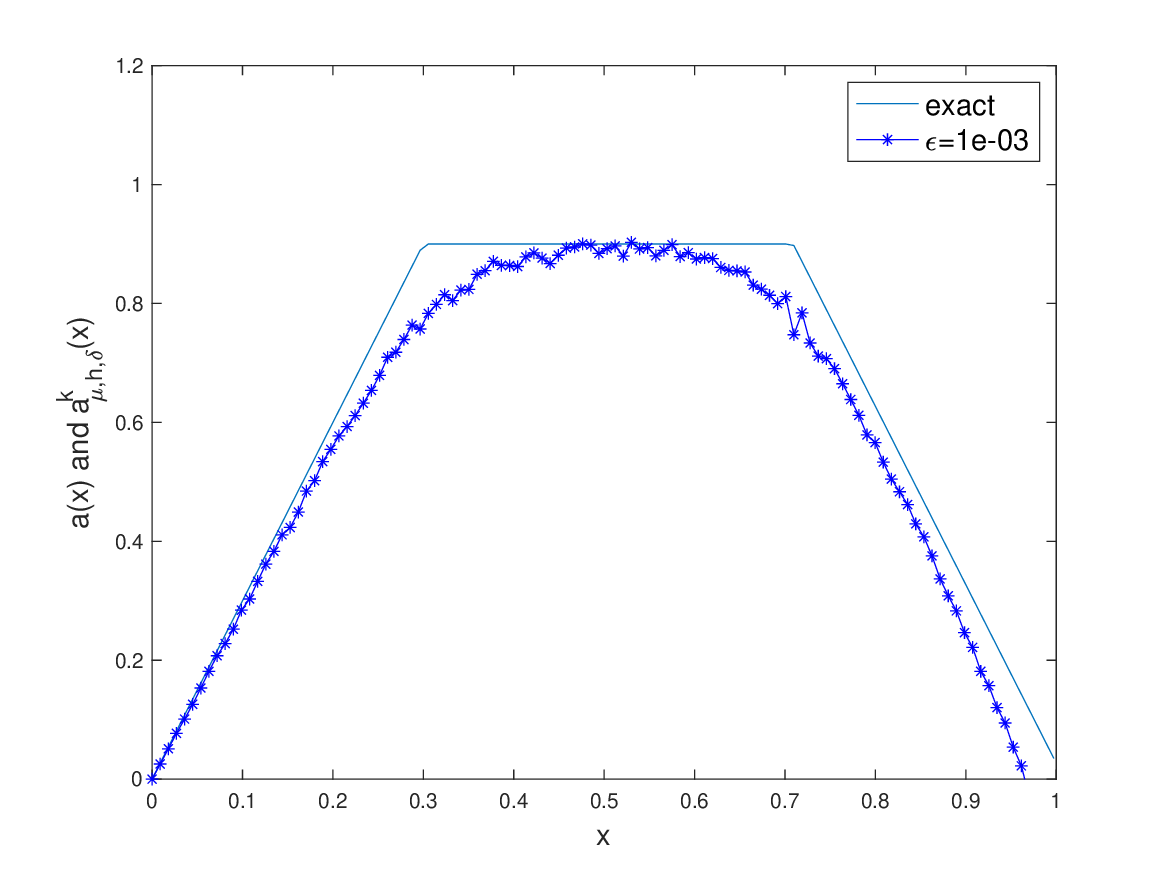}}
  \subfigure{
    \includegraphics[width=0.35\textwidth,height=1.5in]{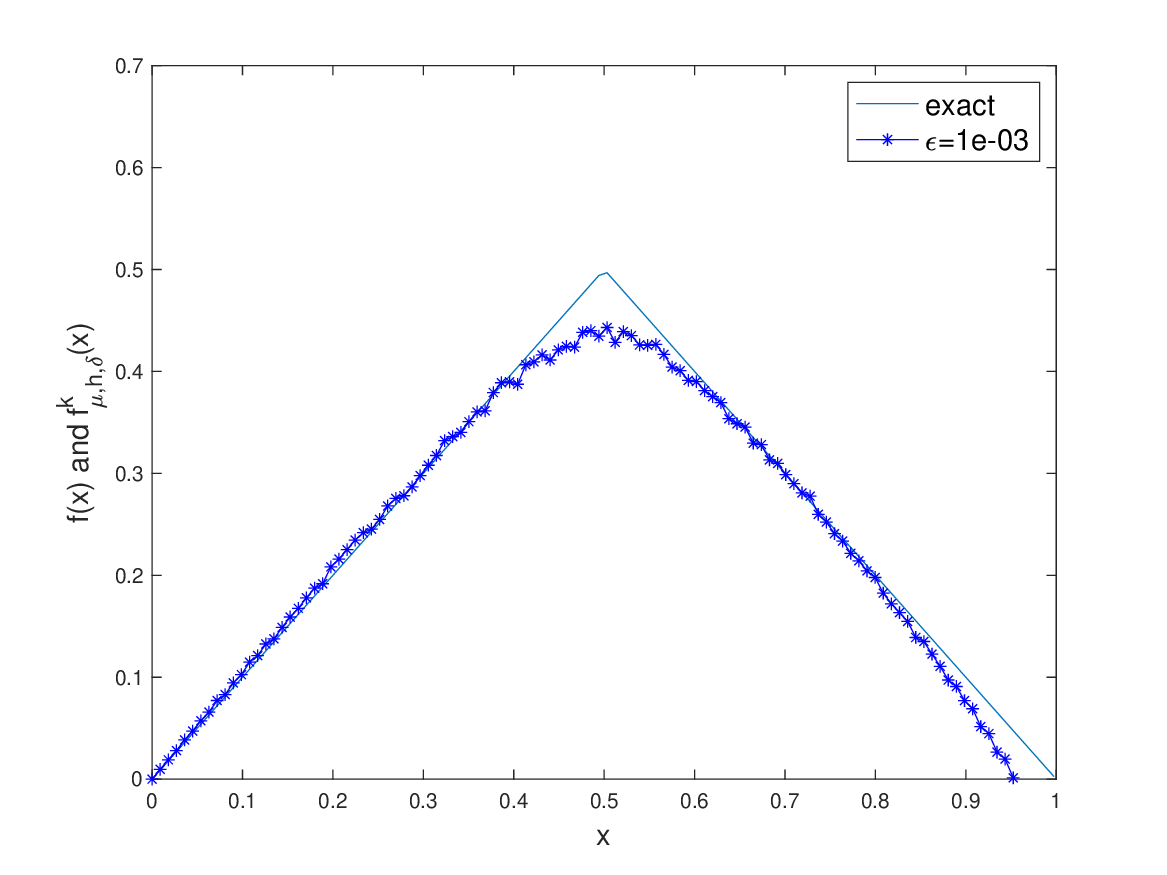}}\\
      \subfigure{
    \includegraphics[width=0.35\textwidth,height=1.5in]{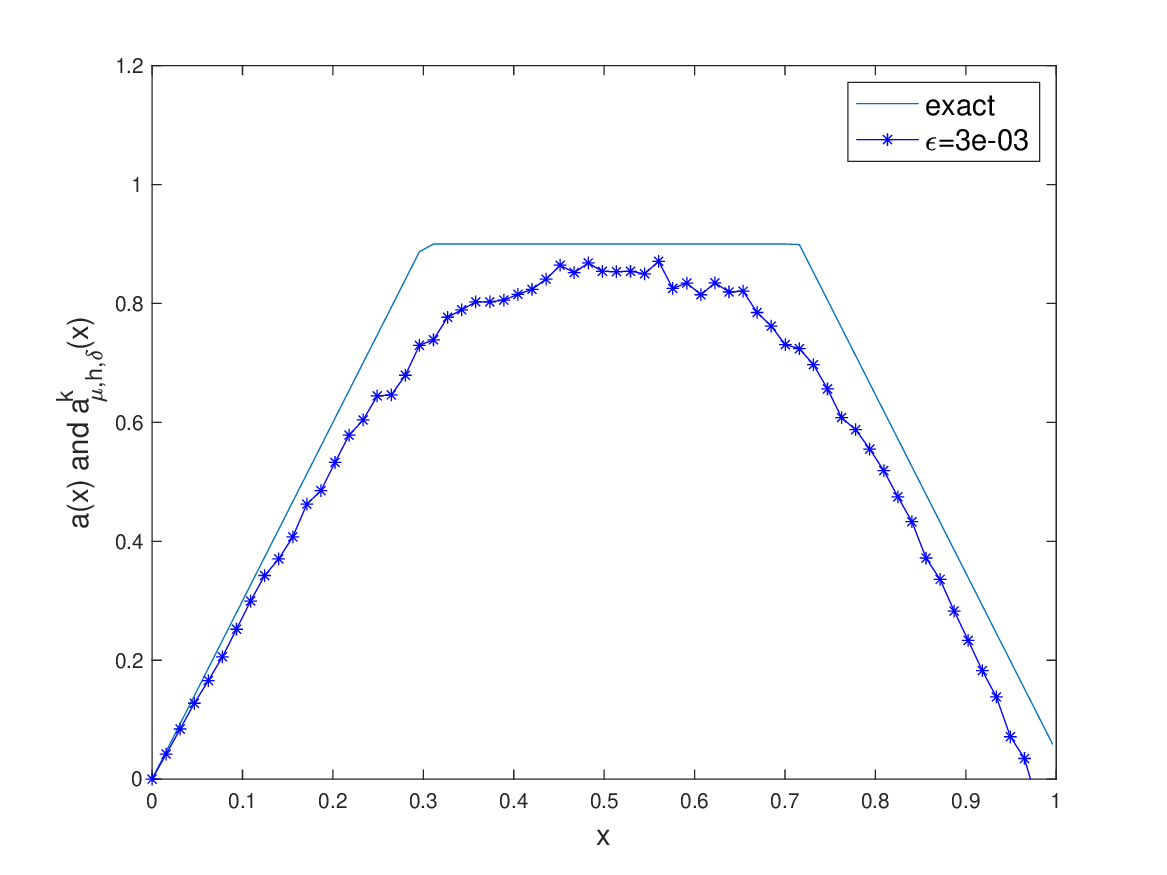}}
      \subfigure{
    \includegraphics[width=0.35\textwidth,height=1.5in]{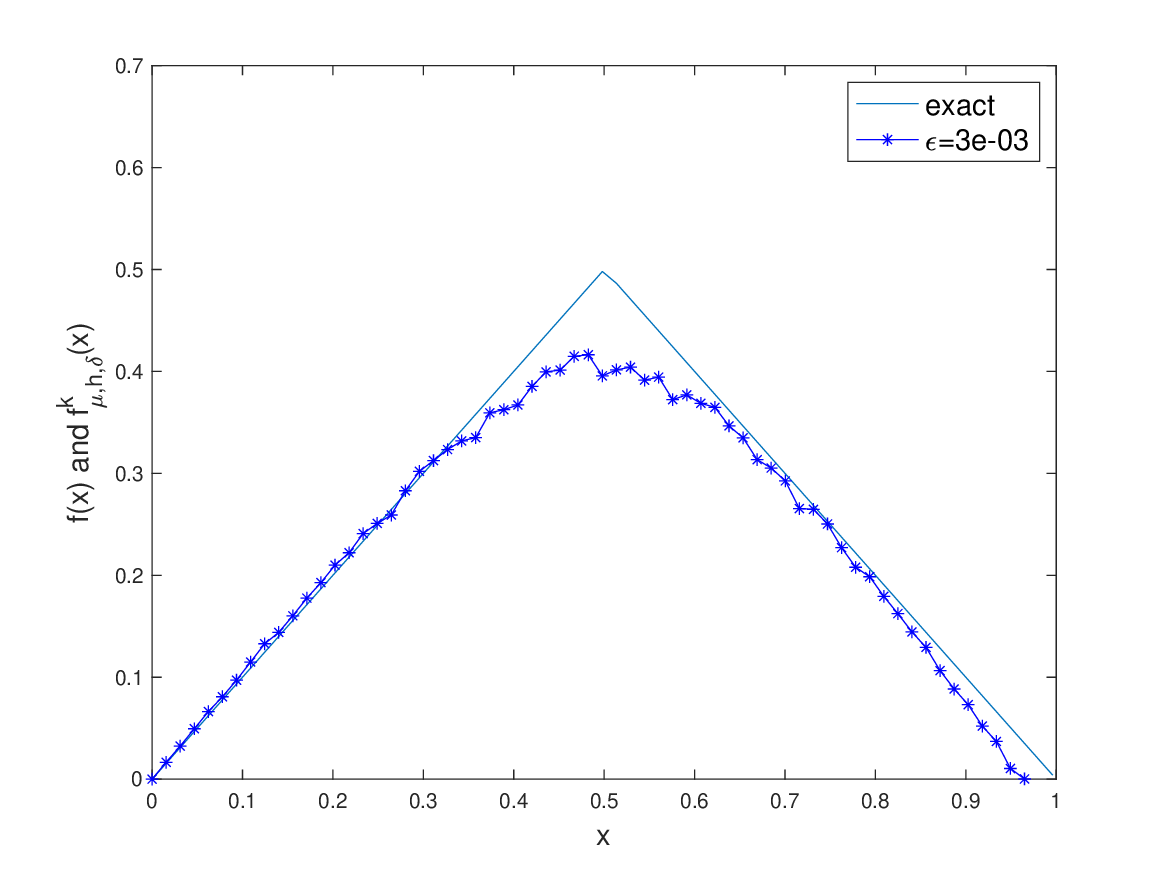}}\\
  \subfigure{
    \includegraphics[width=0.35\textwidth,height=1.5in]{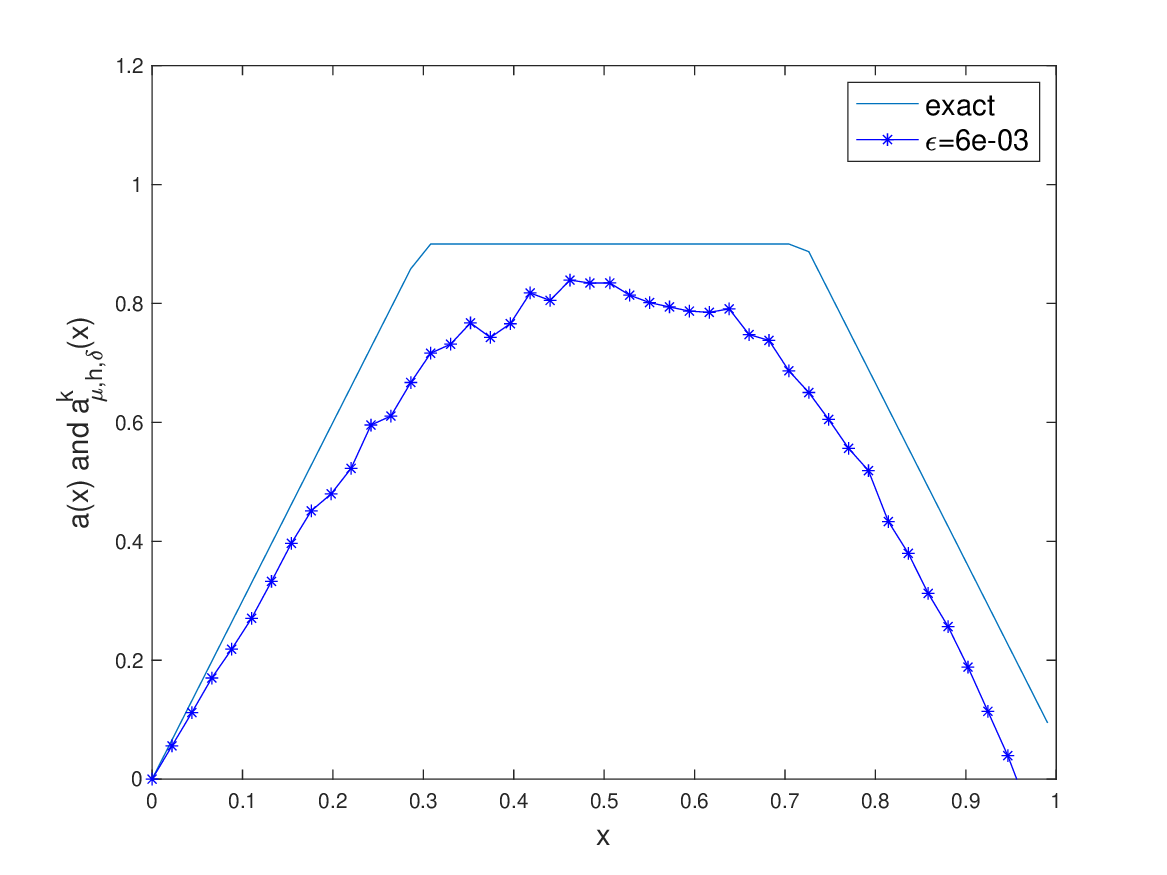}}
  \subfigure{
    \includegraphics[width=0.35\textwidth,height=1.5in]{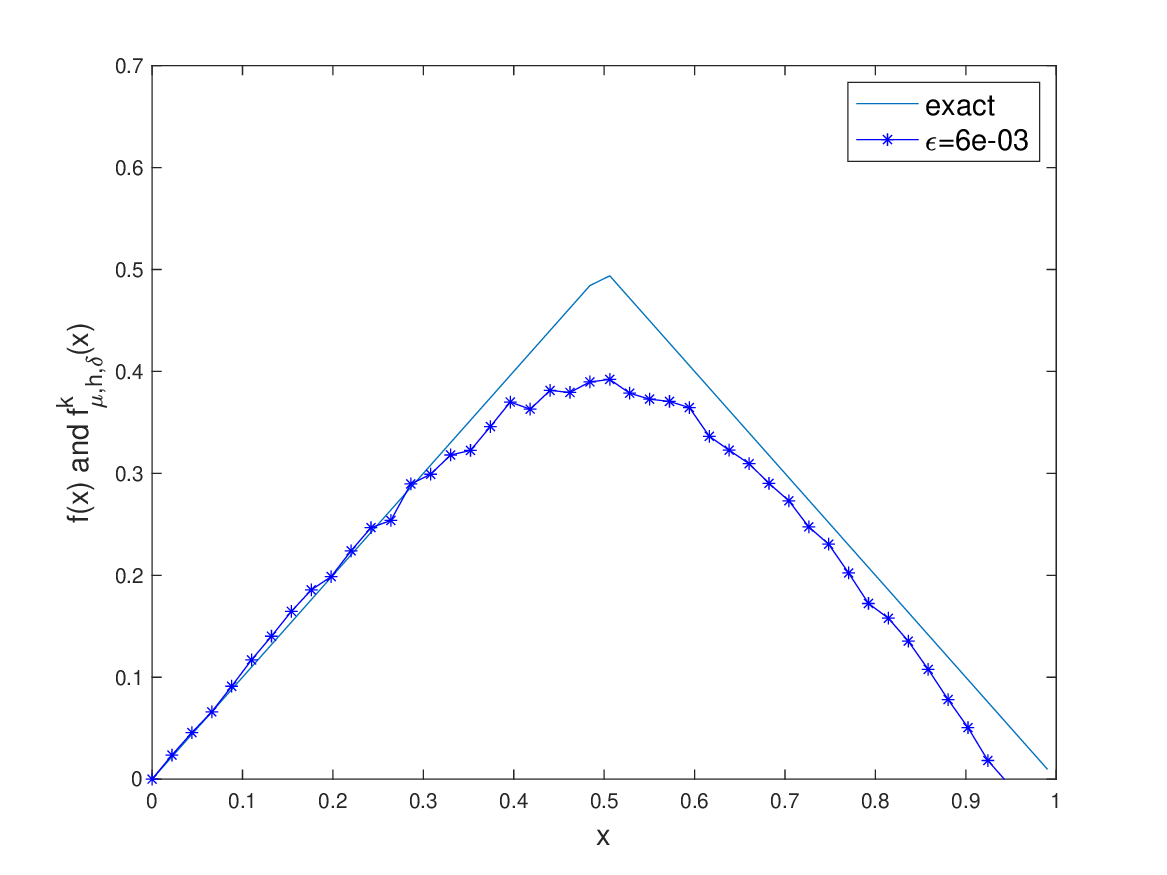}}\\
      \subfigure{
    \includegraphics[width=0.35\textwidth,height=1.5in]{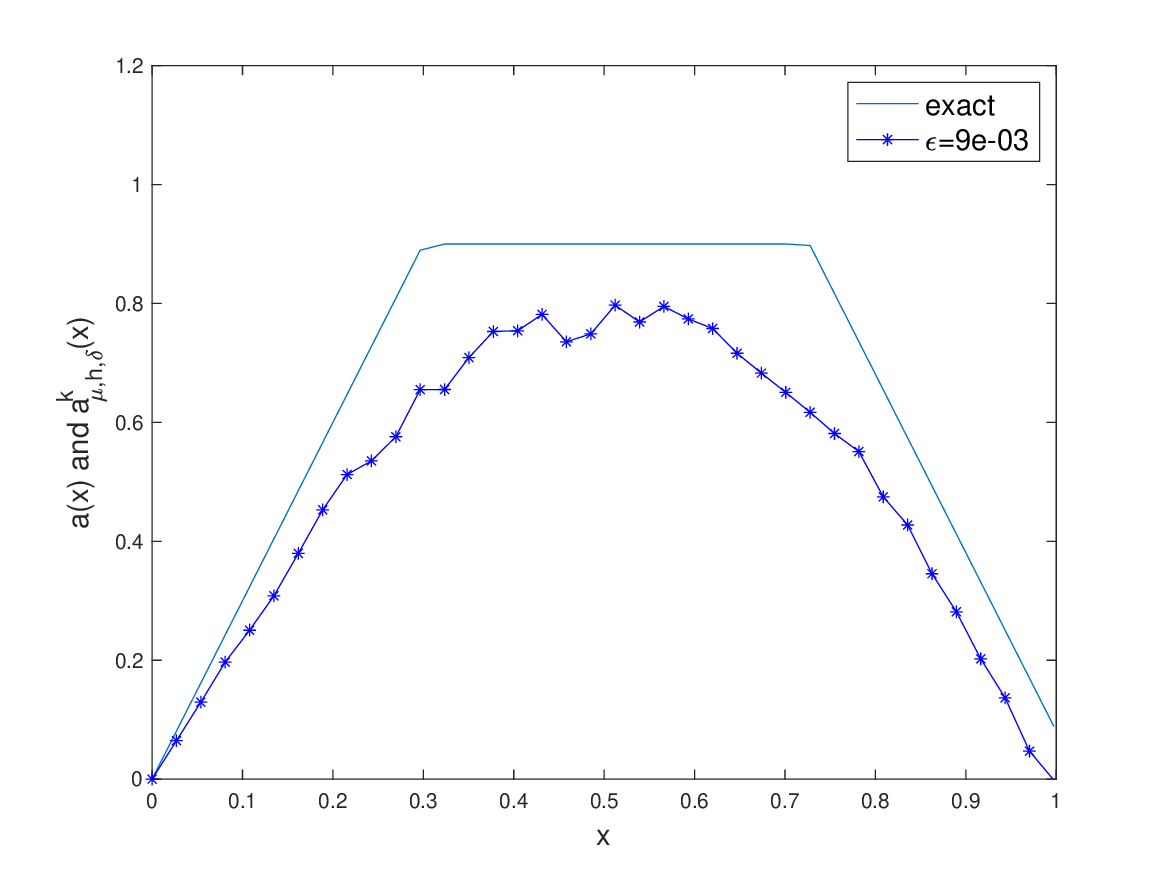}}
      \subfigure{
    \includegraphics[width=0.35\textwidth,height=1.5in]{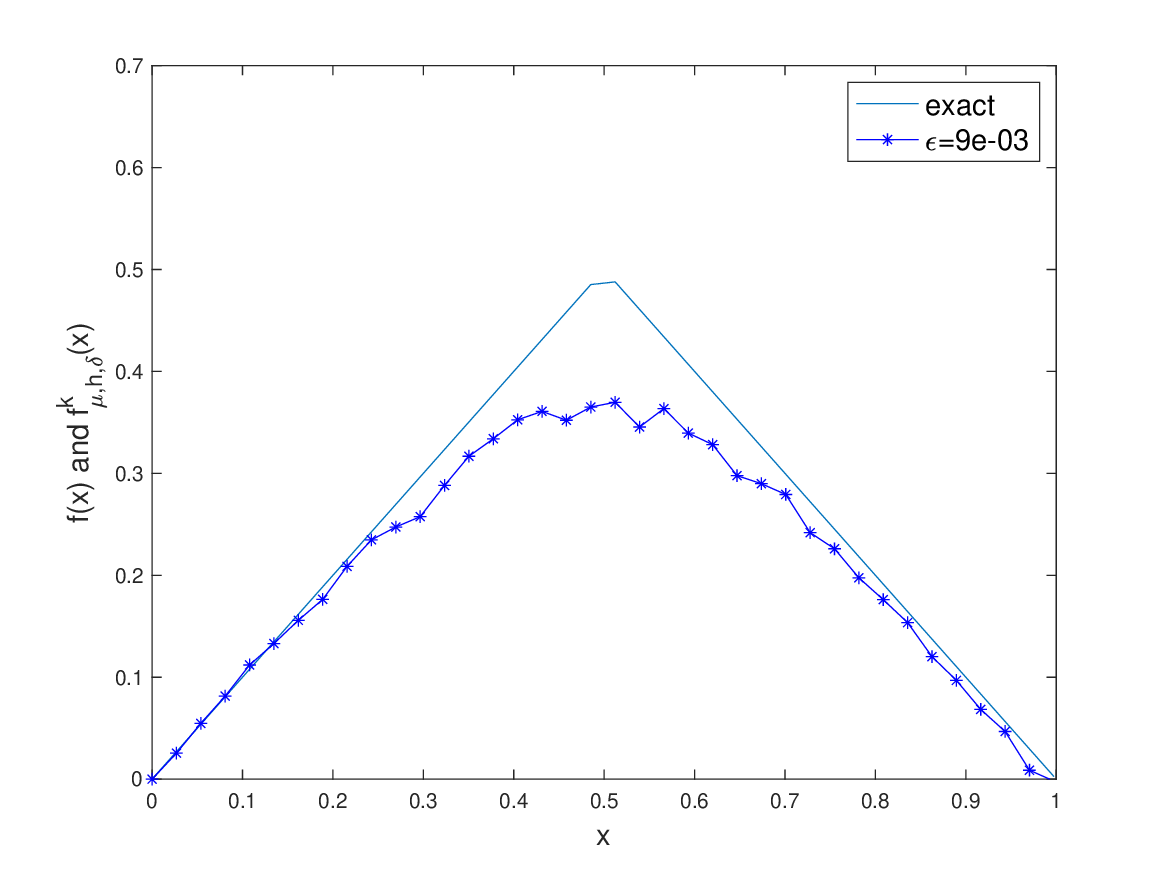}}      
    \subfigure{
    \includegraphics[width=0.35\textwidth,height=1.5in]{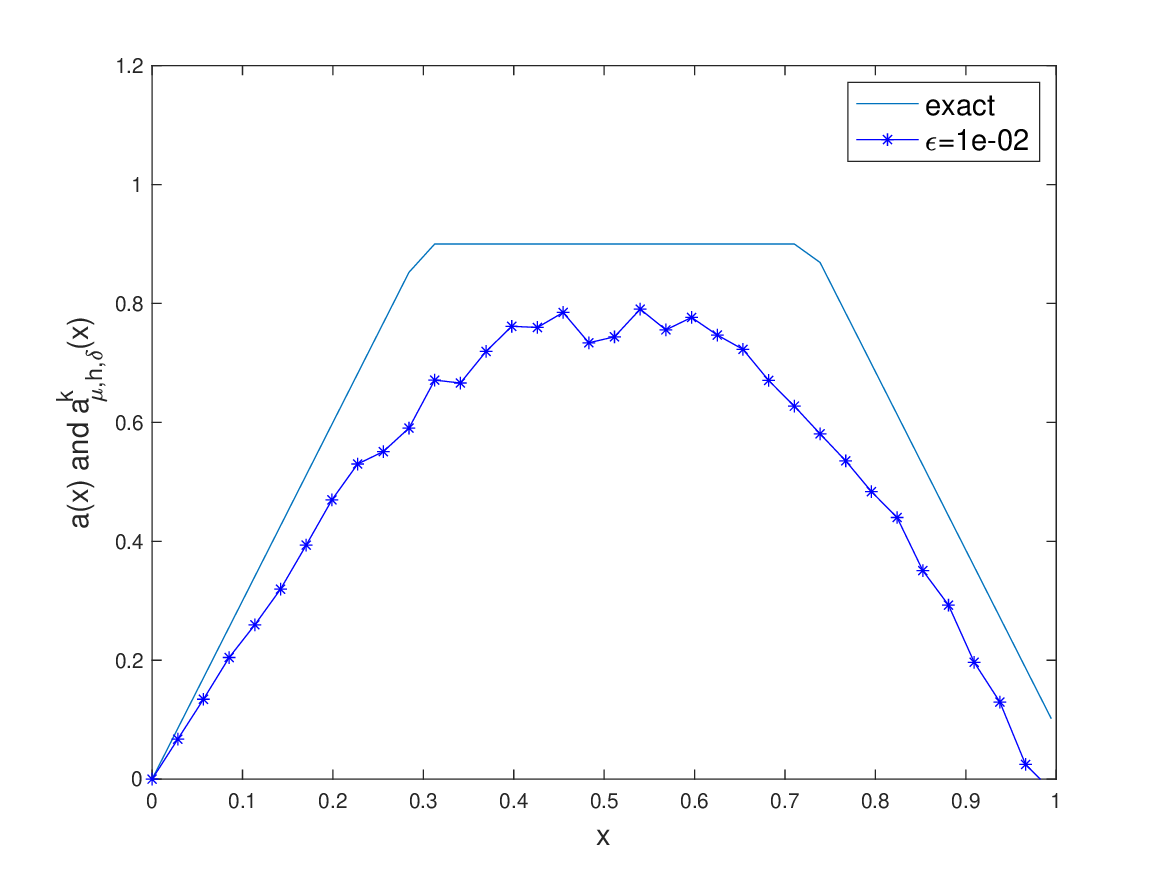}}
     \subfigure{
    \includegraphics[width=0.35\textwidth,height=1.5in]{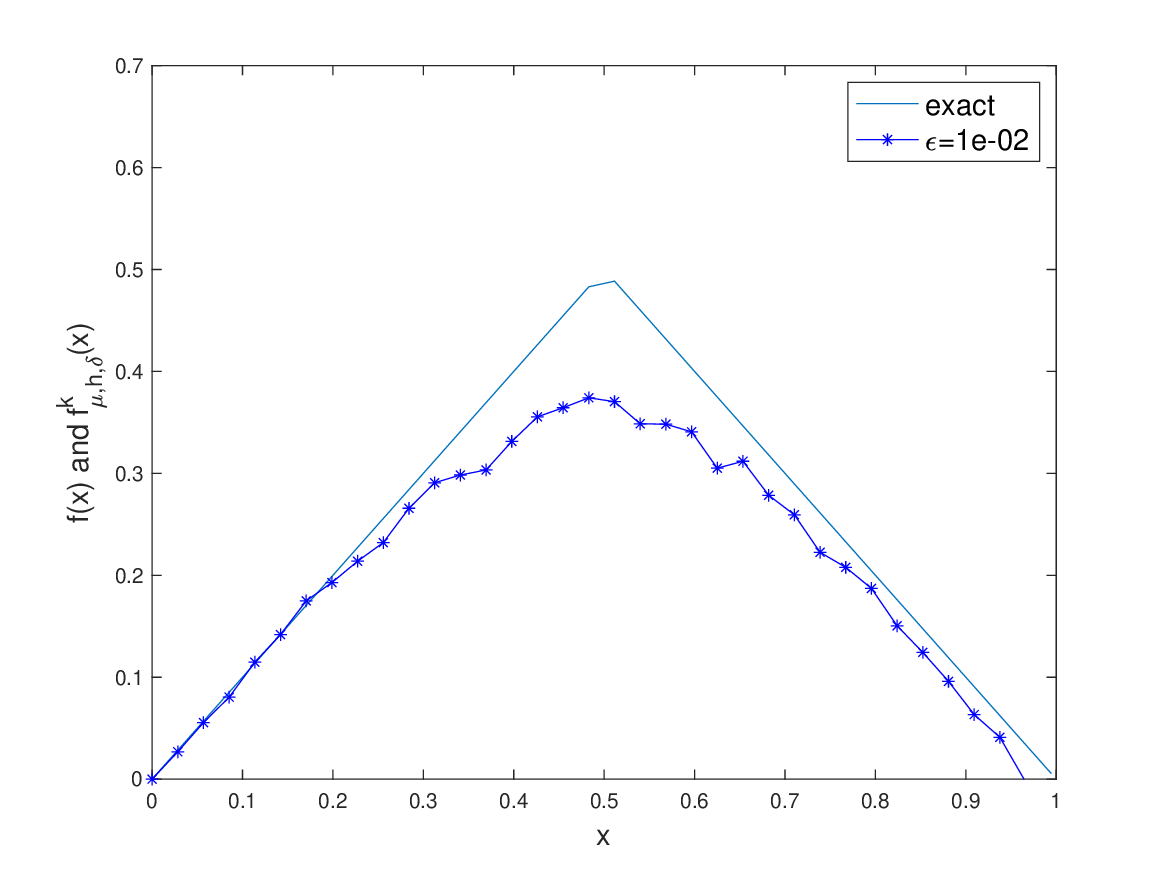}}
  \caption{Approximate solutions $a_{\mu,h,\delta}^{(k)}$ (left column) and $f_{\mu,h,\delta}^{(k)}$(right column) of Example 2 based on standard Galerkin method with various $\epsilon$.}\label{fig2}
\end{figure}

In this example, we consider the initial value and source term to be nonsmooth functions. Consequently, the exact solution of problem (\ref{IP}) can be given by an infinite series in the form of (\ref{solution2}). To obtain the additional measurements, we apply the finite difference method proposed in \cite{Du+Cao+Sun-2010} for solving problem (\ref{IP}). Since $a(x),f(x)\in D((-\Delta)^{\frac12})$, we set $h=0.8\delta^{\frac12}$ and $\mu=3\delta^{\frac23}$, and expect a convergence order of $O(\delta^{\frac13})$. Then we apply formulae (\ref{sak}) and (\ref{sfk}) to compute the approximate solutions of the exact solutions $(a,f)$. This numerical experiment is carried out across various noise levels of $\epsilon=0.001,0.003,0.006,0.009,0.01$, respectively.  As illustrated in Figure \ref{fig2}, the reconstructed solutions exhibit strong agreement with the exact solutions across different noise levels. The relative errors and convergence orders reported in Table \ref{tab2} align with the theoretical expectations, affirming the reliability of the proposed method for nonsmooth unknown initial values and source terms.

\begin{figure}[htbp]
  \centering
  \subfigure{
    \includegraphics[width=0.35\textwidth,height=1.5in]{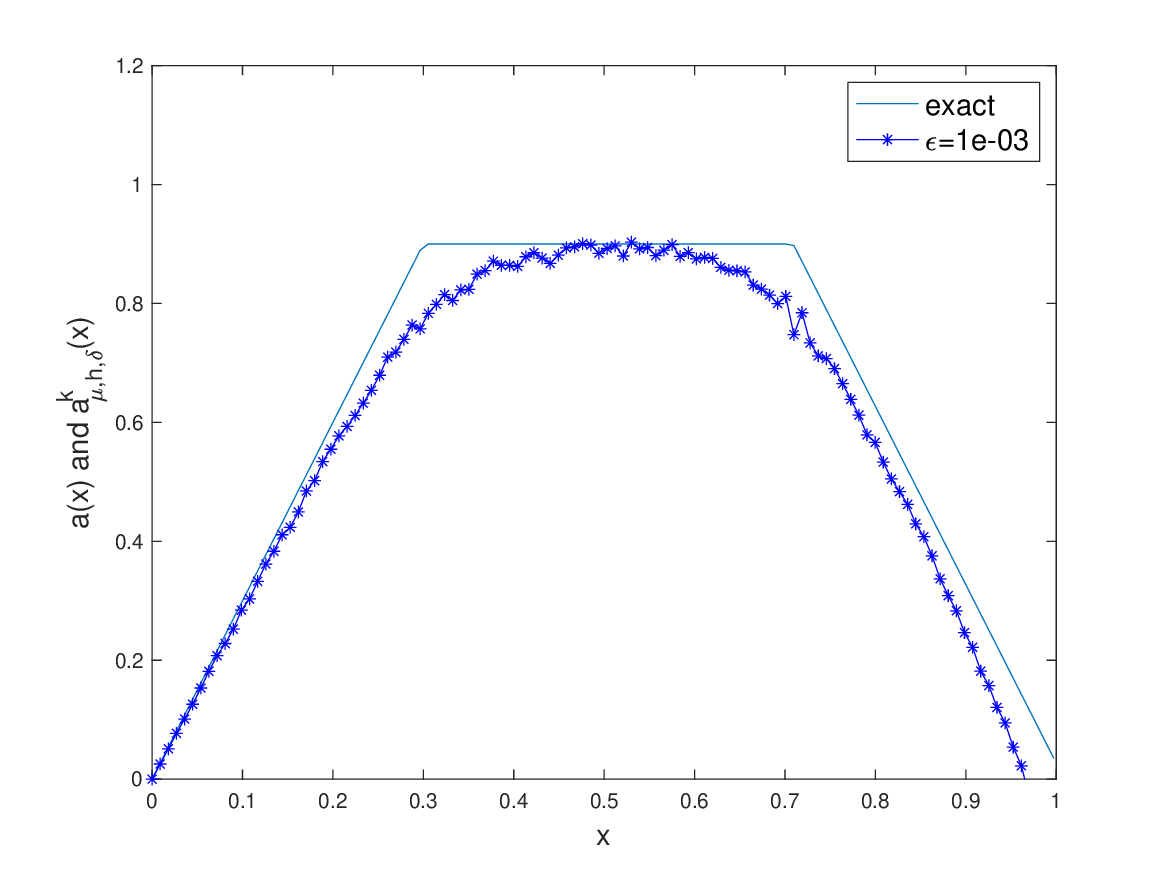}}
  \subfigure{
    \includegraphics[width=0.35\textwidth,height=1.5in]{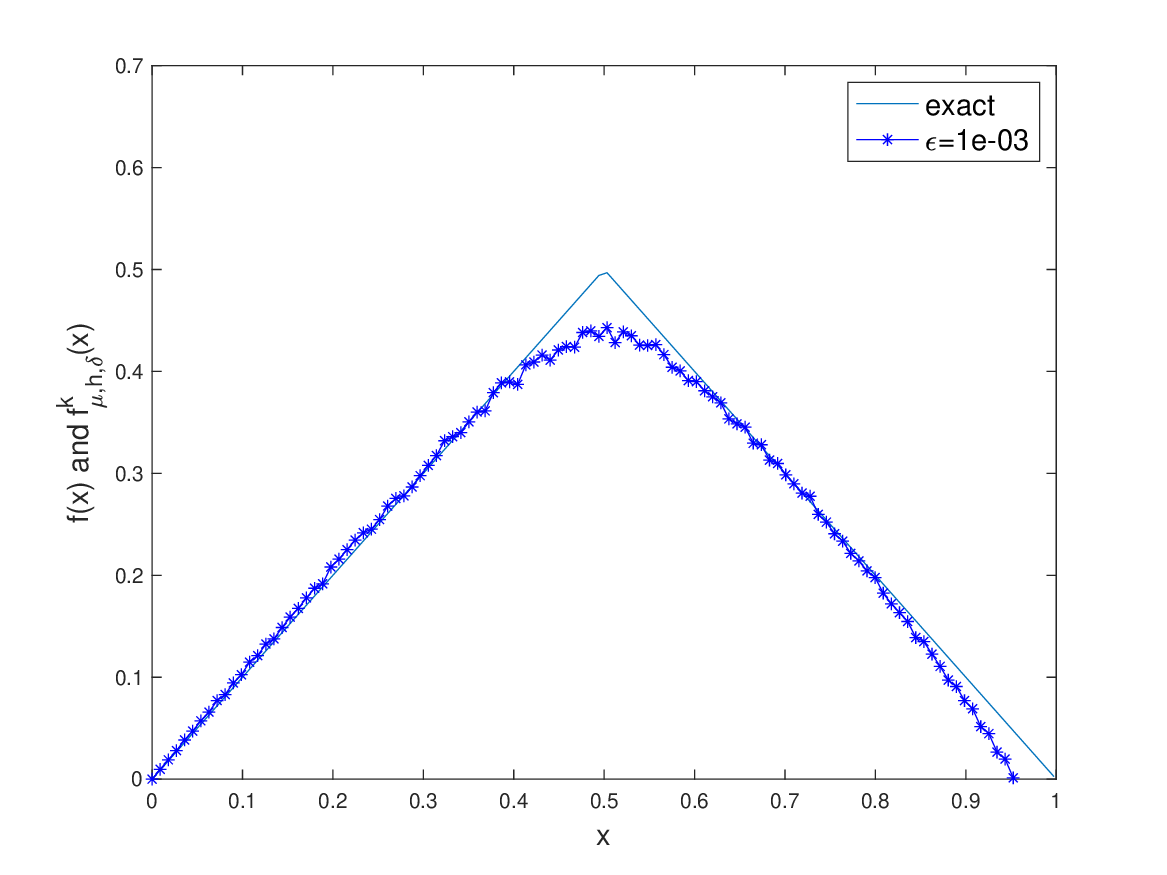}}\\
      \subfigure{
    \includegraphics[width=0.35\textwidth,height=1.5in]{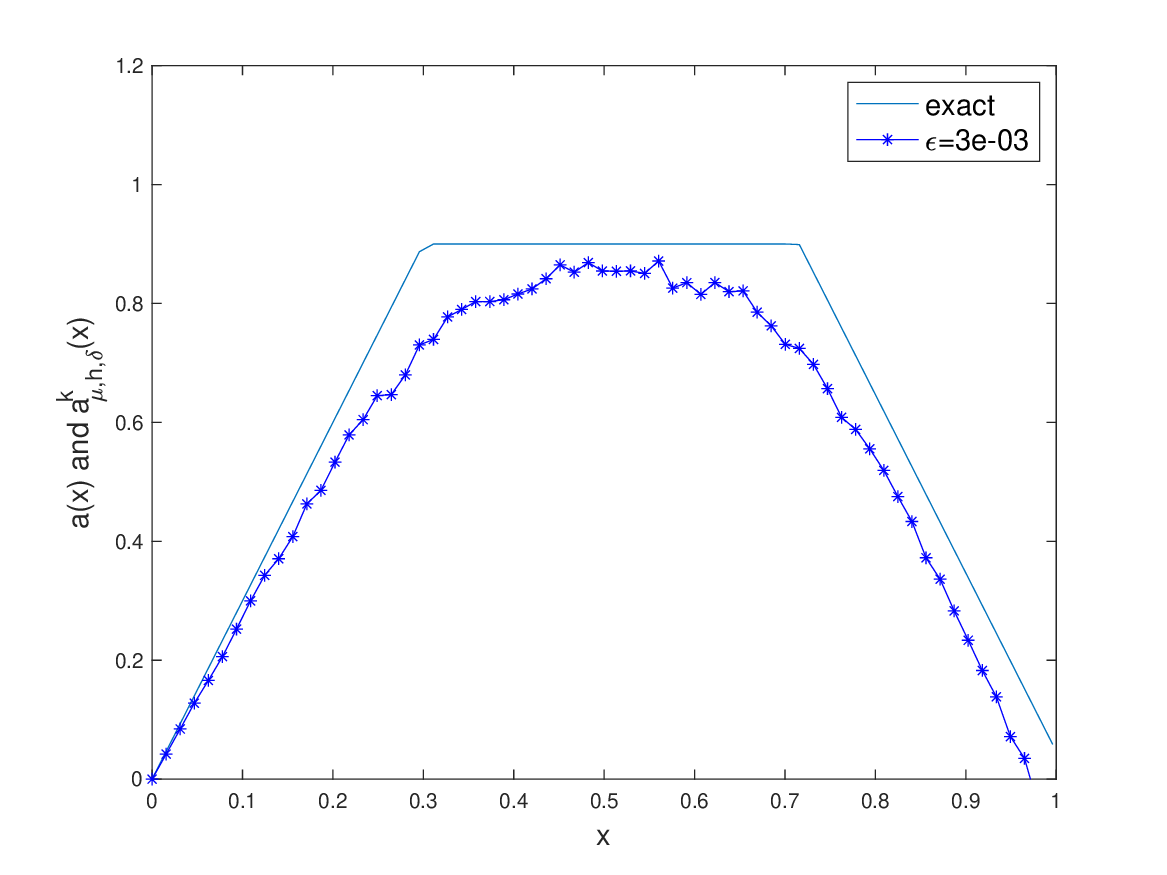}}
      \subfigure{
    \includegraphics[width=0.35\textwidth,height=1.5in]{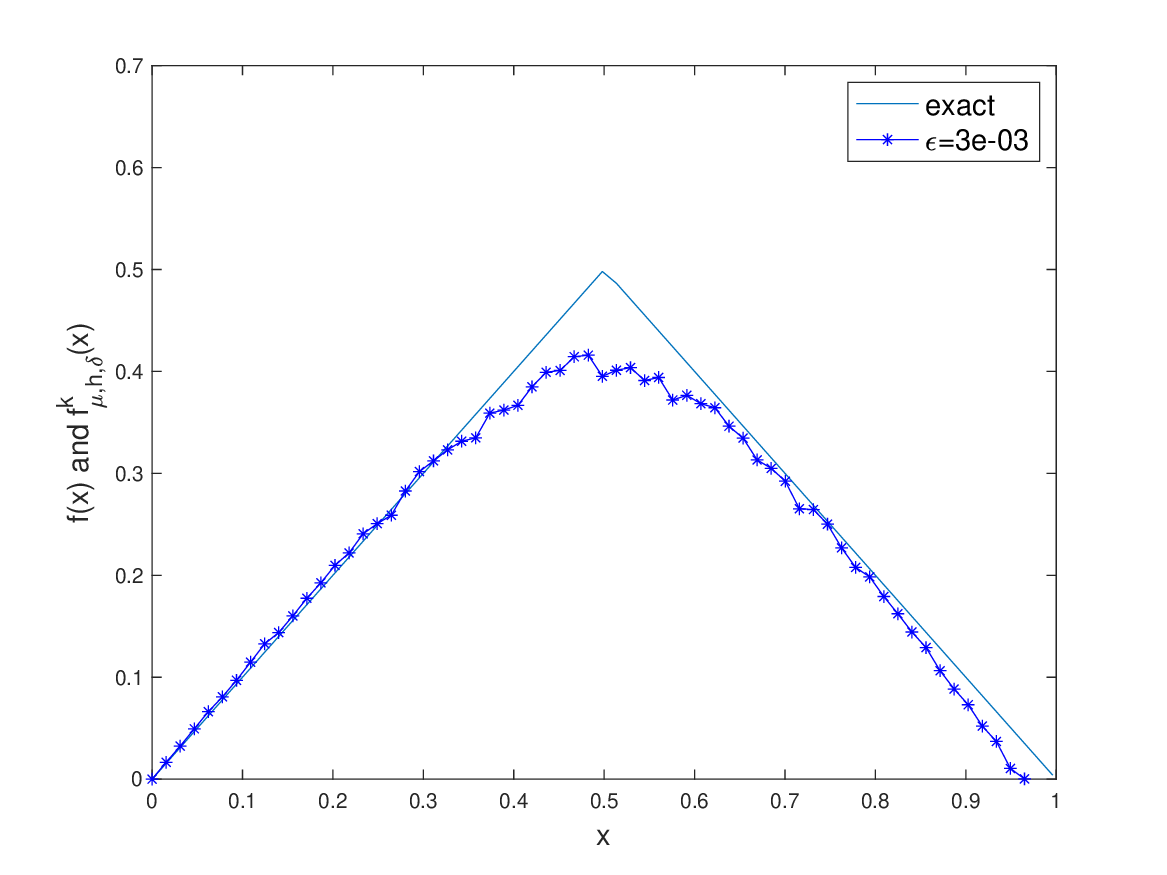}}\\
  \subfigure{
    \includegraphics[width=0.35\textwidth,height=1.5in]{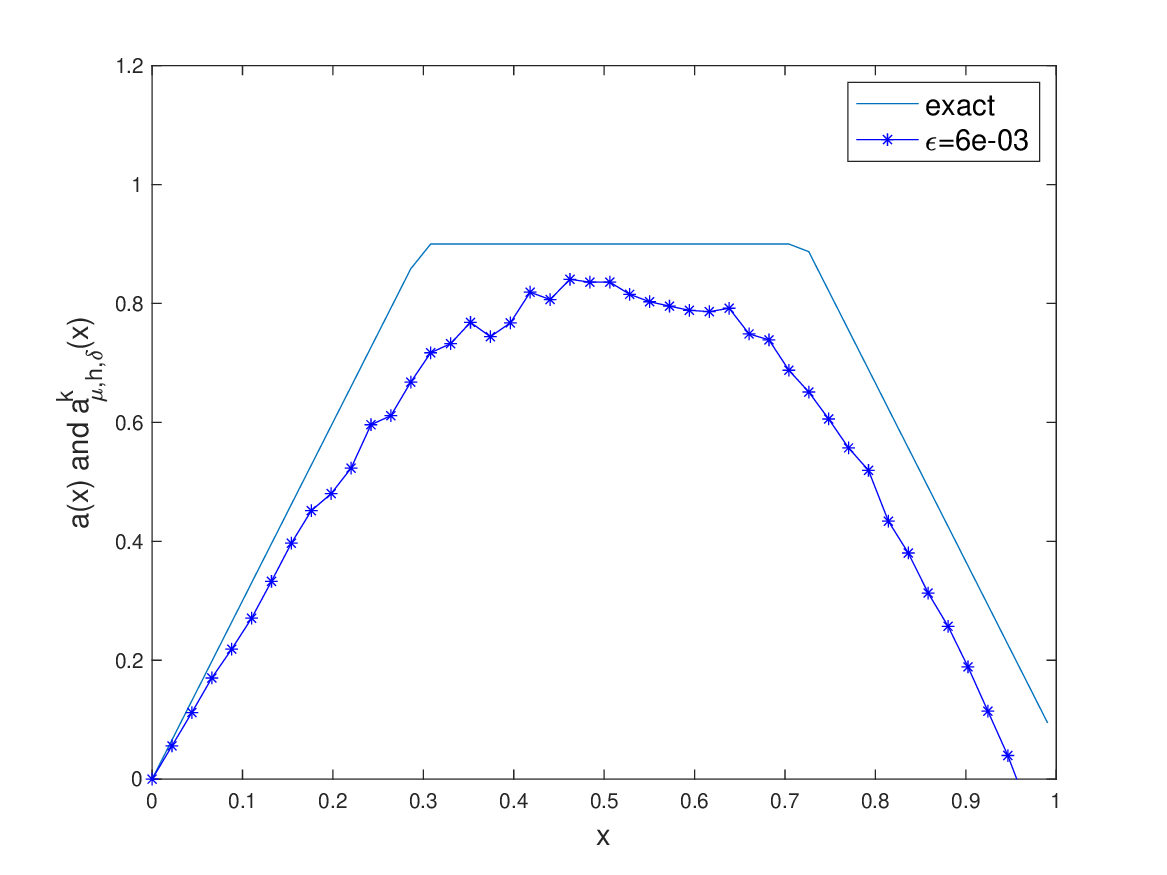}}
  \subfigure{
    \includegraphics[width=0.35\textwidth,height=1.5in]{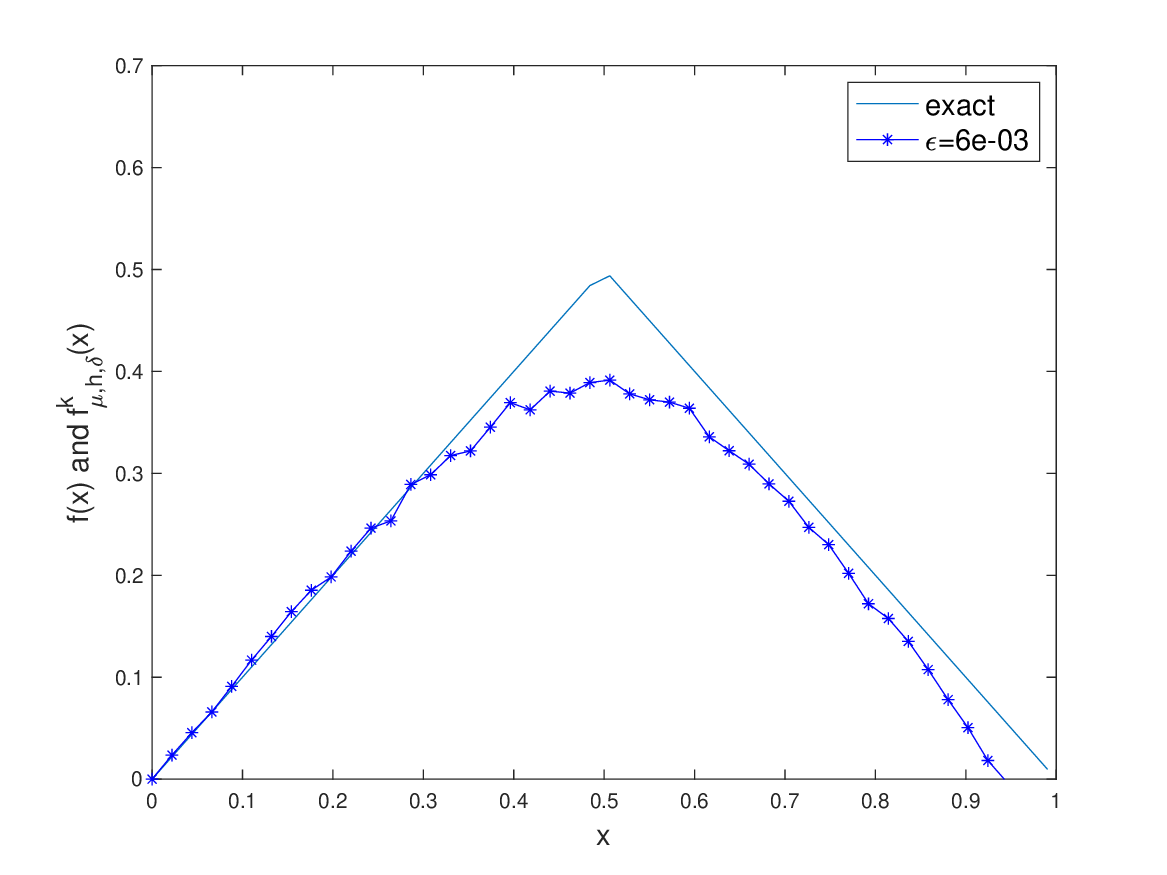}}\\
      \subfigure{
    \includegraphics[width=0.35\textwidth,height=1.5in]{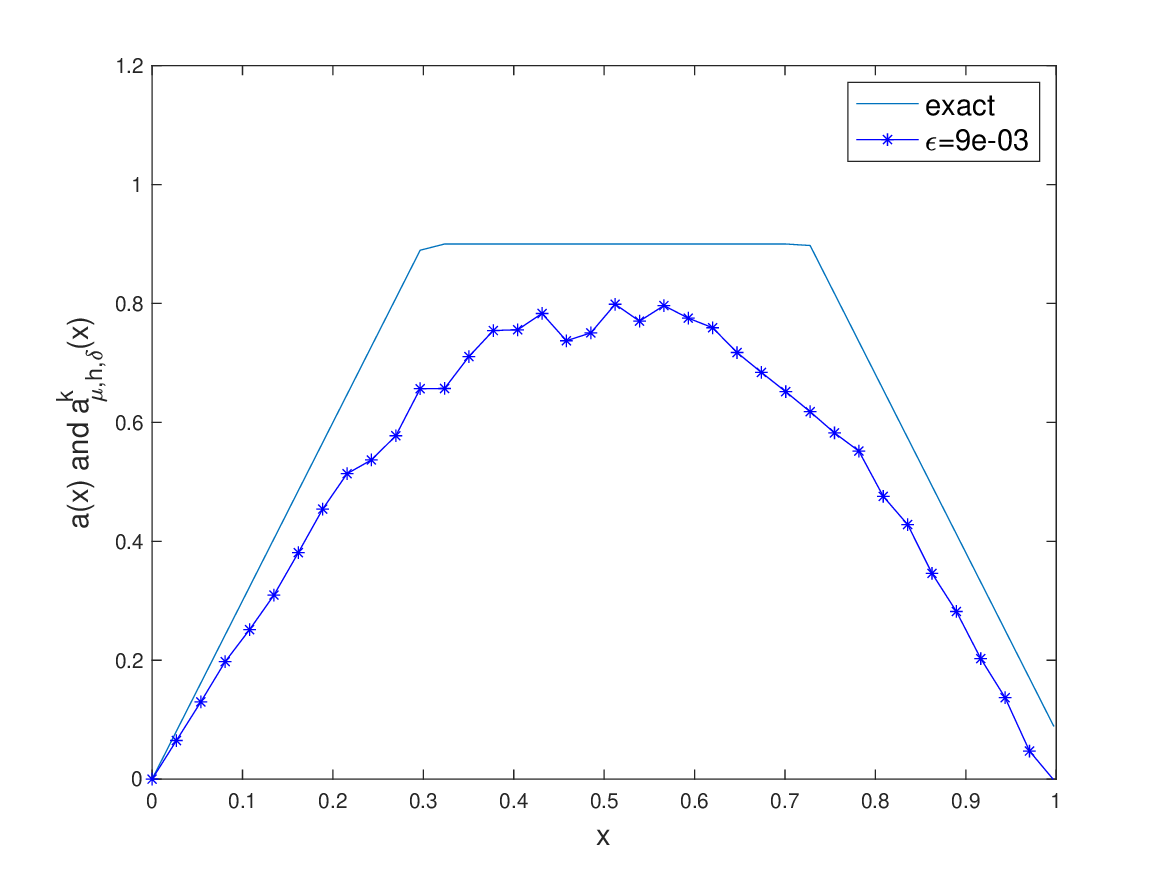}}
      \subfigure{
    \includegraphics[width=0.35\textwidth,height=1.5in]{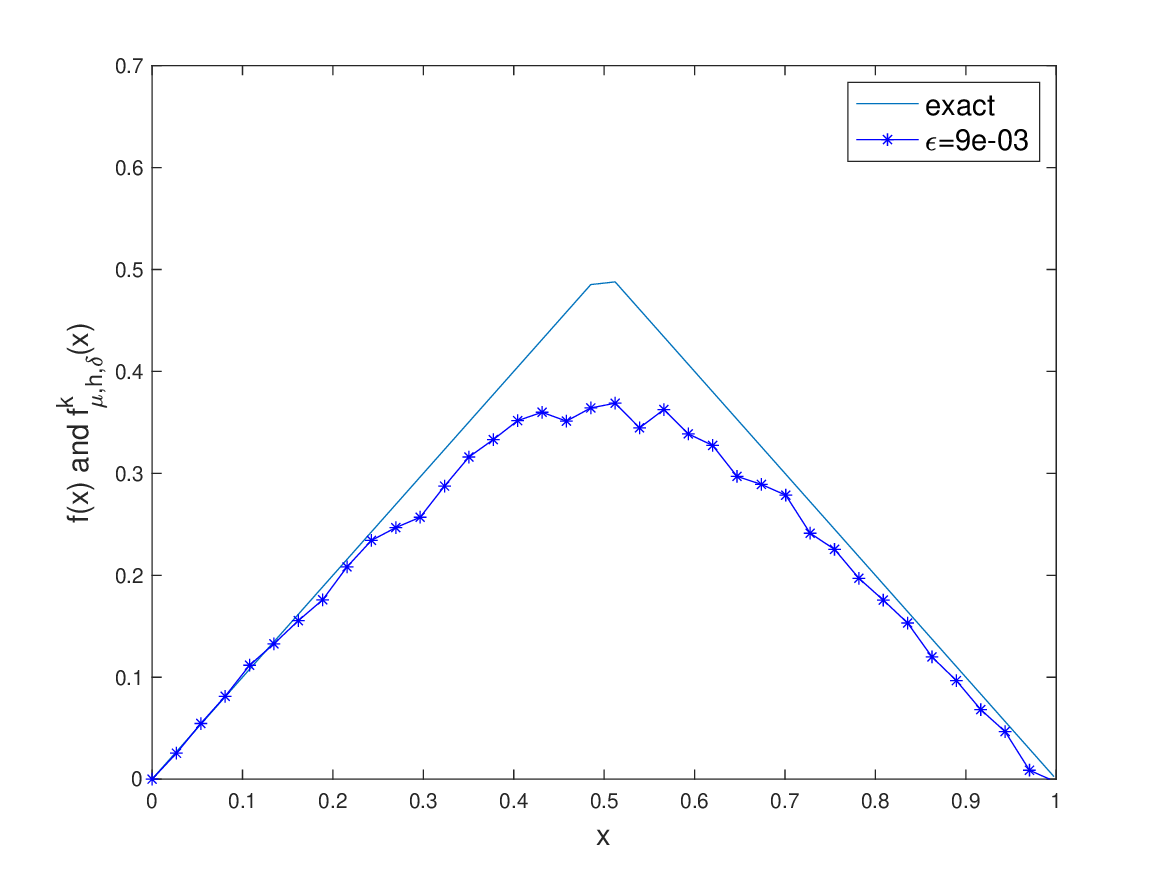}}      
    \subfigure{
    \includegraphics[width=0.35\textwidth,height=1.5in]{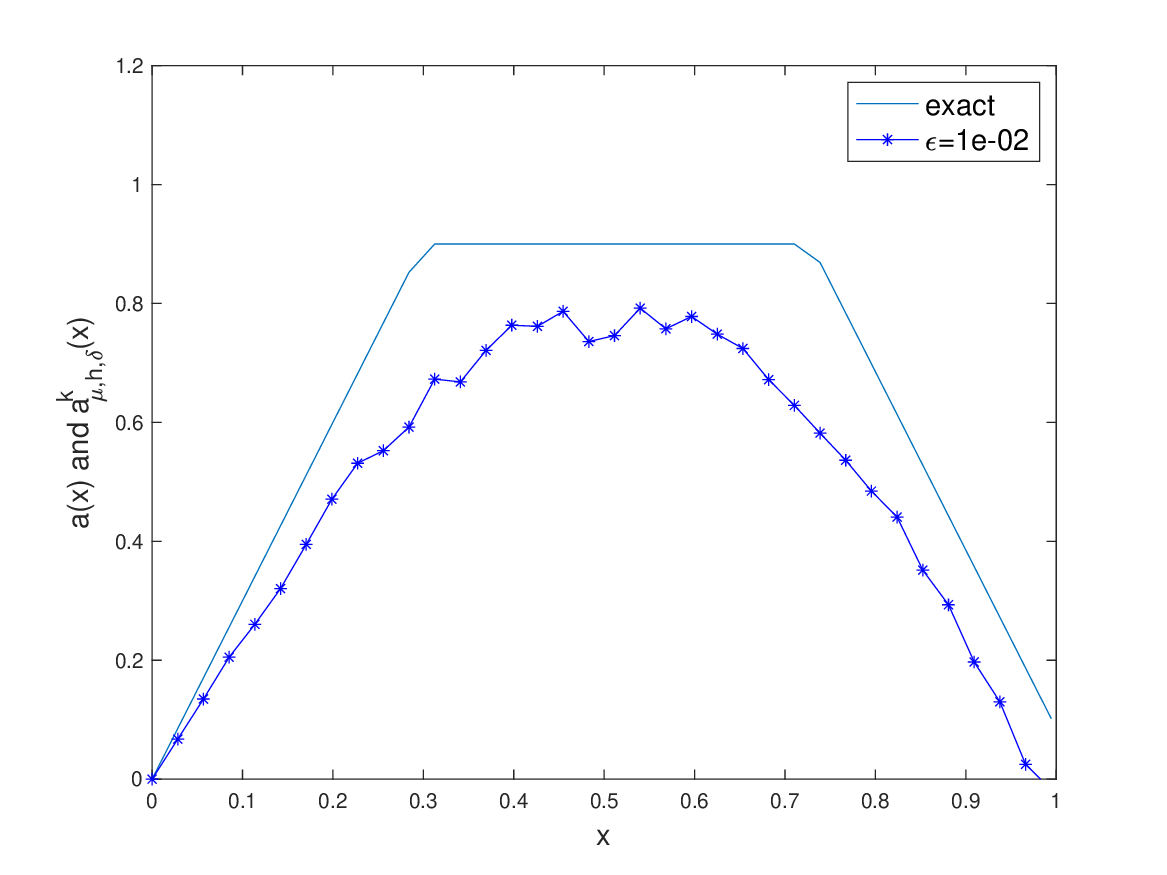}}
     \subfigure{
    \includegraphics[width=0.35\textwidth,height=1.5in]{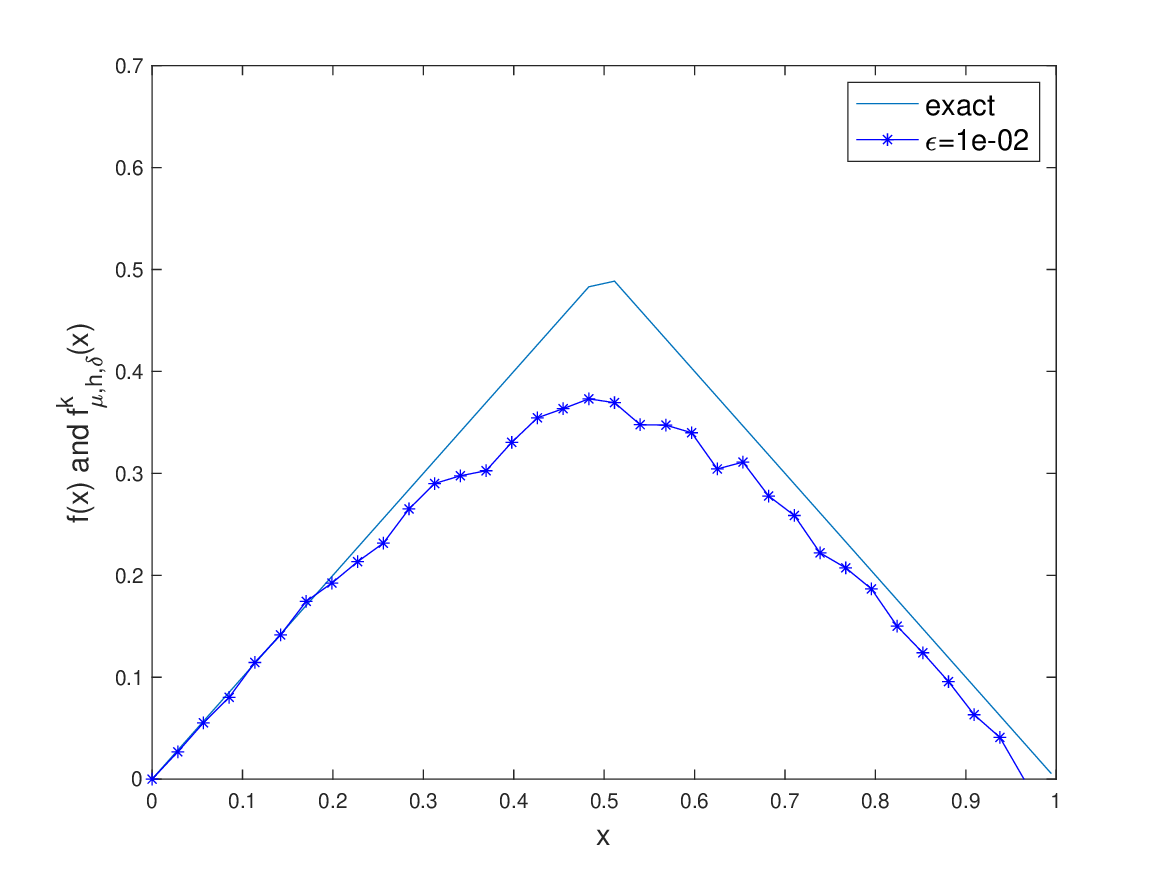}}
  \caption{Approximate solutions $a_{\mu,h,\delta}^{(k)}$ (left column) and $f_{\mu,h,\delta}^{(k)}$(right column) of Example 2 based on lumped mass method with various $\epsilon$.}\label{fig3}
\end{figure}

\begin{table}[h]
\centering
\begin{tabular}{ccccc}
\hline
$\epsilon$ & $re_a$ & Order(a) & $re_f$ & Order(f) \\
%\hline
%$0.00095$ & 0.0843 & $\delta=9.5279e-05$ & 0.0636 & $\delta=9.5279e-05$ \\
\hline
$0.001$ & 0.0955 & 0.3752 & 0.0753 & 0.3715 \\
\hline
$0.003$ & 0.1386 & 0.3527 & 0.1065 & 0.3812 \\
\hline
$0.006$ & 0.2019 & 0.2909 & 0.1519 & 0.3298 \\
\hline
$0.009$ & 0.2234 & 0.2847 & 0.1683 & 0.3257 \\
\hline
$0.01$ & 0.2394 & 0.2764 & 0.1753 & 0.3251 \\
\hline
\end{tabular}
\caption{The relative errors and convergence orders based on standard Galerkin method of Example 2.}\label{tab2}
\end{table}

The second objective of this example is to compare the performance of proposed method with two different semidiscrete schemes. Then we apply formulae (\ref{saklum}) and (\ref{sfklum}) to compute the approximate solutions of the exact solutions $(a,f)$. It can be seen from Figure \ref{fig3} and Table \ref{tab3} that the reconstructed solutions exhibit well agreement with the exact solutions overall. By comparison of results presented in Tables \ref{tab2} and \ref{tab3}, it is evident that the relative errors and convergence orders are essentially identical. This observation confirms that both the standard Galerkin method and lumped mass method yield equivalent convergence rates.

\begin{table}[h]
\centering
\begin{tabular}{ccccc}
\hline
$\epsilon$ & $re_a$ & Order(a) & $re_f$ & Order(f) \\
%\hline
%$0.00095$ & 0.0841 & $\delta=1.1984e-04$ & 0.0637 & $\delta=1.1984e-04$ \\
\hline
$0.001$ & 0.0953 & 0.3752 & 0.0754 & 0.3715 \\
\hline
$0.003$ & 0.01381 & 0.3534 & 0.1071 & 0.3798 \\
\hline
$0.006$ & 0.2010 & 0.2917 & 0.1530 & 0.3283 \\
\hline
$0.009$ & 0.2219 & 0.2861 & 0.1701 & 0.3234  \\
\hline
$0.01$ & 0.2378 & 0.2778 & 0.1772 & 0.3227  \\
\hline
\end{tabular}
\caption{The relative errors and convergence orders based on lumped mass method of Example 2.}\label{tab3}
\end{table}

\subsection{Two-dimensional case}\label{sec5.2}

We take $\Omega=(0,1)^2$ as a unit square. To discretize the problem, we partition the interval $(0,1)$ into $K$ subintervals of equal width, where the mesh size is given by $h=\frac{1}{K}$. Consequently, the domain $\Omega$ is subdivided into $K^2$ smaller squares. By connecting the diagonals of each of these small squares, we obtain a symmetric triangulation of $\Omega$. Then it follows from \cite{Jin+Lazarov+Pasciak+Zhou-2015} that the eigensystem of lumped mass method-based discrete Laplacian $-\bar{\Delta}_h$ can be given by 
\begin{equation*}
\bar{\lambda}_{n,m}^h=\frac{4}{h^2}\left(\sin^2(\frac{n\pi h}{2})+\sin^2(\frac{m\pi h}{2})\right),\quad\bar{\varphi}_{n,m}^h=2\sin(n\pi x_i)\sin(m\pi y_j),\quad 1\leq n,m\leq K-1,
\end{equation*}
where $\left(x_i,y_j\right)(i,j=1\cdots,K-1)$ are mesh points. 

\textbf{Example 3.} Let $a(x,y)=\frac12e^{-(x+y)}\sin(\pi x)\sin(\pi y)$, $f(x,y)=\frac12xy(1-x)(1-y)$, $\epsilon=0.001$ and $\mu=\frac16\delta^{\frac12}$.

\begin{figure}[htbp]
  \centering
  \subfigure{
    \includegraphics[width=0.35\textwidth,height=1.5in]{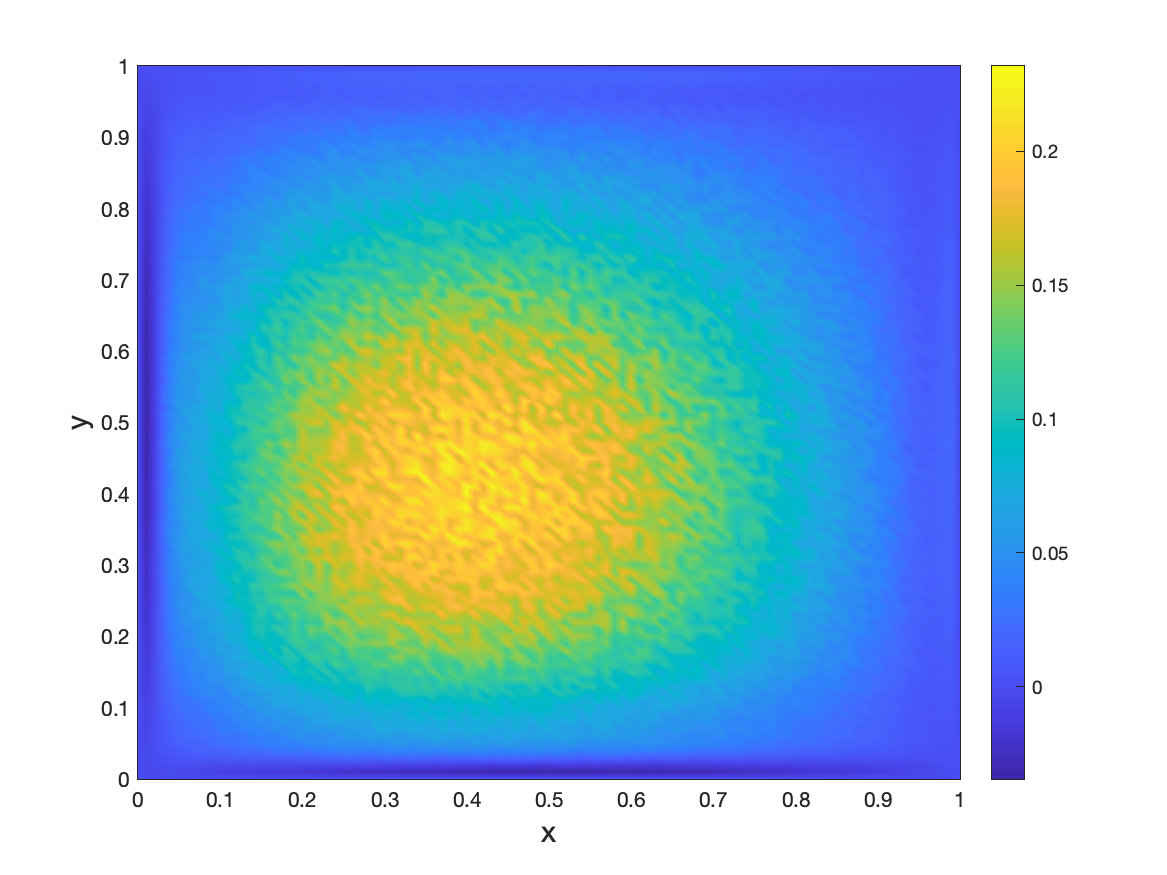}}
  \subfigure{
    \includegraphics[width=0.35\textwidth,height=1.5in]{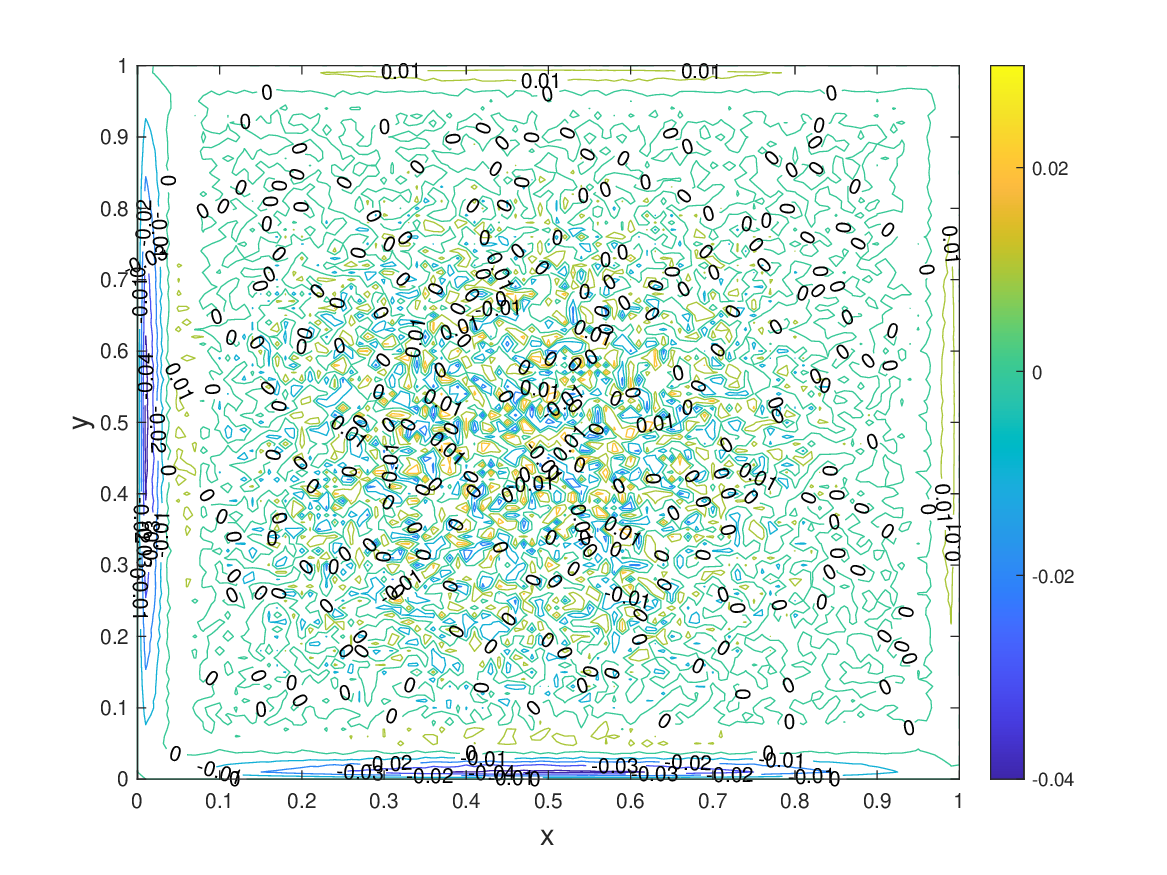}}\\
      \subfigure{
    \includegraphics[width=0.35\textwidth,height=1.5in]{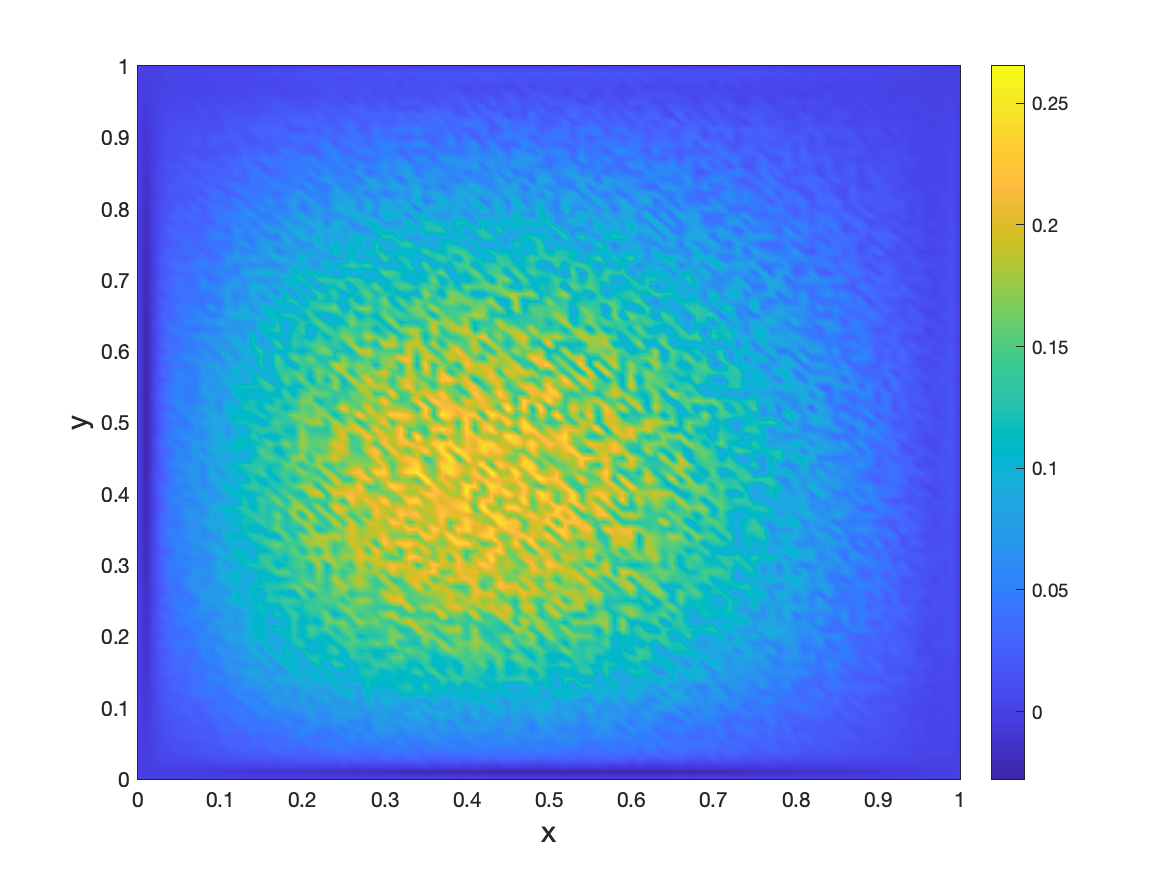}}
      \subfigure{
    \includegraphics[width=0.35\textwidth,height=1.5in]{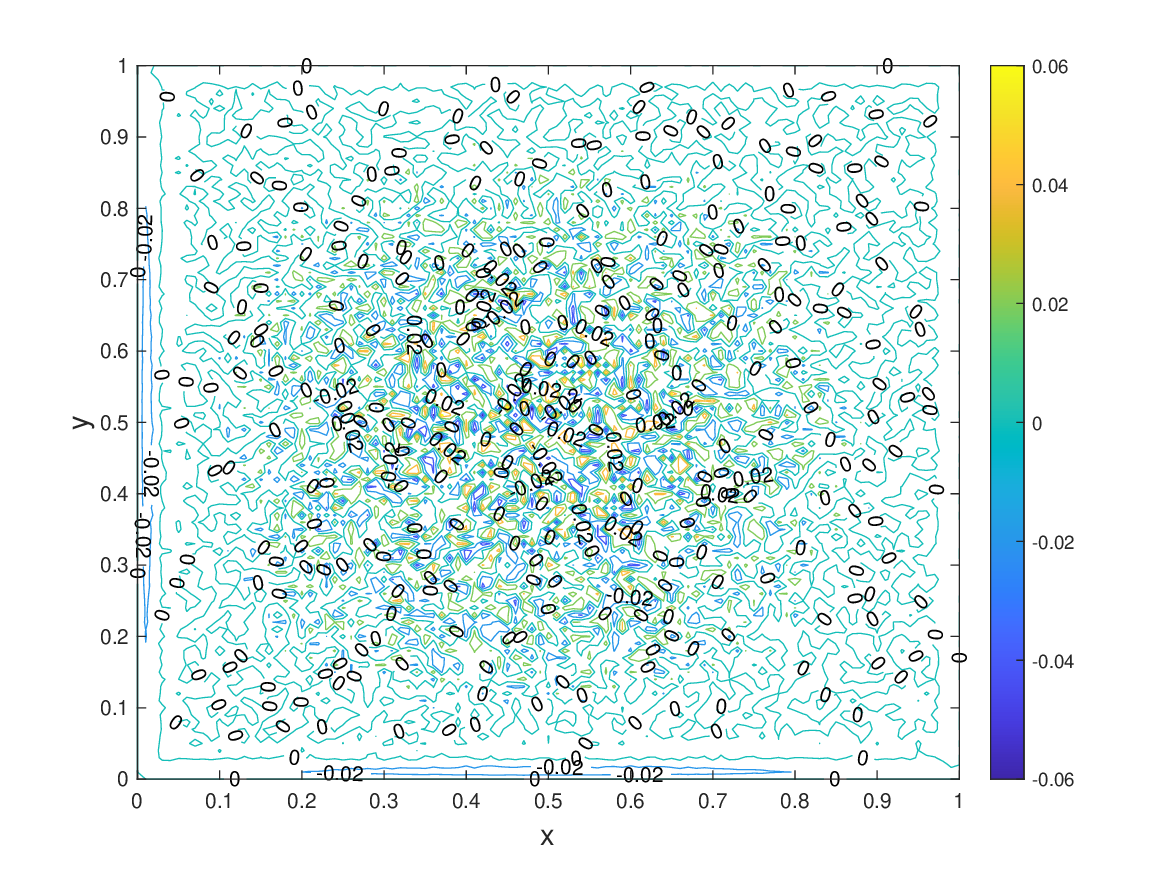}}\\
      \subfigure{
    \includegraphics[width=0.35\textwidth,height=1.5in]{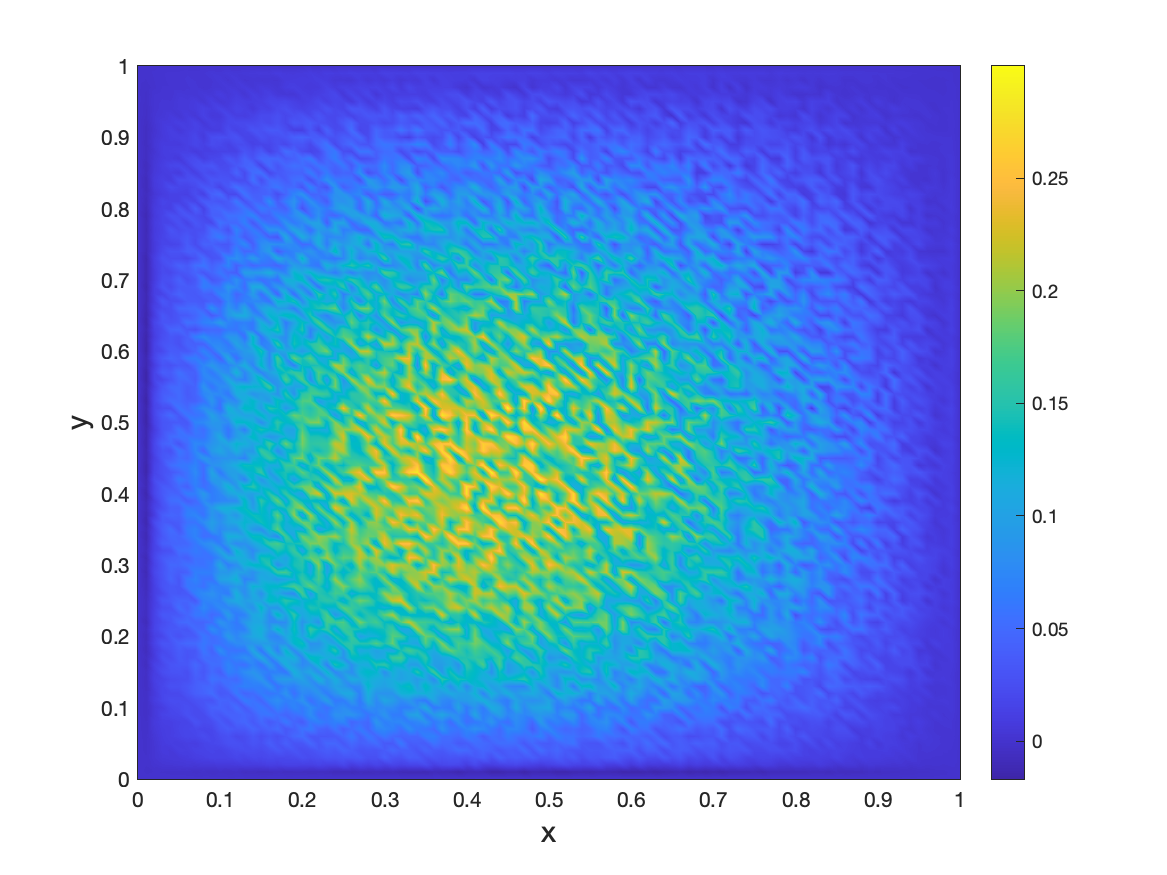}}
      \subfigure{
    \includegraphics[width=0.35\textwidth,height=1.5in]{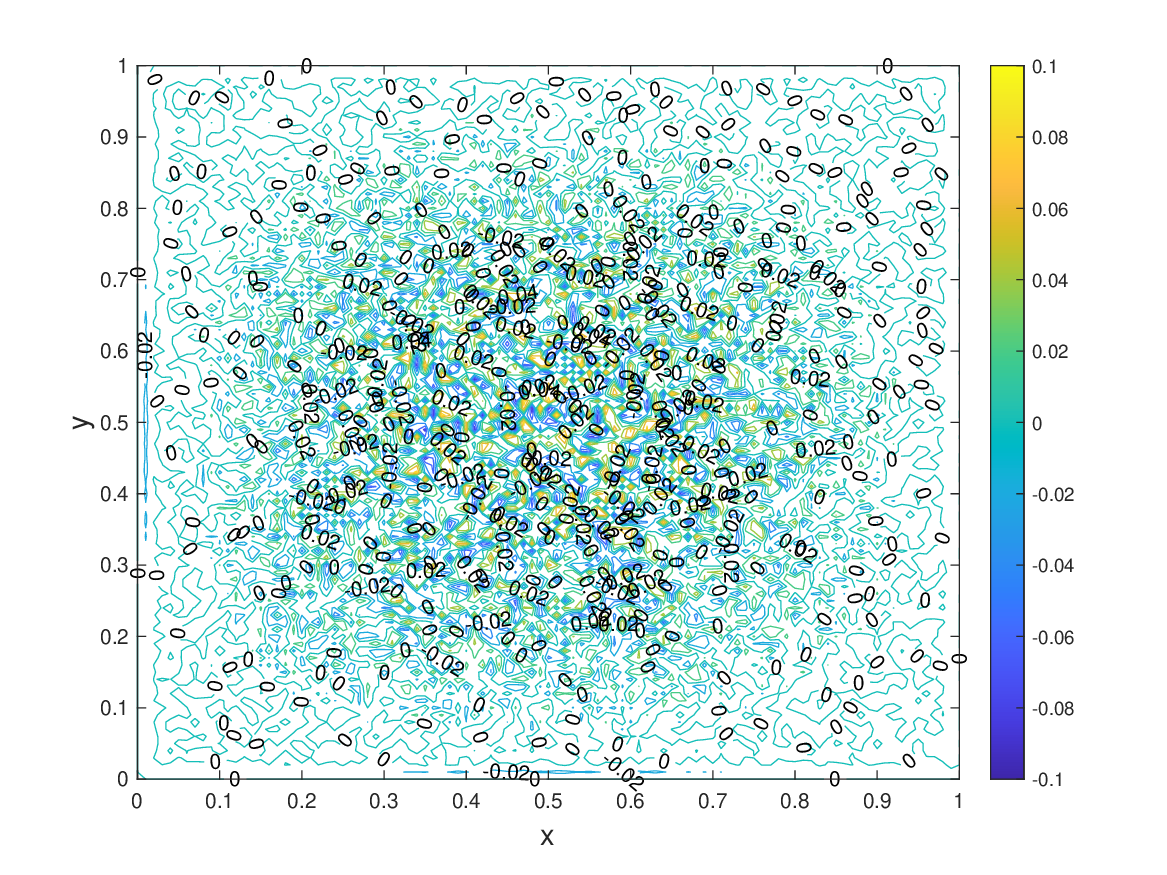}}\\
  \subfigure{
    \includegraphics[width=0.35\textwidth,height=1.5in]{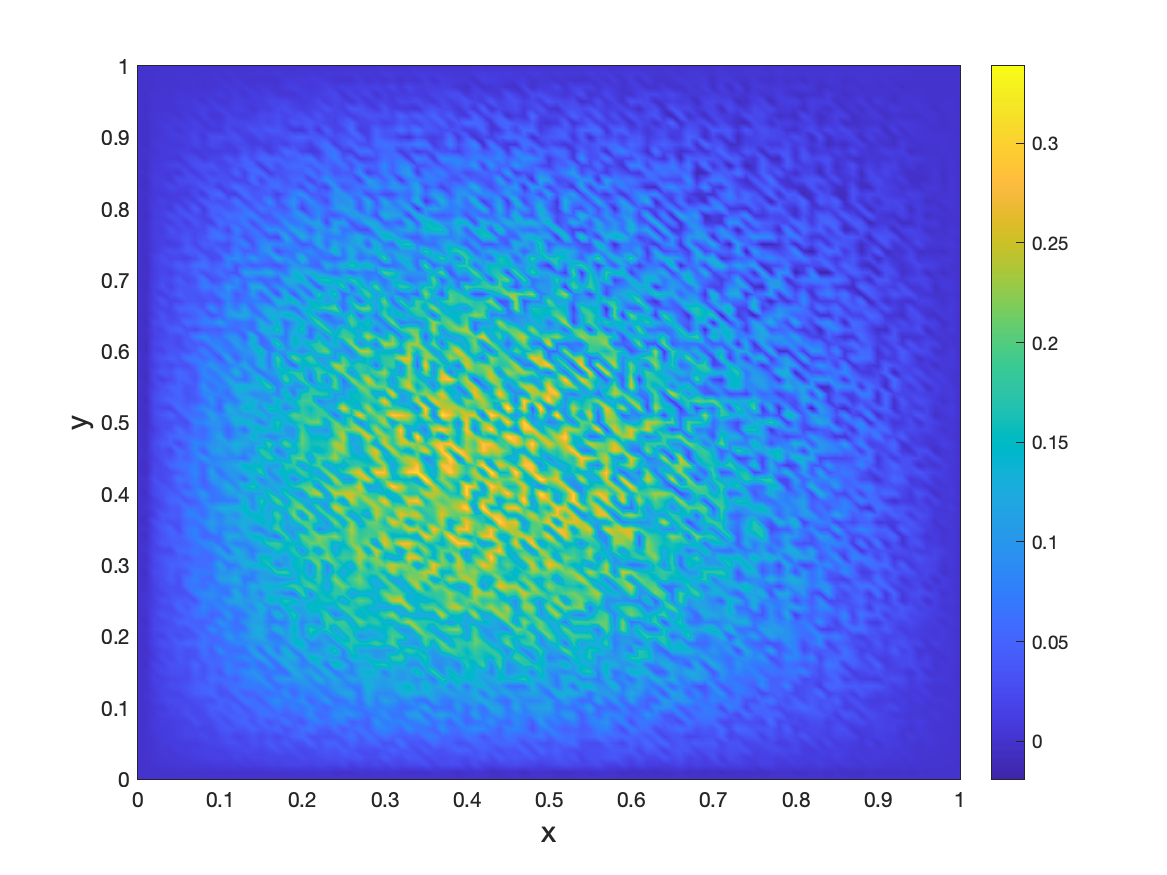}}
      \subfigure{
    \includegraphics[width=0.35\textwidth,height=1.5in]{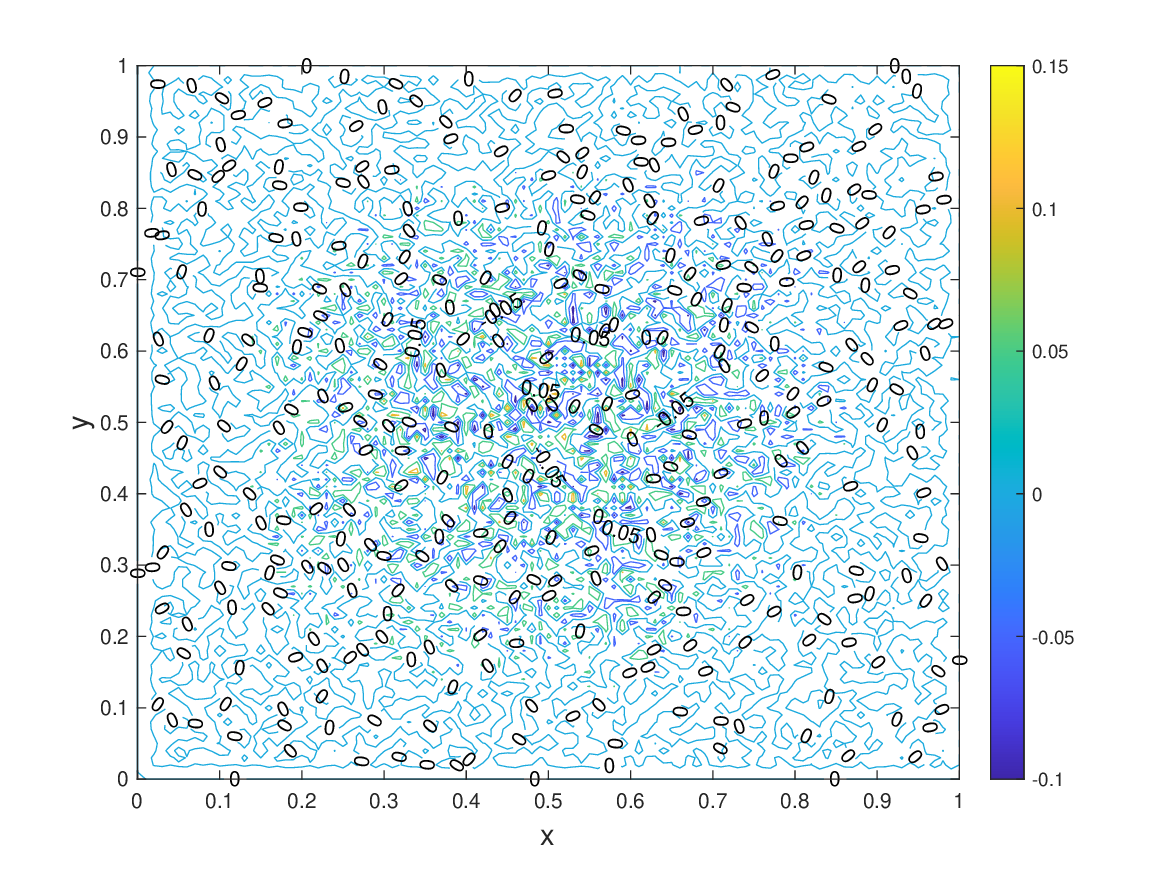}}\\
      \subfigure{
    \includegraphics[width=0.35\textwidth,height=1.5in]{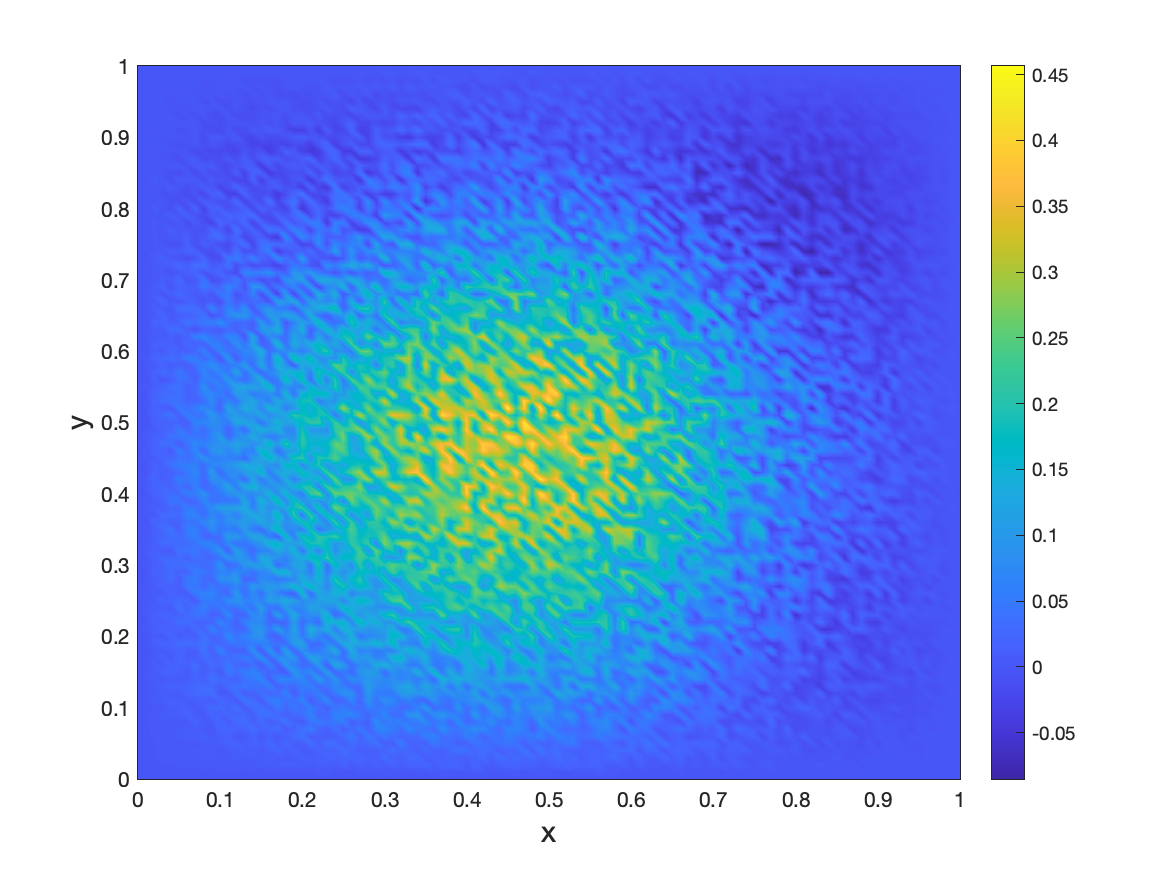}}
      \subfigure{
    \includegraphics[width=0.35\textwidth,height=1.5in]{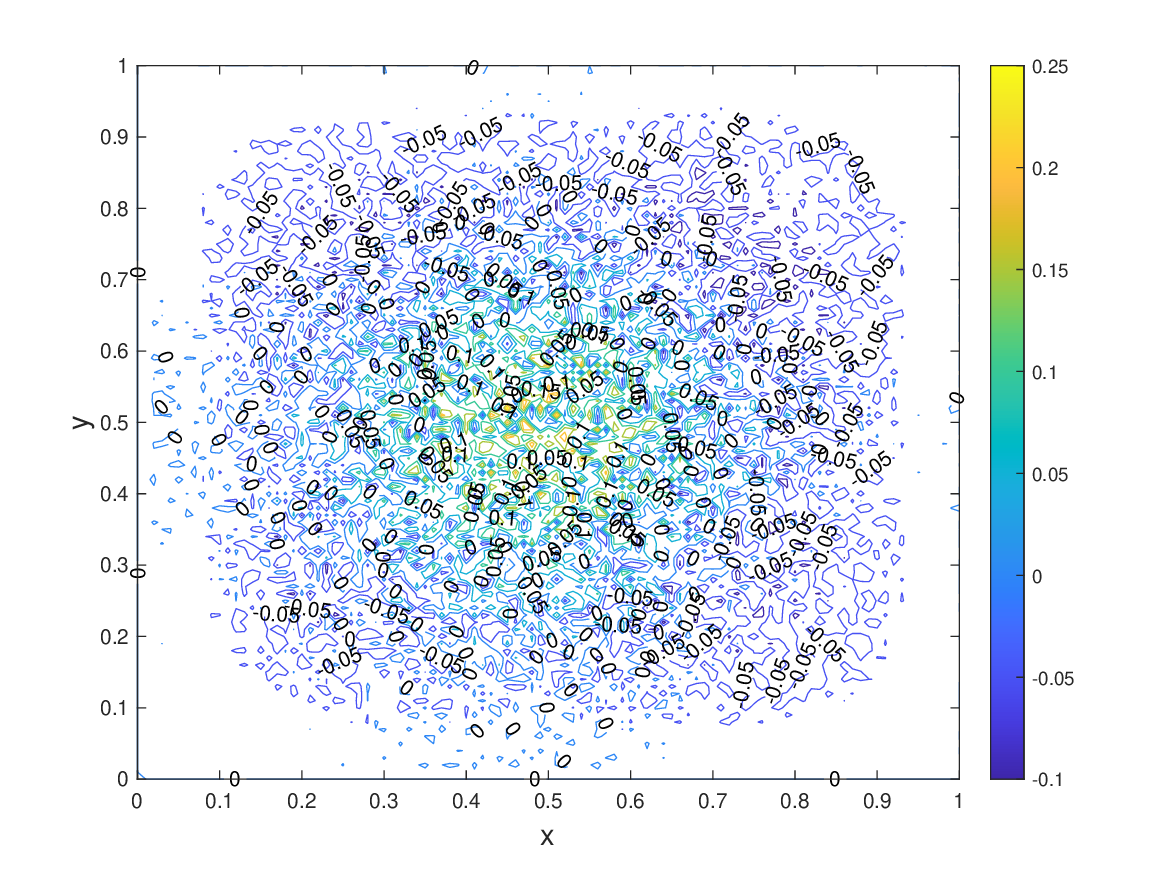}}
  \caption{Approximate solutions $a_{\mu,h,\delta}^{(k)}$ (left column) and the absolute errors (right column) of Example 3 with $\alpha=1.1,1.2,1.3,1.4,1.5$ (from top to bottom).}\label{fig4}
\end{figure}

\begin{figure}[htbp]
  \centering
      \subfigure{
    \includegraphics[width=0.35\textwidth,height=1.5in]{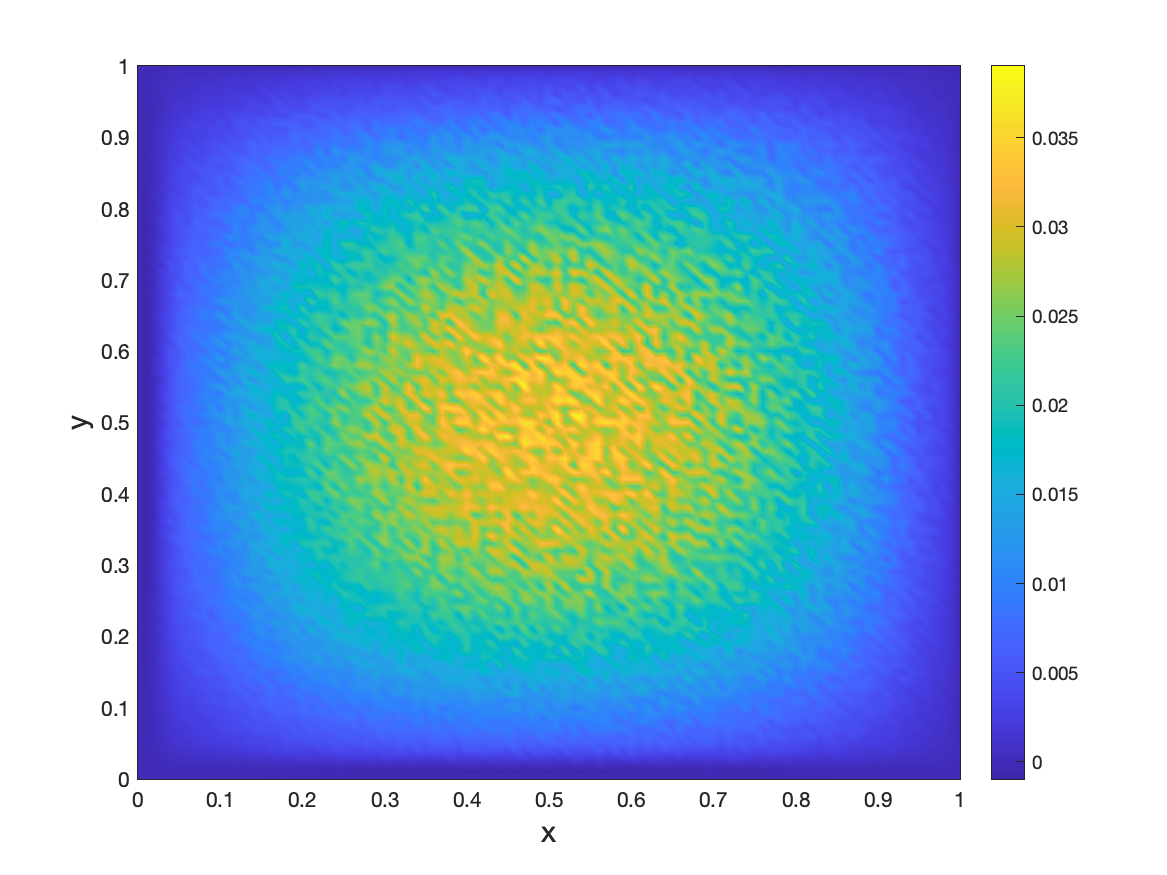}}
      \subfigure{
    \includegraphics[width=0.35\textwidth,height=1.5in]{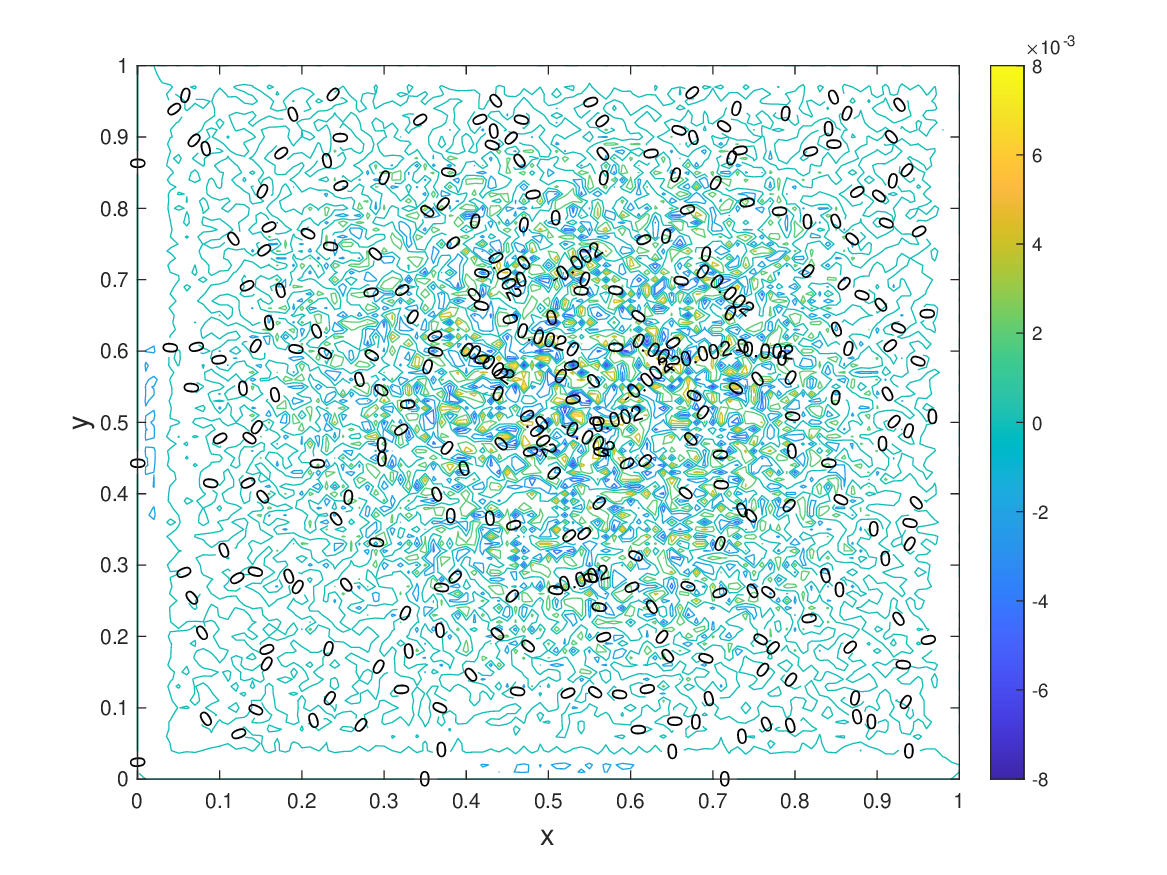}}\\
      \subfigure{
    \includegraphics[width=0.35\textwidth,height=1.5in]{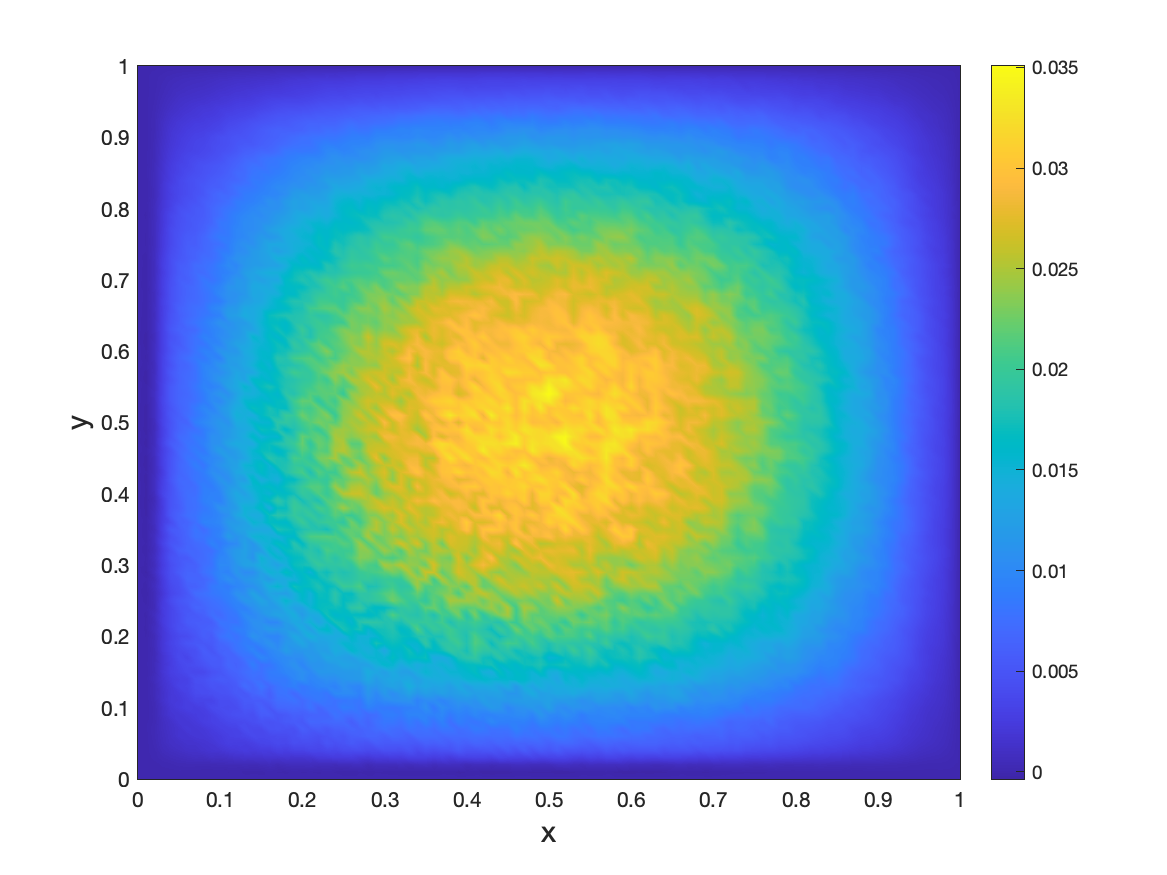}}
      \subfigure{
    \includegraphics[width=0.35\textwidth,height=1.5in]{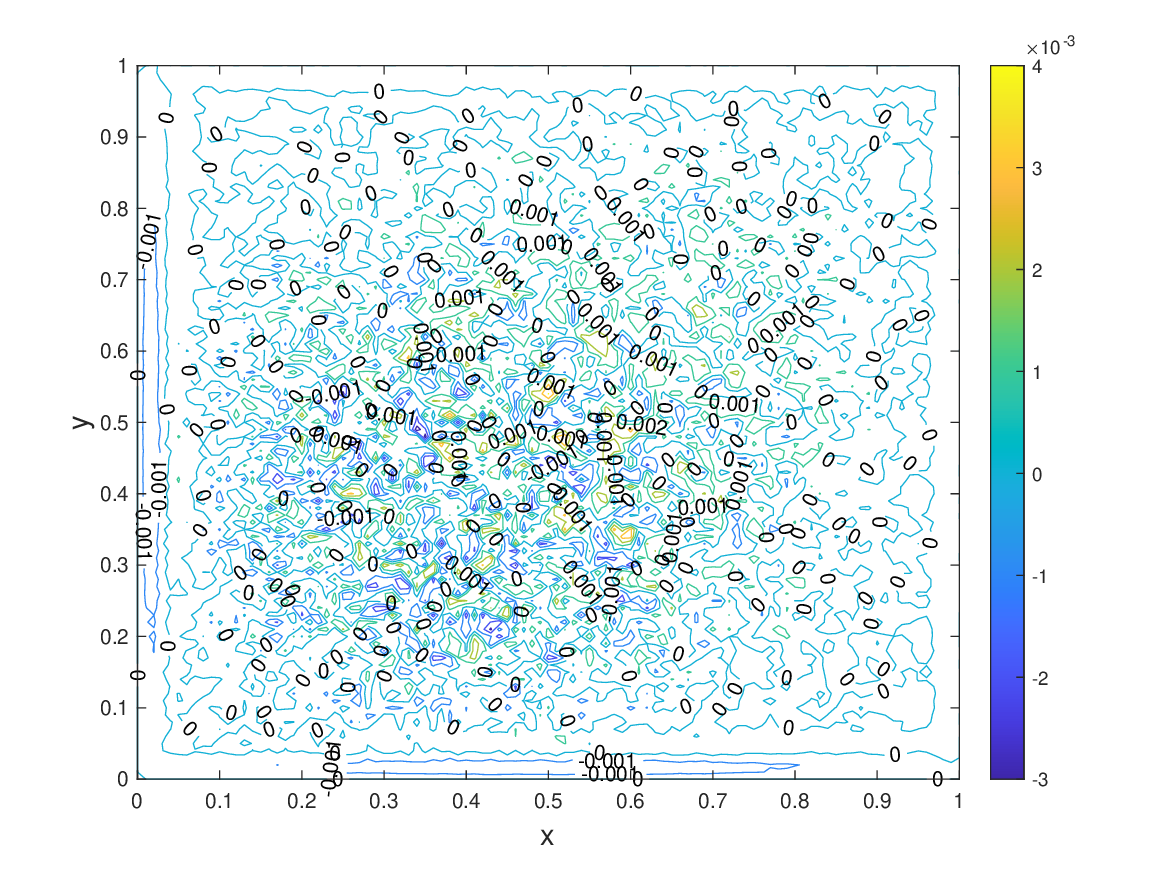}}
      \subfigure{
    \includegraphics[width=0.35\textwidth,height=1.5in]{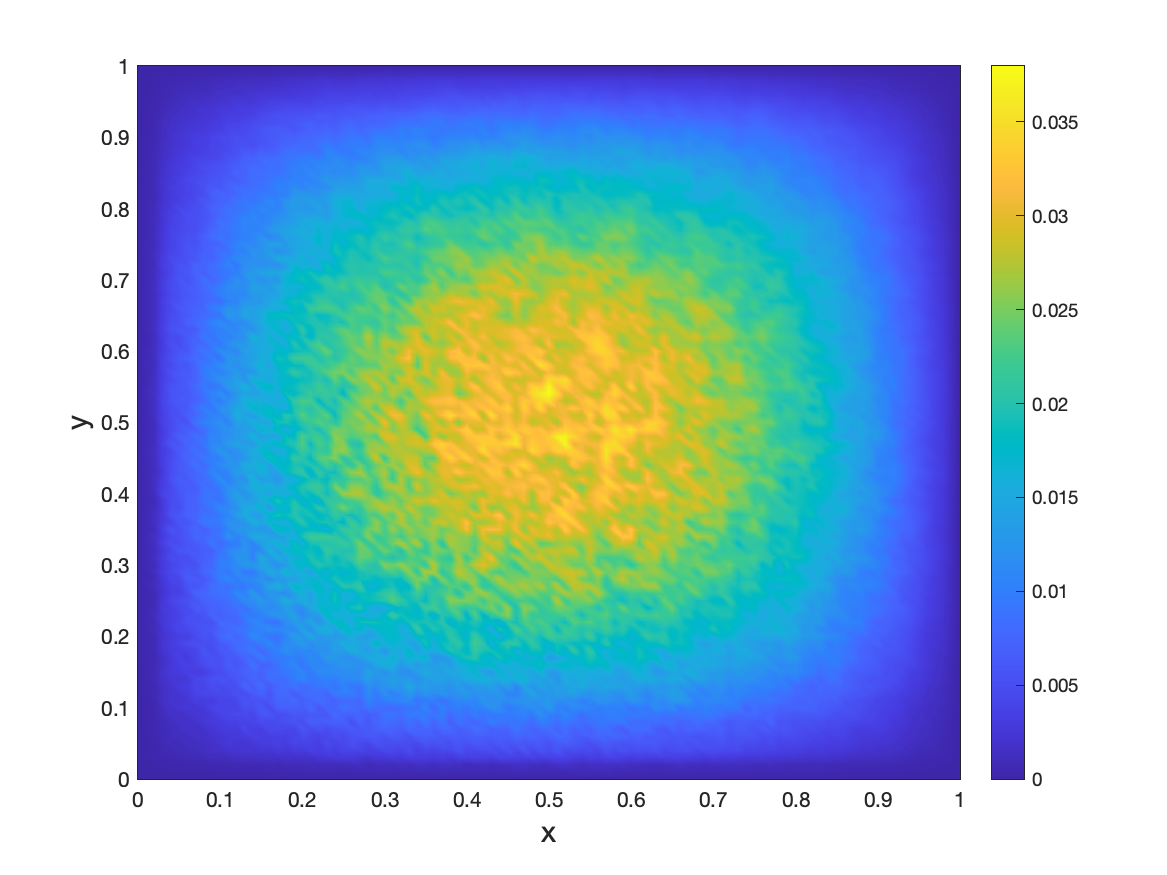}}
      \subfigure{
    \includegraphics[width=0.35\textwidth,height=1.5in]{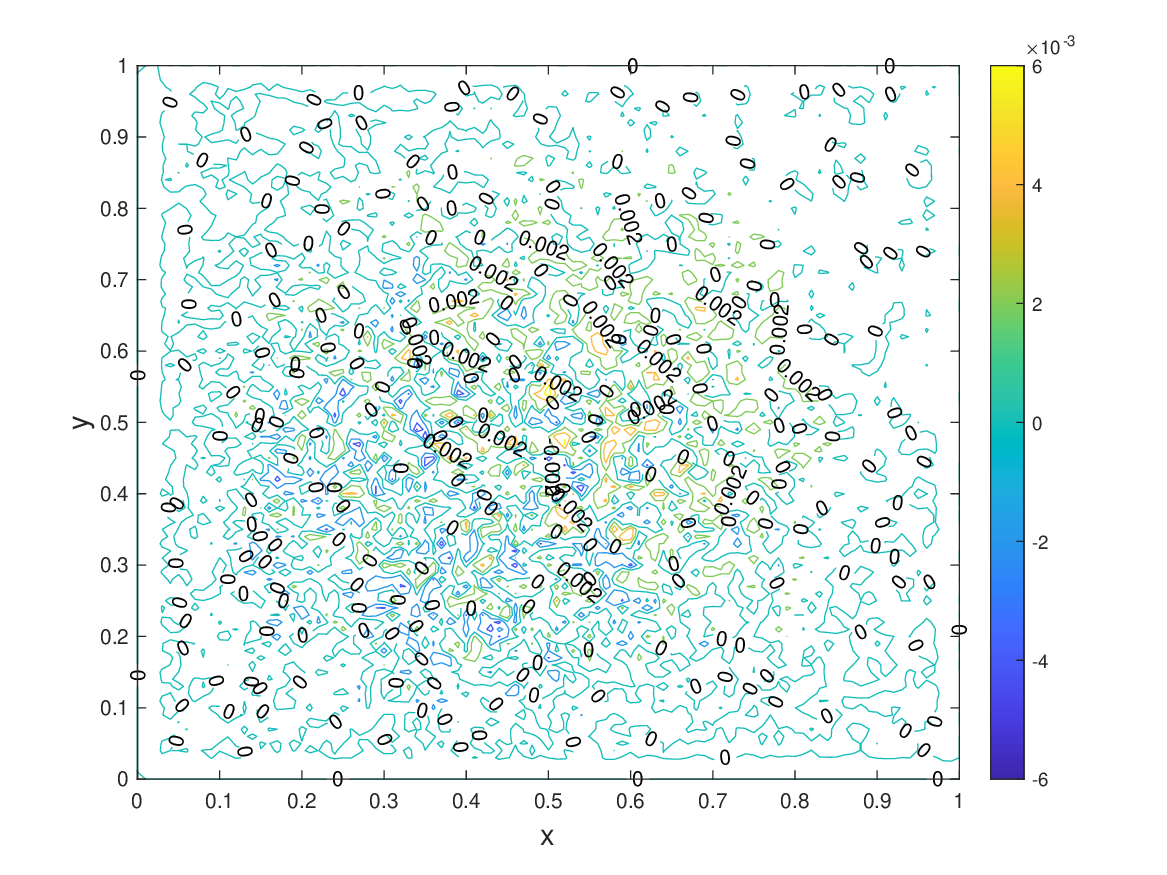}}\\
      \subfigure{
    \includegraphics[width=0.35\textwidth,height=1.5in]{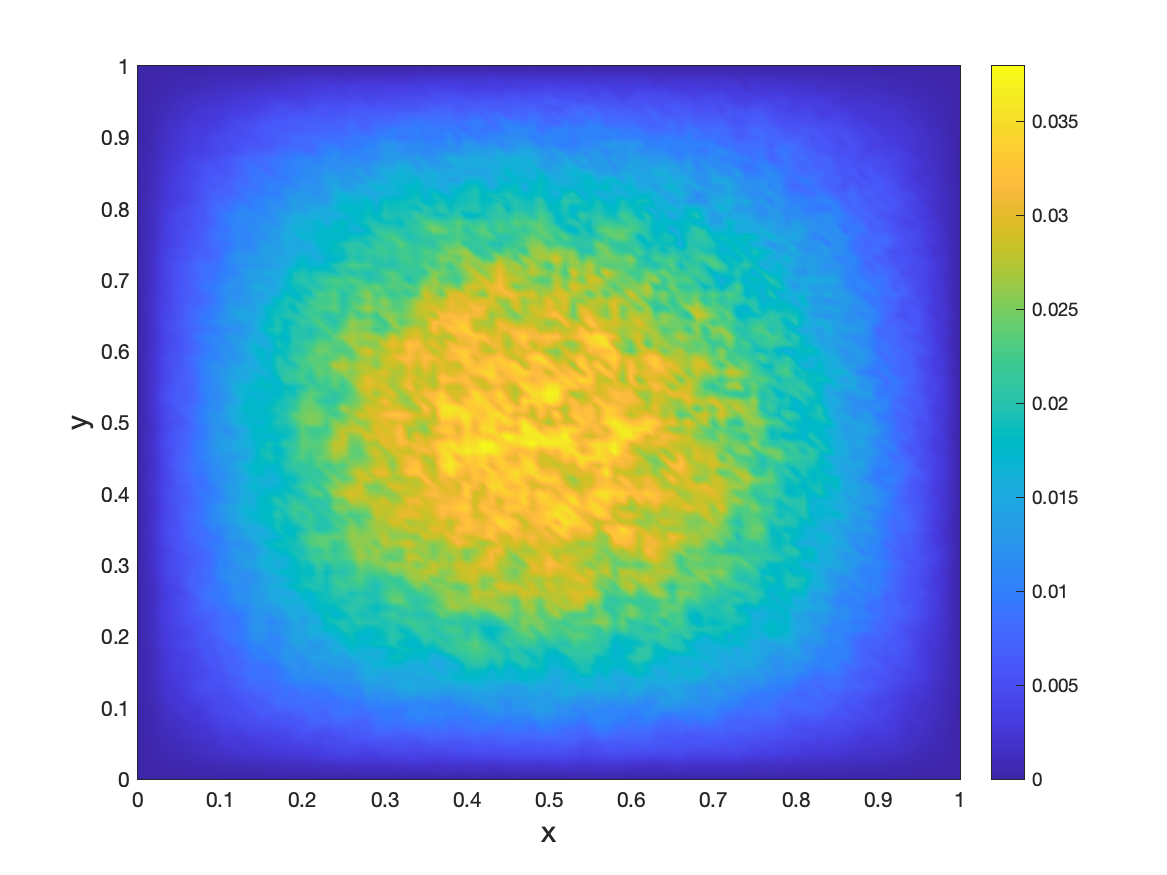}}
      \subfigure{
    \includegraphics[width=0.35\textwidth,height=1.5in]{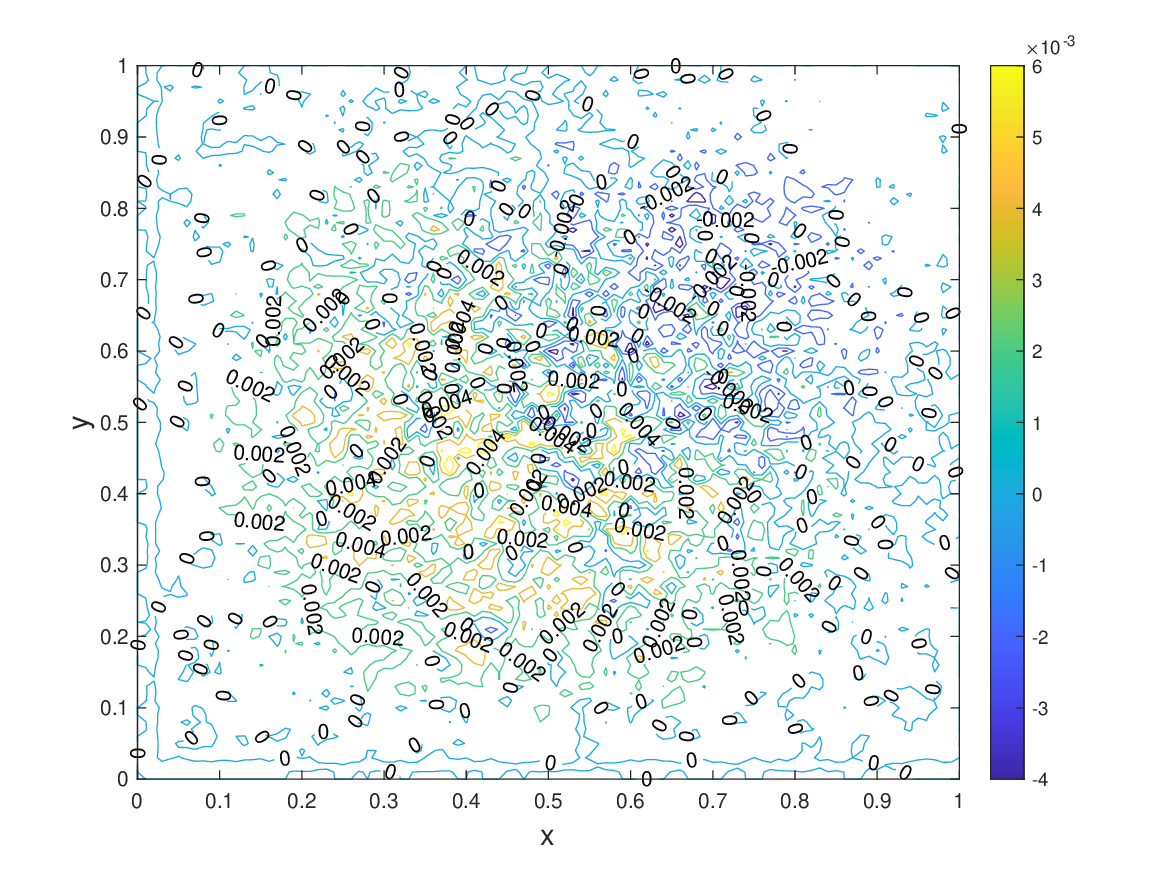}}\\
      \subfigure{
    \includegraphics[width=0.35\textwidth,height=1.5in]{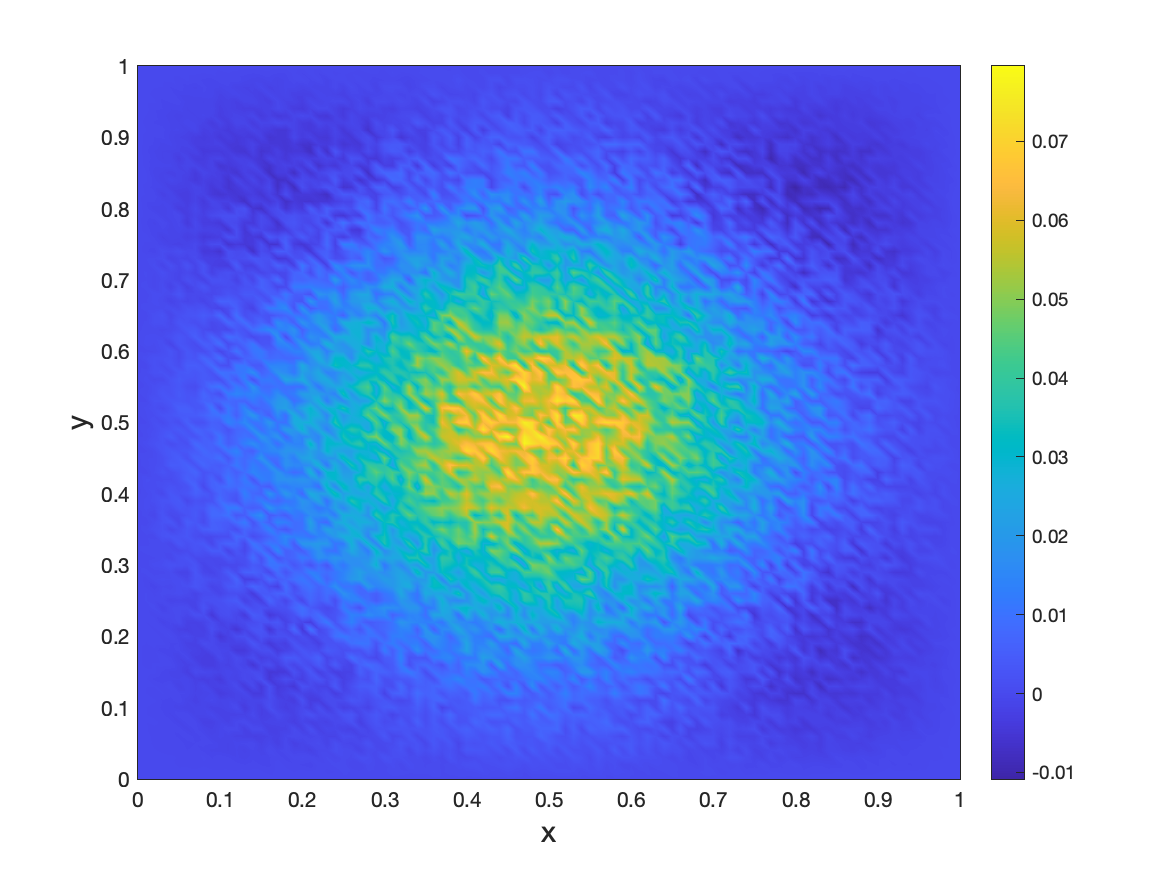}}
      \subfigure{
    \includegraphics[width=0.35\textwidth,height=1.5in]{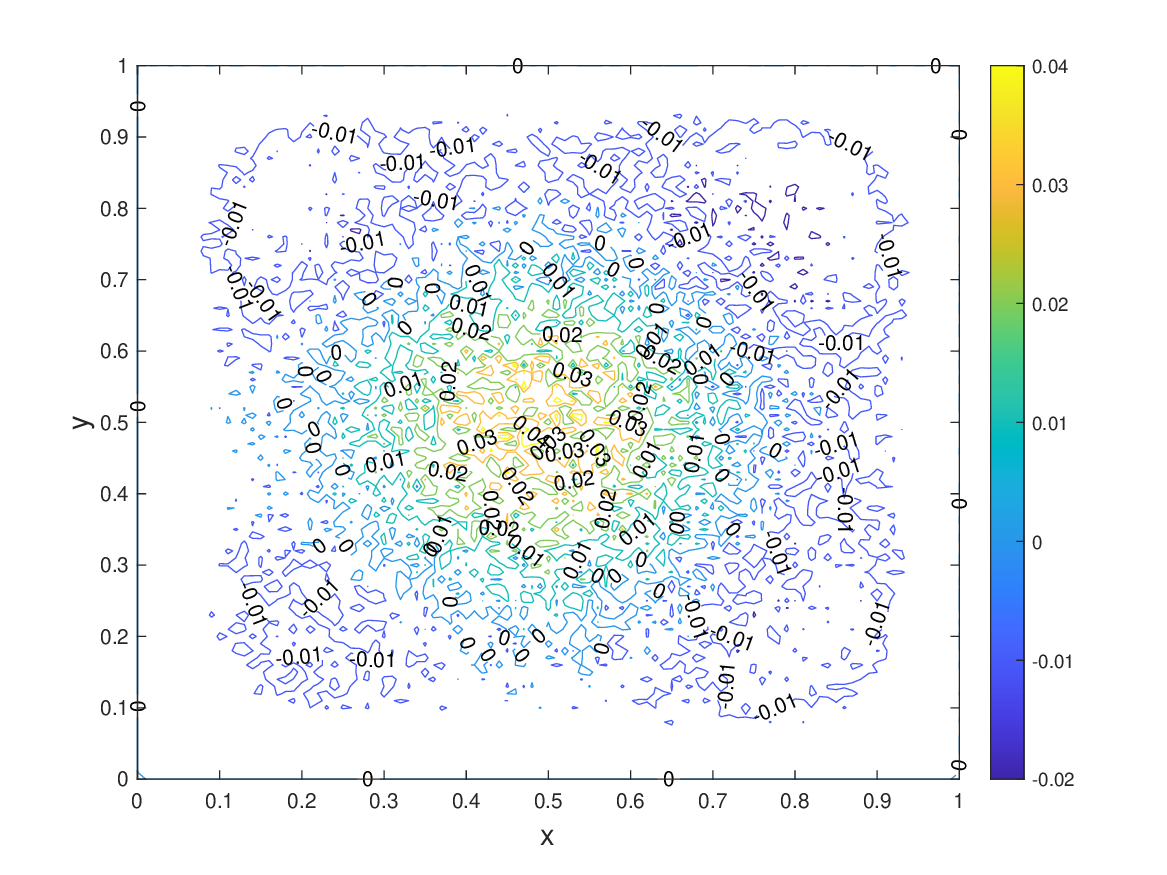}}
  \caption{Approximate solutions $f_{\mu,h,\delta}^{(k)}$ (left column) and the absolute errors (right column) of Example 3 with $\alpha=1.1,1.2,1.3,1.4,1.5$ (from top to bottom)}\label{fig5}
\end{figure}

To obtain the additional measurements, we apply the alternating direction finite difference method proposed in \cite{Zhang+Sun-2011} to solve problem (\ref{IP}). The purpose of this example is to examine the behavior of the numerical results under different fractional orders or the values of final time. Figures \ref{fig4} and \ref{fig5} report the approximate solutions and absolute errors of initial value and source term, respectively. It can be learned from Figure \ref{fig4} and Table \ref{tab4} that the reconstructed solutions fit the exact solutions more accurately with smaller fractional order $\alpha$. And it begins to deteriorate significantly when $\alpha$ is larger than $1.5$. The main reason is that the Mittag-Leffler functions exhibit a pronounced oscillatory behavior with larger $\alpha$. Next, we intend to explore the impact of final time $T_2$ on the numerical results. We fix $\alpha=1.5$ and take $T_2=1,2,5,8,10$, respectively.  It is evident that the relative errors exhibit a decreasing trend followed by stabilization as 
$T_2$ increases. This phenomenon validates the conclusion in Lemma \ref{SQQS_bound} that the final time should be taken large sufficiently.

\begin{table}[h]
\centering
\begin{tabular}{cccccc}
\hline
$\alpha$ & $1.1$ & $1.2$ & $1.3$ & $1.4$ & $1.5$ \\
\hline
$re_a$ & $0.0455$ & $0.0505$ & $0.0758$ & $0.1030$ & $0.2755$ \\
\hline
$re_f$ & $0.0396$ & $0.0195$ & $0.0307$ & $0.0642$ & $0.4696$ \\
\hline
\end{tabular}
\caption{The relative errors of approximate solutions of Example 3 with various fractional order $\alpha$ when $T_2=1$.}\label{tab4}
\end{table}

\begin{table}[h]
\centering
\begin{tabular}{cccccc}
\hline
$T_2$ & $1$ & $2$ & $5$ & $8$ & $10$ \\
\hline
$re_a$ & $0.2755$ & $0.7581$ & $0.0637$ & $0.0327$ & $0.0565$ \\
\hline
$re_f$ & $0.4696$ & $0.3533$ & $0.0944$ & $0.0950$ & $0.0947$ \\
\hline
\end{tabular}
\caption{The relative errors of approximate solutions of Example 3 with various final time $T_2$ when $\alpha=1.5$.}\label{tab5}
\end{table}

\section{Conclusion}\label{sec6}
This paper considers an inverse problem of recovering the initial value and source term simultaneously in a time-fractional diffusion-wave equation. There are two major difficulties when we deal with the inverse problem. On one hand, both the IVVP and ISP are ill-posed, some regularization methods need to be adopted. On the other hand, the initial value and source term are coupled, which increases the difficulty of the inverse problem. To efficiently solve the problem, we come up with an alternating approach for decoupling the inverse problem. Then we employ the quasi-boundary value regularization method to solve each decoupled inverse problem. Moreover, for the numerical implementation of the proposed method, we utilize standard Galerkin method and lumped mass method for discretization the problem in space. Comprehensive error estimates considering the noise level, regularization parameter and discretization parameter are given. Numerical experiments and error analysis are reported to support the theoretical results.

While this work addresses several aspects, there remain many open problems. Notably, the paper does not consider time discretization schemes. Leveraging existing numerical methods for computing Mittag-Leffler functions, such as method  outlined in\cite{Garrappa-2015}, we can directly obtain the numerical solutions of the inverse problem. However, for more intricate models like fractional diffusion equations with multi-term fractional derivatives, distributed order and variable order,  some time discretization schemes need to be applied. There is existing literature on the inverse problems of the mentioned models. However, error estimates for fully discrete schemes related to these inverse problems, especially nonlinear inverse problems, are scarce. These issues will be addressed in our ongoing and future work.

\bigskip
\noindent{\bf Acknowledgements} 
This research is supported by the Natural Science Basic Research Program of Shaanxi (Nos. 2024JC-YBQN-0050, 2024JC-YBMS-034, 2023-JC-YB-054) and the Fundamental Research Funds for the Central Universities (No. ZYTS24077).\\

\end{document}